\documentclass{article}

\usepackage[colorlinks,linkcolor=magenta,citecolor=blue, pagebackref=true]{hyperref}
\usepackage[margin=1.2in]{geometry}

\usepackage[T1]{fontenc}    
\usepackage{url}            
\usepackage{booktabs}       
\usepackage{nicefrac}       
\usepackage{microtype}      
\usepackage{tcolorbox}
\usepackage[numbers]{natbib}
\usepackage[T1]{fontenc}    
\usepackage{url}            
\usepackage{booktabs}       
\usepackage{nicefrac}       
\usepackage{microtype}      
\usepackage{algpseudocode}
\usepackage{placeins}
\usepackage{titlesec}
\usepackage{authblk}
\usepackage{graphicx}

\usepackage{setspace}

\usepackage[shortlabels]{enumitem}
\usepackage{amsthm,amsmath,amsfonts,mathtools,algorithmicx,dsfont,amssymb,mathabx}
\usepackage[capitalize]{cleveref}
\usepackage{mathrsfs}
\usepackage{array}
\usepackage{cancel}
\usepackage{tikz}
\usetikzlibrary{positioning}
\usetikzlibrary{calc} 

\usepackage{autonum}
\theoremstyle{plain}

\newtheorem{condition}{Condition}
\newtheorem{claim}{Claim}[section]
\newtheorem{theorem}{Theorem}[section]
\newtheorem{proposition}[theorem]{Proposition}
\newtheorem{corollary}[theorem]{Corollary}
\newtheorem{lemma}[theorem]{Lemma}

\theoremstyle{remark}
\newtheorem{definition}[theorem]{Definition}

\newtheorem{remark}{Remark}

\DeclareMathOperator{\essinf}{\mathrm{essinf}}

\DeclareMathOperator{\Var}{\mathrm{var}}

\DeclareMathOperator{\Cov}{\mathrm{cov}}

\DeclareMathOperator{\Ind}{\mathds 1}

\newcommand{\eqas}{=}

\newcommand{\samplesplit}{\mathrm{split}}
\newcommand{\splitSet}{{\mathcal D}}
\newcommand{\splitSetInfty}{\splitSet_\infty}
\newcommand{\whichSplitSet}{\mathcal S}

\DeclareMathOperator*{\argmin}{\arg\!\min}

\newcommand{\train}{\mathrm{trn}}
\newcommand{\eval}{\mathrm{eval}}

\newcommand{\SA}{\Gamma^\mathrm{SA}}
\newcommand{\Bias}{\Gamma^\mathrm{B}}
\newcommand{\EP}{\Gamma^\mathrm{EP}}
\newcommand{\SAtilde}{\widetilde S^\mathrm{SA}}
\newcommand{\Biastilde}{\widetilde S^\mathrm{B}}
\newcommand{\EPtilde}{\widetilde S^\mathrm{EP}}
\newcommand{\simiid}{\sim}
\newcommand{\iid}{{\text{i.i.d.}}}
\newcommand{\Teval}{{T}}
\newcommand{\Ttrain}{{T'}}
\newcommand{\as}[0]{\mathrm{a.s.}}
\newcommand{\oas}[1]{{o_\mathrm{a.s.}( #1 )}}%
\newcommand{\Oas}[1]{{O_\mathrm{a.s.} ( #1  )}}

\newcommand{\asympcs}{AsympCS}

\newcommand{\boundary}{{\mathfrak B}}
\newcommand{\LP}{{L_2(\PP)}}
\newcommand{\eps}{\varepsilon}

\newcommand{\PP}{\mathbb P}
\newcommand{\NN}{\mathbb N}
\newcommand{\QQ}{\mathbb Q}
\newcommand{\EE}{\mathbb E}
\newcommand{\RR}{\mathbb R}

\newcommand{\Fcal}{\mathcal F}
\newcommand{\Gcal}{\mathcal G}

\newcommand{\Lcal}{\mathcal L}

\newcommand{\Ncal}{\mathcal N}

\newcommand{\Tcal}{\mathcal T}

\newcommand{\seq}[4]{(#1_{#2})_{{#2}={#3}}^{#4}}
\newcommand{\infseq}[3]{\seq{#1}{#2}{#3}{\infty}}
\newcommand{\infseqt}[1]{(#1)_{t=1}^\infty}
\newcommand{\infseqtcal}[1]{(#1)_{t\in \Tcal}}
\newcommand{\infseqtm}[1]{(#1)_{t=m}^\infty}

\newcommand{\cs}{CS}
\newcommand{\ci}{CI}

\newcommand{\1}{\mathds 1}

\newcommand{\frob}{{\mathrm{F}}}

\newcommand{\ball}{{\mathrm{B}_d}}

\newcommand{\fluid}{\mathscr L^\mathrm{24h}}

\newcommand{\opnorm}[1]{{\left\vert\kern-0.25ex\left\vert\kern-0.25ex\left\vert #1 
    \right\vert\kern-0.25ex\right\vert\kern-0.25ex\right\vert}}

\newcommand{\symdiff}{\triangle}

\newcommand{\blind}{1}

\newcommand\sequentialCausalLink{
    \ifnum\blind=1{
        \href{https://github.com/wannabesmith/drconfseq}{\texttt{\textup{drconfseq}}}
    }
    \else{\href{omitted for blinded version}{\texttt{\textup{omitted for blinded version}}}} \fi
}

\newcommand\reproduceDataLink{
  \ifnum\blind=1{
    \unskip\href{https://github.com/WannabeSmith/drconfseq/tree/main/paper_plots/sepsis}{\texttt{\textup{github.com/WannabeSmith/drconfseq/tree/main/paper\_plots/sepsis}}}
  }
  \else{\href{omitted for blinded version}{\texttt{\textup{omitted for blinded version}}}} \fi
}

\newcommand{\dd}{\mathrm{d}}

\newcommand{\DS}{\mathrm{DS}}

\newcommand{\bvone}{\mathrm{Bv1}}
\newcommand{\bvtwo}{\mathrm{Bv2}}

\theoremstyle{plain}

\newtheorem{conditiongeneral}{Condition}

\newtheorem{conditionlind}{Condition}

\newtheorem{conditiontvate}{Condition}

\newtheorem{conditioncustom}{Condition}[condition]

\newenvironment{conditionalt}[1]{
  \renewcommand\theconditioncustom{#1}
  \conditioncustom
}{\endconditioncustom}

\title{Time-uniform central limit theory\\and asymptotic confidence sequences}
\author{%
Ian Waudby-Smith$^1$, David Arbour$^2$, Ritwik Sinha$^2$, Edward H.\ Kennedy$^1$, and Aaditya Ramdas$^1$\vspace{0.05in}\\
  $^1$Carnegie Mellon University\\
  $^2$Adobe Research \vspace{0.05in}\\
  \texttt{ianws@cmu.edu}, \texttt{arbour@adobe.com}, \texttt{risinha@adobe.com},\\ \texttt{edward@stat.cmu.edu}, \texttt{aramdas@cmu.edu}
}
\date{}

\begin{document}
\maketitle
\setcounter{tocdepth}{2}
\makeatletter
\renewcommand\tableofcontents{%
    \@starttoc{toc}%
}
\makeatother
\begin{abstract}
Confidence intervals based on the central limit theorem (CLT) are a cornerstone of classical statistics. Despite being only asymptotically valid, they are ubiquitous because they permit statistical inference under weak assumptions and can often be applied to problems even when nonasymptotic inference is impossible. 
  This paper introduces time-uniform analogues of such asymptotic confidence intervals, adding to the literature on confidence sequences (CS) --- sequences of confidence intervals that are uniformly valid over time --- which provide valid inference at arbitrary stopping times and incur no penalties for ``peeking'' at the data, unlike classical confidence intervals which require the sample size to be fixed in advance.
  Existing CSs in the literature are nonasymptotic, enjoying finite-sample guarantees but not the aforementioned broad applicability of asymptotic confidence intervals.
  This work provides a definition for ``asymptotic CSs'' and a general recipe for deriving them. Asymptotic CSs forgo nonasymptotic validity for CLT-like versatility and (asymptotic) time-uniform guarantees. While the CLT approximates the distribution of a sample average by that of a Gaussian for a fixed sample size, we use strong invariance principles (stemming from the seminal 1960s work of Strassen) to uniformly approximate the entire sample average process by an implicit Gaussian process. 
  As an illustration, we derive asymptotic CSs for the average treatment effect in observational studies (for which nonasymptotic bounds are essentially impossible to derive even in the fixed-time regime) as well as randomized experiments, enabling causal inference in sequential environments.
\end{abstract}


\tableofcontents
\section{Introduction}\label{section:introduction}
The central limit theorem (CLT) is arguably the most widely used result in applied statistical inference, due to its ability to provide large-sample confidence intervals (CI) and $p$-values in a broad range of problems under weak assumptions. Examples include (a) nonparametric estimation of means, such as population proportions, (b) maximum likelihood and other M-estimation problems \citep{van2000asymptotic}, and (c) modern semiparametric causal inference methodology involving (augmented) inverse propensity score weighting \citep{robins1994estimation,van2011targeted,kennedy2016semiparametric,chernozhukov2017double}. Crucially, in some of these problems such as doubly robust estimation in observational studies, nonasymptotic inference is typically not possible, and hence the CLT yields asymptotic \ci{}s for an otherwise unsolvable inference problem. 

While the CLT makes efficient statistical inference possible in a broad array of problems, the resulting CIs are only valid at a prespecified sample size $n$, invalidating any inference that occurs at data-dependent stopping times, for example under continuous monitoring.
CIs that retain validity in sequential environments are known as \emph{confidence sequences} (CS) \citep{darling1967confidence,robbins1970statistical} and can be used to make decisions at arbitrary stopping times (e.g.~while adaptively sampling, continuously peeking at the data, etc.). CSs are an inherently nonasymptotic notion, and thus essentially every published CS is nonasymptotic, including various recent state-of-the-art constructions in different settings~\citep{howard2018uniform,howard2019sequential,waudby2020estimating,wang2022catoni}. 

This paper presents a new notion: an ``asymptotic confidence sequence''. For the familiar reader, this might at first sound like an oxymoron. Further, it is not obvious how to posit a definition that is simultaneously sensible and tractable, meaning whether it is possible to develop such asymptotic CSs (whatever it may mean). We believe that we have formulated the ``right'' definition, because we accompany it with a universality result that parallels the CLT --- a universal asymptotic CS that is valid under the exact same moment assumptions required by the CLT, and exploits certain time-uniform central limit theory to arrive at boundaries that one would use if the data were Gaussian. This enables the construction of asymptotic CSs in a myriad of new situations where the distributional assumptions are weak enough to remain out of the reach of nonasymptotic techniques even in fixed-time settings. The width of this universal asymptotic CS scales with the variance of the data, just like the empirical variance used in the CLT --- such variance-adaptivity is only achievable for nonasymptotic methods in very specialized settings~\cite{waudby2020estimating}.

Before proceeding, let us first briefly review some notation and key facts about CSs.

\subsection{Time-uniform confidence sequences (CSs)}
Consider the problem of estimating the population mean $\mu = \EE (Y_1)$ from a sequence of \iid{} data $(Y_t)_{t=1}^\infty \equiv (Y_1, Y_2, \dots)$ that are observed sequentially over time. A nonasymptotic $(1-\alpha)$-CI for $\mu$ is a set\footnote{We use overhead dots $\dot C_n$ to denote fixed-time (pointwise) \ci{}s and overhead bars $\bar C_t$ for time-uniform \cs{}s.} $\dot C_n \equiv \dot C(Y_1, \dots, Y_n)$ with the property that
\begin{equation}
  \label{eq:ci}
  \forall n \in \NN^+,\ \PP(\mu \in \dot C_n) \geq 1-\alpha, ~~~~\text{or equivalently,}~~~~ \forall n \in \NN^+,\ \PP(\mu \notin \dot C_n) \leq \alpha.
\end{equation}
The coverage guarantee~\eqref{eq:ci} of a CI is only valid at some \emph{prespecified} sample size $n$, which must be decided in advance of seeing any data --- peeking at the data in order to determine the sample size is a well-known form of ``$p$-hacking''. However, it may be restrictive to fix $n$ beforehand, and even if clever sample size calculations are carried out based on prior knowledge, it is impossible to know a priori whether $n$ will be large enough to detect some signal of interest: after collecting the data, one may regret collecting too little data or collecting much more than necessary.

\begin{figure}[!htbp]
    \centering
    \includegraphics[width=\columnwidth]{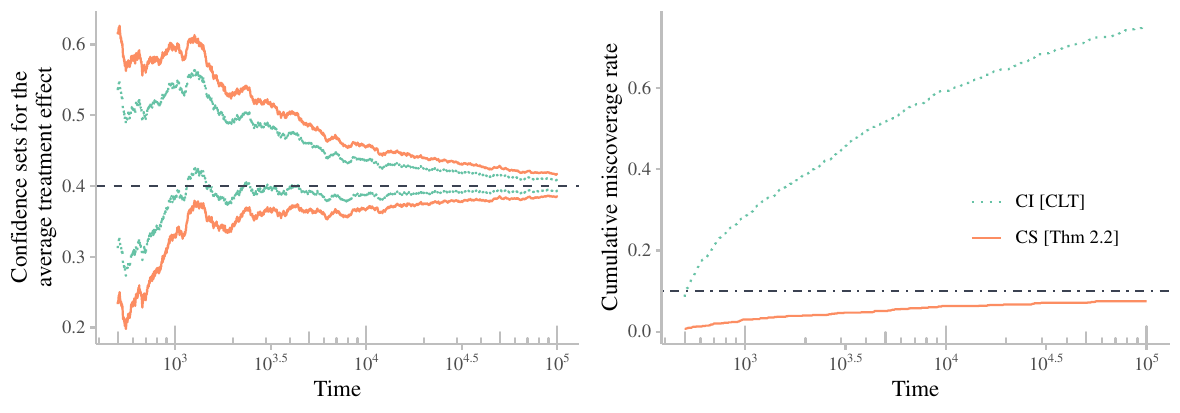}
    \caption{The left plot shows one run of a single experiment: an asymptotic CS alongside an asymptotic CI for a parameter of interest (in this case, the average treatment effect (ATE) of 0.4, an example we expand on in \cref{section:csate}). The true value of the ATE is covered by the CS simultaneously from time 30 to 10000. On the other hand, the CI fails to cover the true ATE at several points in time. By repeating such an experiment hundreds of times, one obtains the right plot which displays the cumulative probability of miscoverage  --- i.e.\ the probability of the \cs{} or CI failing to capture the true ATE at any time up to $t$. Notice that the CI error rate begins at $\alpha=0.1$ and quickly grows, while the CS error rate never exceeds $\alpha=0.1$.}
    \label{fig:CSvsCI}
  \end{figure}

CSs provide the flexibility to choose sample sizes data-adaptively while controlling the type-I error rate (see \cref{fig:CSvsCI}). Formally, a CS is a sequence of CIs $(\bar C_t)_{t=1}^\infty$ such that
\begin{equation}
  \label{eq:cs}
  \PP(\forall t \in \NN^+,\ \mu \in \bar C_t) \geq 1-\alpha, ~~~~\text{or equivalently,}~~~~ \PP(\exists t \in \NN^+: \mu \notin \bar C_t) \leq \alpha.
\end{equation}
The statements \eqref{eq:ci} and \eqref{eq:cs} look similar but are markedly different from the data analyst's or experimenter's perspective. In particular, employing a CS has the following implications:
\begin{enumerate}[(a)]
    \item The \cs{} can be (optionally) updated whenever new data become available;
    \item Experiments can be continuously monitored, adaptively stopped, or continued;
    \item The type-I error is controlled at all stopping times, including data-dependent times.
\end{enumerate}
In fact, \cs{}s may be equivalently defined  as \ci{}s that are valid at arbitrary stopping times, i.e.
\begin{equation}
  \label{eq:cs-stopping-time}
  \PP(\mu \in \bar C_\tau) \geq 1-\alpha ~~~\text{for any stopping time $\tau$}.
\end{equation}
A proof of this equivalence can be found in \citet[Lemma 3]{howard2018uniform}.

As mentioned before, while nonparametric CSs have been developed for several problems, they have thus far been \emph{nonasymptotic}.
Nonasymptotic inference for means of random variables \emph{requires} strong assumptions on the distribution of the data \citep{bahadur1956nonexistence}. These assumptions often take the form of a parametric likelihood \citep{wasserman2020universal,howard2019sequential}, known bounds on the random variables themselves \citep{howard2018uniform,waudby2020estimating}, on their moments \citep{wang2022catoni}, or on their moment generating functions~\citep{howard2018uniform}. 

These added distributional assumptions make existing CSs quite unlike CLT-based \ci{}s which (a) are universal, meaning they take the same form --- up to a change in influence functions --- and are computed in the same way for most problems, and (b) are often applicable even when no nonasymptotic \ci{} is known, such as in doubly robust inference of causal effects in observational studies.
Our work bridges this gap, bringing properties (a) and (b) to the anytime-valid sequential regime by making one simple modification to the usual \ci{}s. Just as CLT-based \ci{}s yield approximate inference for a wide variety of problems in  fixed-$n$ settings, our paper yields the same for \emph{sequential} settings.

\subsection{Contributions and outline}
We begin by rigorously defining ``asymptotic confidence sequences'' (\asympcs{}s) in \cref{definition:asympcs} and providing a general recipe to derive explicit \asympcs{}s 
that are as easy to implement and apply as CLT-based \ci{}s in \cref{section:general-framework}. Using this recipe, we develop a Lindeberg-type \asympcs{} that is able to capture time-varying means under martingale dependence (\cref{section:lindeberg-type-asympcs}). Furthermore, in \cref{section:type-I-error}, we give a definition of asymptotic time-uniform coverage (akin to coverage of asymptotic \ci{}s) and show how \emph{sequences} of our \asympcs{}s enjoy this property.
In \cref{section:csate} we illustrate how the \asympcs{}s of \cref{section:asymptoticCS} enable asymptotically anytime-valid semiparametric inference for causal effects in both randomized experiments and observational studies (\cref{section:csate}). To be clear, we are not focused on deriving new semiparametric estimators; we simply demonstrate how semiparametric causal inference --- a problem for which no known \cs{}s exist in the observational setting --- can now be tackled in fully sequential environments using the existing state-of-the-art estimators combined with our \asympcs{}s (Theorems~\ref{theorem:csate_randomized} and~\ref{theorem:csate_observational}). In \cref{section:simulations-width-coverage}, we provide a simulation study to illustrate empirical widths and miscoverage rates of \asympcs{}s and compare them to some existing (nonasymptotic) \cs{}s in the literature. Finally, we apply the \asympcs{}s of \cref{section:csate} to a real observational data set by sequentially estimating the effects of fluid intake on 30-day mortality in sepsis patients.
In sum, this work expands the scope of anytime-valid inference by tackling sequential estimation problems under CLT-like moment assumptions and guarantees.


\section{Asymptotic confidence sequences}
\label{section:time-uniform-clt}

We first define what it means for a sequence of intervals to form an asymptotic confidence sequence (\asympcs{}). Then, we derive a ``universal'' \asympcs{} in the sense that the \asympcs{} does not depend on any features of the distribution beyond its mean and variance.\footnote{We use ``universal'' in the same way that the CLT and law of large numbers are considered universal~\citep{tao2012pluribus}, as they describe macroscopic behaviors that are independent of most microscopic details of the system.} Much like classical asymptotic confidence intervals based on the CLT, this universal \asympcs{} fundamentally relies on Gaussian approximation. However, in this setting, the particular type of central limit theory being invoked is that of strong invariance principles, where an implicit Gaussian process is coupled with a partial sum with probability one (more details are provided in \cref{section:strassen-strong-approx}). Finally, similar to \ci{}s based on martingale CLTs, we derive a Lindeberg-type martingale \asympcs{} that can track a moving average of conditional means.

\subsection{Defining asymptotic confidence sequences}\label{section:asymptoticCS}

Here, we define and present ``asymptotic confidence sequences'' as time-uniform analogues of CLT-based asymptotic \ci{}s, making similarly weak moment assumptions and providing a universal closed-form boundary.

The term ``asymptotic confidence sequence'' may at first seem paradoxical. Indeed, ever since their introduction by Robbins and collaborators \citep{darling1967confidence,lai1976boundary,lai1976confidence}, \cs{}s have been defined nonasymptotically, satisfying the time-uniform guarantee in equation~\eqref{eq:cs}. So how could a bound be both time-uniform and asymptotically valid? We clarify this critical point soon, with an analogy to classical asymptotic \ci{}s. Similar to asymptotic \ci{}s, \asympcs{}s trade nonasymptotic guarantees for (a) simplicity and universality, and (b) the ability to tackle a much wider variety of problems, especially those for which there is no known nonasymptotic \cs{}. Said differently, \asympcs{}s trade finite sample validity for versatility (exemplified in \cref{section:csate} with a particular emphasis on modern causal inference).

Indeed, there is a clear desire for (asymptotically) time-uniform methods with CLT-like simplicity and versatility, especially in the context of causal inference. For example,~\citet[Section 4.3]{johari2015always} use a Gaussian mixture sequential probability ratio test (SPRT) to conduct A/B tests (i.e.~randomized experiments) for data coming from (non-Gaussian) exponential families and mentions that CLT approximations hold at large sample sizes. Similarly,~\citet{yu2020new} develop a mixture SPRT for causal effects in generalized linear models, where they say that their likelihood ratio forms an ``approximate martingale'', meaning its conditional expectation is constant up to a factor of $\exp \left \{ o_{\PP}(1) \right \}$. Moreover, \citet{pace2020likelihood} suggest using Robbins' Gaussian mixture \cs{} as a closed-form ``approximate \cs{}'' and they demonstrate through simulations that the time-uniform coverage guarantee tends to hold in the asymptotic regime. However, all of these approaches justify time-uniform inference with $o_{\PP}(\cdot)$ approximations that only hold at a~\emph{fixed, pre-specified sample size}, and yet inferences are being carried out at \emph{data-dependent sample sizes}. This section remedies the tension between fixed-$n$ approximations and time-uniform inference by defining \asympcs{}s such that Gaussian approximations must hold almost surely for \emph{all sample sizes simultaneously}. The \asympcs{}s we define will also be valid in a wide range of nonparametric scenarios (beyond exponential families, parametric models, and so on).


To motivate the definition of an \asympcs{} that follows, let us briefly review the CLT in the batch (non-sequential) setting. Suppose $Y_1, \dots, Y_n \simiid \PP$ with mean $\EE(Y_1) = \mu$ and variance $\Var(Y_1) = \sigma^2$. Then the standard CLT-based \ci{} for $\mu$ (with known variance $\sigma$) takes the form
\begin{equation}
  \label{eq:CLT}
  \dot C_n := [ \widehat \mu_n \pm \dot \boundary_n ]\equiv \left [\widehat \mu_n \pm \sigma \cdot \frac{\Phi^{-1}(1-\alpha/2)}{\sqrt{n}}\right ] ,
\end{equation}
where $\widehat \mu_n$ is the sample mean and $\Phi^{-1}(1-\alpha/2)$ is the $(1-\alpha/2)$-quantile of a standard Gaussian $N(0, 1)$ (e.g. for $\alpha =0.05$, we have $\Phi^{-1}(0.975) \approx 1.96$).
The classical notion of ``asymptotic validity'' is 
\begin{equation}\label{eq:asymp-CI}
\liminf_{n\to\infty} \PP(\mu \in \dot C_n) \geq 1- \alpha.
\end{equation} 
While the above is the standard definition of an asymptotic CI, one could have arrived at an alternative definition by noting the following rather strong statement that can be made under the same conditions: there exist\footnote{Technically, writing \eqref{eq:strong-coupling-in-probability} may require enriching the probability space so that both $Y$ and $Z$ can be defined (but without changing their laws). See \citet[Equation (1.2)]{einmahl2009new} for a precise statement.} \iid{} standard Gaussians $Z_1, \dots, Z_n \sim N(0, 1)$ such that
\begin{equation}
  \label{eq:strong-coupling-in-probability}
  \frac{1}{n}\sum_{i=1}^n(Y_i - \mu)/\sigma = \frac{1}{n}\sum_{i=1}^n Z_i + o_{\PP}(1/\sqrt{n}).
\end{equation}
From this vantage point, we note that the \ci{} in \eqref{eq:CLT} has the additional guarantee that there in fact exists an (unknown)  \emph{nonasymptotic} $(1-\alpha)$-\ci{} $[\widehat \mu_n \pm \dot \boundary_n^\star]$ such that
\begin{equation}
  \label{eq:CI-convergence-in-prob}
  \dot \boundary_n^\star / \dot \boundary_n \xrightarrow{\PP} 1.
\end{equation}
We deliberately highlight the above property of asymptotic \ci{}s because it ends up serving as a natural starting point for defining asymptotic confidence \emph{sequences}.
In particular, we will define \asympcs{}s so that an analogous approximation to~\eqref{eq:CI-convergence-in-prob} holds uniformly over time, almost surely.
Statements like \eqref{eq:strong-coupling-in-probability} are known as ``couplings'' and appear in the literature on strong approximations and invariance principles where similar guarantees can indeed be shown to hold almost surely and at faster rates under additional moment assumptions \citep{einmahl2009new,komlos1975approximation,komlos1976approximation}.
\begin{definition}[Asymptotic confidence sequences]
  \label{definition:asympcs}
  Let $\Tcal$ be a totally ordered infinite set (denoting time) that has a minimum value $t_0 \in \Tcal$.
  We say that the intervals $(\widehat \theta_t - L_t, \widehat \theta_t + U_t)_{t \in \Tcal}$ centered at the estimators $(\widehat \theta_t)_{t \in \Tcal}$ with non-zero bounds $L_t, U_t > 0, \forall t \in \Tcal$ form a $(1-\alpha)$-\emph{asymptotic confidence sequence} (\asympcs{}) for a sequence of real parameters $\infseqtcal{\theta_t}$  if there exists a (typically unknown) nonasymptotic $(1-\alpha)$-\cs{} $(\widehat \theta_t - L_t^\star, \widehat \theta_t + U_t^\star)_{t \in \Tcal}$ for $\infseqtcal{\theta_t}$ --- i.e.~satisfying
  \begin{equation}
    \PP \left (\forall t \in \Tcal,\ \theta_t \in [ \widehat \theta_t - L_t^\star, \widehat \theta_t + U_t^\star] \right ) \geq 1-\alpha,
  \end{equation}
  and such that $L_t, U_t$ become arbitrarily precise almost-sure approximations to $L_t^\star$ and $U_t^\star$:   \begin{equation}
    \label{eq:asympcs}
    L_t^\star / L_t \xrightarrow{\as} 1\quad\text{and}\quad U_t^\star / U_t\xrightarrow{\as} 1.
  \end{equation}
\end{definition}
In words, Definition~\ref{definition:asympcs} says that an \asympcs{} $(C_t)_{t\in \Tcal}$ centered at $(\widehat \theta_t)_{t\in \Tcal}$ is an arbitrarily precise approximation of some nonasymptotic \cs{} $(C_t^\star)_{t\in \Tcal}$ centered at $(\widehat \theta_t)_{t\in \Tcal}$ as $t \to \infty$. Throughout the paper, we will mostly focus on the case where $\Tcal = \NN$ with $t_0 = 0$.

It is important to note that alternate definitions fail to be coherent in different ways.
As one example, we could have hypothetically defined a sequence of intervals $(C_t(\alpha))_{t\in \Tcal}$ to be a $(1-\alpha)$-\asympcs{} if $\limsup_{m \to \infty} \PP(\exists t \geq m: \mu \notin C_t(\alpha)) \leq \alpha$, analogously to asymptotic \ci{}s which satisfy $\limsup_{n \to \infty} \PP(\mu \notin \dot C_n(\alpha)) \leq \alpha$. In words, we could have posited that if we just start peeking late enough, then the probability of eventual miscoverage would indeed be below $\alpha$. Unfortunately, even for nonasymptotic \cs{}s constructed at any level $\alpha' \in (0,1)$, the former limit is \emph{zero}, so this inequality would be vacuously true, regardless of $\alpha'$, even if $\alpha' \gg \alpha$.
However, we do show that \emph{sequences} of our \asympcs{}s have a guarantee of this type if peeking starts late enough (see \cref{section:type-I-error}), but we delay these definitions until later as they are slightly more involved. 

  By virtue of being defined in terms of their limiting behavior, one can obviously construct \asympcs{}s (as well as asymptotic confidence intervals) with nonsensical finite-sample behavior. It is thus imperative that if a practitioner decides to employ an asymptotic method, they do so with the understanding that its effectiveness relies on it exploiting some well-approximated nonasymptotic phenomenon, and that its limiting behavior should be viewed as a rigorous manifestation of a guiding principle rather than a panacea.

\begin{remark}[Why almost surely?]\label{remark:why-almost-surely}
  One may wonder why it is necessary to define \asympcs{}s so that $L_t^\star / L_t \to 1$ \emph{almost surely} (rather than in probability, for example). Since \cs{}s are bounds that hold uniformly over time with high probability, convergence in probability $L_t^\star / L_t = 1 + o_\PP(1)$ is not the right notion of convergence, as it only requires that the approximation term $o_\PP(1)$ be small with high probability for sufficiently large \emph{fixed} $t$, but not for all $t$ uniformly. It is natural to try to extend convergence in probability to \emph{time-uniform convergence with high probability} --- i.e. $\sup_{k \geq t} (L_k^\star / L_k) = 1 + o_\PP(1)$ --- but it turns out that this is in fact \emph{equivalent} to almost-sure convergence $L_t^\star / L_t = 1+ o_\as (1)$; see \cref{section:TimeUniformInProb_AlmostSure}. 
\end{remark}

Going forward, we may omit ``a.s.''~from $\oas{\cdot}$ and $\Oas{\cdot}$ and instead simply write $o(\cdot)$ and $O(\cdot)$, respectively to simplify notation.
Now that we have defined \asympcs{}s as time-uniform analogues of asymptotic \ci{}s, we will explicitly derive \asympcs{}s for the mean of \iid{} random variables with finite variances (i.e.~under the same assumptions as the CLT).

\subsection{Warmup:~\asympcs{}s for the mean of \iid{} random variables}\label{section:universal-asympcs-mean}
  
We now construct an explicit \asympcs{} for the mean of \iid{} random variables by combining a variant of Robbins' (nonasymptotic) Gaussian mixture boundary \citep{robbins1970statistical} with Strassen's strong approximation theorem \citep{strassen1964invariance}. Before presenting the result, let us review Robbins' boundary and Strassen's result, and discuss how they can be used in conjunction to arrive at the \asympcs{} in \cref{theorem:acs}.

\subsubsection{Robbins' Gaussian mixture boundary}\label{section:robbins-gaussian-mixture}

The study of \cs{}s began with Herbert Robbins and colleagues \citep{darling1967confidence,robbins1970statistical,robbins1970boundary,lai1976boundary,lai1976confidence}, leading to several fundamental results and techniques including the famous Gaussian mixture boundary for partial sums of \iid{} Gaussian random variables \citep{robbins1970statistical} (see also \citet[\S 3.2]{howard2018uniform}) which we recall here. Suppose $\infseqt{Z_t}$ are \iid{} standard Gaussian random variables. Then for any $\rho > 0$,
\begin{equation}\label{eq:robbins1970}
  \PP \left ( \exists t \geq 1 : \left \lvert  \frac{1}{t}\sum_{i=1}^t Z_i \right \rvert \geq \sqrt{\frac{2(t\rho^2 + 1)}{t^2 \rho^2} \log \left ( \frac{\sqrt{t \rho^2 + 1}}{\alpha} \right )} \right ) \leq \alpha.
\end{equation}
Notice that the above boundary scales as $O(\sqrt{\log t / t})$ for any $\rho > 0$. In fact, \citet[Eq. 11]{robbins1970statistical} noted that \eqref{eq:robbins1970} holds not only for Gaussian random variables, but for those that are $1$-sub-Gaussian, and hence pre-multiplying the boundary by $\sigma$ yields a $\sigma$-sub-Gaussian time-uniform concentration inequality, serving as a time-uniform analogue of Chernoff or Hoeffding inequalities. The connections between these fixed-time and time-uniform concentration inequalities are made explicit in \citet{howard2018uniform}. Nevertheless, \cref{eq:robbins1970} requires \emph{a priori} knowledge of $\sigma > 0$ unlike CLT-based \ci{}s which we aim to emulate in the (asymptotically) time-uniform regime. The following strong Gaussian approximation due to \citet{strassen1964invariance} will serve as a technical tool allowing the nonasymptotic sub-Gaussian bound in \eqref{eq:robbins1970} to be applied to partial sums of arbitrary random variables with finite variances.

\subsubsection{Strassen's strong approximation}\label{section:strassen-strong-approx}

\citet{strassen1964invariance} initiated the study of ``strong approximation'' (also called strong invariance principles or strong embeddings) which blossomed into an active and impactful corner of probability theory research over the subsequent years, culminating in now-classical results such as the Koml\'os-Major-Tusn\'ady embeddings \citep{komlos1975approximation,komlos1976approximation,major1976approximation} and other related works \citep{strassen1967and,einmahl1989extensions,morrow1982almost,chatterjee2012new}.

In his 1964 paper, \citet[\S 2]{strassen1964invariance} used the Skorokhod embedding \citep{skorokhod1961research} (see also \citep[p. 513]{billingsley1995probability}) to obtain a strong invariance principle which connects asymptotic Gaussian behavior with the law of the iterated logarithm. Concretely, let $\infseqt{Y_t}$ be an infinite sequence of \iid{} random variables from a distribution $\PP$ with mean $\mu$ and variance $\sigma^2$. Then, on a potentially richer probability space,\footnote{A richer probability space may be needed to describe Gaussian random variables, if for example, $\infseqt{Y_t}$ are $\{0, 1\}$-valued on a probability space whose probability measure is dominated by the counting measure. This construction of a richer probability space imposes no additional assumptions on $\infseqt{Y_t}$, and is only a technical device used to rigorously couple two sequences of random variables, and appears in essentially all papers on strong invariance principles, not just \citet{strassen1964invariance}.} there exist standard Gaussian random variables $\infseqt{Z_t}$ whose partial sums are almost-surely coupled with those of $\infseqt{Y_t}$ up to iterated logarithm rates, i.e.
\begin{equation}\label{eq:strassen1964}
  \left \lvert \sum_{i=1}^t (Y_i - \mu) / \sigma - \sum_{i=1}^t Z_i \right \rvert = o \left ( \sqrt{t \log \log t} \right ) \quad \text{almost surely.}
\end{equation}
Notice that the law of the iterated logarithm states that $|\sum_{i=1}^t (Y_i - \mu) / \sigma| = O(\sqrt{t \log \log t})$ while \eqref{eq:strassen1964} states that the same partial sum is almost-surely coupled with an implicit Gaussian process --- i.e.~replacing $O(\cdot)$ with $o(\cdot)$. It may be convenient to divide by $t$ and interpret~\eqref{eq:strassen1964} on the level of sample averages rather than partial sums, in which case the right-hand side becomes $o (\sqrt{\log \log t / t})$. 
Let us now describe how Strassen's strong approximation can be combined with Robbins' Gaussian mixture boundary to derive an \asympcs{} under finite moment assumptions akin to the CLT\@.

\subsubsection{The Gaussian mixture asymptotic confidence sequence}

Given the juxtaposition of \eqref{eq:robbins1970} and \eqref{eq:strassen1964}, the high-level approach to the derivation of \asympcs{}s becomes clearer. Indeed, the essential idea behind \cref{theorem:acs} is as follows. By Strassen's strong approximation, we couple the partial sums $S_t := \sum_{i=1}^t (Y_i - \mu) / \sigma$ with implicit partial sums $G_t := \sum_{i=1}^t Z_i$ of Gaussians $\infseqt{Z_t}$, and then use Robbins' mixture boundary to obtain a time-uniform high-probability bound on the deviations of $|G_t|$, noting that the coupling rate $o(\sqrt{t \log \log t})$ is asymptotically dominated by the concentration rate $O(\sqrt{t \log t})$, leading to asymptotic validity in the formal sense of \cref{definition:asympcs}.


\begin{theorem}[Gaussian mixture asymptotic confidence sequence]\label{theorem:acs}
  Suppose $(Y_t)_{t=1}^\infty \simiid \PP$ is an infinite sequence of \iid{} observations from a distribution $\PP$ with mean $\mu$ and finite variance. Let $\widehat \mu_t := \frac{1}{t} \sum_{i=1}^t Y_i$ be the sample mean, and $\widehat \sigma_t^2 := \frac{1}{t}\sum_{i=1}^t Y_i^2 - (\widehat \mu_t)^2$ the sample variance based on the first $t$ observations. Then, for any prespecified constant $\rho > 0$,
  \begin{equation}
    \label{eq:Gaussian-mixture-asympcs}
    \widebar C_t^\Gcal \equiv (\widehat \mu_t \pm \widebar \boundary_t^\Gcal) := \left ( \widehat \mu_t \pm \widehat \sigma_t \sqrt {\frac{2(t\rho^2 + 1)}{t^2 \rho^2} \log\left (\frac{\sqrt{t\rho^2 + 1}}{\alpha} \right)} \right )
  \end{equation}
  forms a $(1-\alpha)$-\asympcs{} for $\mu$.

\end{theorem}
The proof of \cref{theorem:acs} is in Appendix~\ref{proof:acs}.
We can think of $\rho > 0$ as a user-chosen tuning parameter which dictates the time at which \eqref{eq:Gaussian-mixture-asympcs} is tightest, and we discuss how to easily tune this value in \cref{section:optimizingMixture}. A one-sided analogue of~\eqref{eq:Gaussian-mixture-asympcs} can be found in \cref{section:one-sided}.

While \eqref{eq:Gaussian-mixture-asympcs} may look visually similar to Robbins' (sub)-Gaussian mixture CS \citep{robbins1970statistical} --- written explicitly in \citet[Eq. (14)]{howard2018uniform} --- it is worth pausing to reflect on how they are markedly different. Firstly, Robbins' CS is a nonasymptotic bound that is only valid for $\sigma$-sub-Gaussian random variables, meaning $\EE \exp \{ \lambda (Y_1 - \EE Y_1 ) \} \leq \exp \{ \sigma^2\lambda^2/2\}$ for some \emph{a priori} known $\sigma > 0$, while Theorem~\ref{theorem:acs} does not require the existence of a finite MGF at all (much less a known upper bound on it). Secondly, Robbins' CS uses this known (possibly conservative) $\sigma$ in place of $\widehat \sigma_t$ in \eqref{eq:Gaussian-mixture-asympcs}, and thus it cannot adapt to an unknown variance, while \eqref{eq:Gaussian-mixture-asympcs} always scales with $\sqrt{\Var(Y_1)}$. In simpler terms, Theorem~\ref{theorem:acs} is an asymptotically time-uniform analogue of the CLT in the same way that Robbins' CS is a time-uniform analogue of a sub-Gaussian concentration inequality (e.g. Hoeffding's or Chernoff's inequality \citep{hoeffding_probability_1963,howard_exponential_2018}).

It is important not to confuse Theorem~\ref{theorem:acs} with a martingale CLT as the latter still gives fixed-time CIs in the spirit of the usual CLT but under different assumptions on the martingale difference sequence (however, we do present an analogue of \cref{theorem:acs} under martingale dependence in \cref{theorem:lindeberg-martingale-asympcs}).

\subsubsection{An asymptotic confidence sequence with iterated logarithm rates}\label{section:lil}

As a consequence of the law of the iterated logarithm, a confidence sequence for $\mu$ centered at $\widehat \mu_t$ cannot have an asymptotic width smaller than $O( \sqrt{\log \log t / t})$. This is easy to see since
\[ \limsup_{t\to \infty} \frac{\sqrt{t}|\widehat \mu_t - \mu|}{\sigma \sqrt{2 \log \log t}} = 1. \]
This raises the question as to whether $\widebar C_t^\Gcal$ can be improved so that the optimal asymptotic width of $O( \sqrt{ \log \log t / t} )$ is achieved. Indeed, we can replace Robbins' Gaussian mixture boundary with \citet[Eq. (2)]{howard2018uniform} (or virtually any other Gaussian boundary for that matter) in the proof of \cref{theorem:acs} to derive such an \asympcs{}, but as the authors discuss, mixture boundaries such as the one in Theorem~\ref{theorem:acs} may be preferable in practice, because any bound that is tighter ``later on'' (asymptotically) must be looser ``early on'' (at practical sample sizes) due to the fact that all such bounds have a cumulative miscoverage probability $\leq \alpha$. This is formally a concern for nonasymptotic \cs{}s, but only applies to \asympcs{}s insofar as they are asymptotic approximations of nonasymptotic bounds. Nevertheless, we present an \asympcs{} with an iterated logarithm rate here for completeness.
\begin{proposition}[Iterated logarithm asymptotic confidence sequences]
  \label{proposition:acsLIL}
  Under the same conditions as Theorem~\ref{theorem:acs},
  \[ \widebar C_t^\Lcal \equiv (\widehat \mu_t \pm \widebar \boundary_t^\Lcal) := \left ( \widehat \mu_t \pm \widehat \sigma_t \cdot 1.7\sqrt{ \frac{\log \log (2t) + 0.72\log(10.4/\alpha)}{t} } \right ) \]
  forms a $(1-\alpha)$-\asympcs{} for $\mu$.
\end{proposition}
We omit the proof of~\cref{proposition:acsLIL} as it proceeds in a similar fashion to that of~\cref{theorem:acs}. In fact, both of these \asympcs{}s are simply instantiations of a more general recipe for deriving \asympcs{}s by combining strong approximations with time-uniform boundaries for the approximating process, an approach that we discuss in the following section.

\subsection{A general recipe for deriving asymptotic confidence sequences}\label{section:general-framework}

The proofs of both \cref{theorem:acs} and \cref{proposition:acsLIL} follow the same general structure, combining strong approximations with time-uniform boundaries along with some other almost-sure asymptotic behavior. Abstracting away the details specific to these particular results, we provide the following four general conditions under which many \asympcs{}s can be derived, including those from the previous section but also Lyapunov- and Lindeberg-type \asympcs{} that we will state in \cref{section:lindeberg-type-asympcs}).

In what follows, let $\Tcal$ be a totally ordered infinite set that includes a minimum value $t_0 \in \Tcal$ (for example, one may think about $\Tcal$ as $\RR^{\geq 0}$ or $\NN$ with $t_0 = 0$) and let $\infseqtcal{\widehat \theta_t}$ be a sequence of estimators for the real-valued parameters $\infseqtcal{\theta_t}$. Then, consider the following four conditions where we use the ``Condition \underline{G}-\underline{X}'' enumeration as a mnemonic for the \underline{X}$^\text{th}$ condition in the section on \underline{G}eneral recipes for \asympcs{}s).

\begin{conditiongeneral}[Strong approximation]\label{condition:strong-approx}
  On a potentially enriched probability space, there exists a process $(Z_t)_{t \in \Tcal}$ starting at $Z_{t_0} \equiv 0$ that strongly approximates $(\widehat \theta_t - \theta_t)_{t \in \Tcal}$ up to a rate of $(r_t)_{t \in \Tcal}$, i.e.
  \begin{equation}\label{eq:abstract-strong-approx}
    (\widehat \theta_t - \theta_t) - Z_t = O \left ( r_t \right )\quad\text{almost surely}.
  \end{equation}
\end{conditiongeneral}

\begin{conditiongeneral}[Boundary for the approximating process]\label{condition:boundary-for-approximating-process}
  There exist $\widehat L_t > 0$ and $\widehat U_t > 0$ for each $t \in \Tcal$ so that $[- \widehat L_t, \widehat U_t]_{t \in \Tcal}$ forms a $(1-\alpha)$-boundary for the process $(Z_t)_{t \in \Tcal}$ given in \eqref{eq:abstract-strong-approx}:
  \begin{equation}\label{eq:abstract-boundary}
    \PP \left ( \forall t \in \Tcal,\ Z_t \in [- \widehat L_t, \widehat U_t] \right ) \geq 1-\alpha.
  \end{equation}
\end{conditiongeneral}

\begin{conditiongeneral}[Strong approximation rate]\label{condition:strong-approx-rate}
  The approximation rate $(r_t)_{t \in \Tcal}$ in \eqref{eq:abstract-strong-approx} is faster than both $(\widehat L_t)_{t \in \Tcal}$ and $(\widehat U_t)_{t \in \Tcal}$ in~\eqref{eq:abstract-boundary}, i.e.
  \begin{equation}
    r_t = o\left(\widehat L_t \land \widehat U_t \right )\quad\text{almost surely.}
  \end{equation}
\end{conditiongeneral}

\begin{conditiongeneral}[Almost-sure approximate boundary]\label{condition:boundary-approximation}
  The $(1-\alpha)$-boundary $[-\widehat L_t,\widehat U_t]_{t \in \Tcal}$ for $(Z_t)_{t\in \Tcal}$ is almost-surely approximated by the sequences $[-L_t, U_t]_{ t\in \Tcal}$, i.e.
  \begin{equation}
    L_t / \widehat L_t \xrightarrow{\as} 1 \quad \text{and} \quad U_t / \widehat U_t \xrightarrow{\as} 1.
  \end{equation}
\end{conditiongeneral}

Deriving new \asympcs{}s then reduces to the conceptually simpler but nevertheless nontrivial task of satisfying the requisite conditions above. For example, in the previous section, we satisfied \cref{condition:strong-approx} via \citet{strassen1964invariance}, \cref{condition:boundary-for-approximating-process} via \citet{robbins1970statistical}, \cref{condition:strong-approx-rate} via the combination of \citet{strassen1964invariance} and \citet{robbins1970statistical}, and \cref{condition:boundary-approximation} via the strong law of large numbers (SLLN). The only difference between \cref{theorem:acs} and \cref{proposition:acsLIL} was in what boundaries were being used for $[L_t, U_t]_{t \in \Tcal}$ and $[\widehat L_t, \widehat U_t]_{t\in \Tcal}$. More generally under \namecref{condition:boundary-approximation}s~\labelcref{condition:strong-approx}--\labelcref{condition:boundary-approximation}, we have the following abstract \namecref{theorem:general-asympcs} for \asympcs{}s.

\begin{theorem}[An abstract \asympcs{} for well-approximated processes]\label{theorem:general-asympcs}
  Let $\Tcal$ be a totally ordered infinite set containing a minimal element $t_0\in \Tcal$ and let $(\widehat \theta_t)_{t \in \Tcal}$ be a real-valued process. Under \namecref{condition:strong-approx}s~\ref{condition:strong-approx}--\ref{condition:boundary-approximation}, 
  \begin{equation}
    \left [ \widehat \theta_t - L_t, \widehat \theta_t + U_t \right ]
  \end{equation}
  forms a $(1-\alpha)$-\asympcs{} for $\theta_t$ meaning there exists (on a potentially enriched probability space) some nonasymptotic $(1-\alpha)$-\cs{} $[\widehat \theta_t - L_t^\star, \widehat \theta_t +U_t^\star]_{t \in \Tcal}$ for $(\theta_t)_{t \in \Tcal}$, i.e. 
  \begin{equation}\label{eq:abstract-nonasymptotic-cs}
    \PP \left (\forall t \in \Tcal,\ \theta_t \in \left [\widehat \theta_t - L_t^\star, \widehat \theta_t + U_t^\star \right ]\right ) \geq 1-\alpha
  \end{equation}
  such that
  \begin{equation}
    L_t^\star/  L_t \xrightarrow{\as} 1\quad\text{and}\quad U_t^\star/  U_t \xrightarrow{\as} 1.
  \end{equation}
\end{theorem}
We provide a short proof of \cref{theorem:general-asympcs} in \cref{proof:general-asympcs}. Note that the lower boundaries given by $L_t^\star$ and $\widehat L_t$ are not the same, but rather $L_t^\star$ is constructed from $\widehat L_t$ (and similarly for $U_t^\star$ and $\widehat U_t$.
 In the following section, we will use the general recipe of \cref{theorem:general-asympcs} to obtain \asympcs{}s for time-varying means from non-\iid{} random variables under martingale dependence akin to the Lindeberg CLT \citep{lindeberg1922neue,billingsley1995probability}.

\subsection{Lindeberg- and Lyapunov-type \asympcs{}s for time-varying means}
\label{section:lindeberg-type-asympcs}
\begin{figure}[!htbp]
  \centering
  \includegraphics[width=0.9\columnwidth]{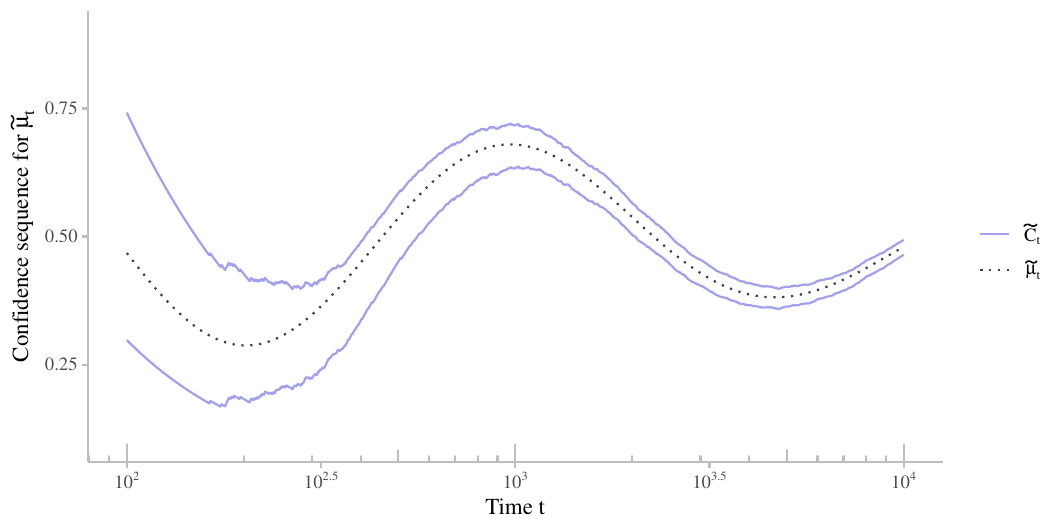}
  \caption{A $90\%$-\asympcs{} for the time-varying mean $\widetilde \mu_t$ using~\cref{theorem:lindeberg-martingale-asympcs} with $\rho$ optimized for $t^\star = 500$ based on the exact solution of~\cref{section:optimizingMixture}. Here, we have set $\mu_t := \frac{1}{2}(1 - \sin(2 \log(e + 10t)) / \log(e + 0.01t))$ to produce the sinusoidal behavior of $\widetilde \mu_t$. Notice that $\widetilde C_t$ uniformly captures $\widetilde \mu_t$, adapting to its non-stationarity.}
  \label{fig:cs_wavy}
\end{figure}

The results in \cref{theorem:acs} and \cref{proposition:acsLIL} focused on the situation where the observed random variables are independent and identically distributed, as this is one of the most commonly studied regimes in statistical inference. One may also be interested in the case where means and variances do not remain constant over time, or where observations are dependent. We will now show that an analogue of Theorem~\ref{theorem:acs} holds for random variables with time-varying \emph{means and variances} under \emph{martingale dependence}. In this case, rather than the \asympcs{} covering some fixed $\mu$, it covers the average conditional mean thus far: $\widetilde \mu_t := \frac{1}{t} \sum_{i=1}^t \mu_i$ --- to be made precise shortly.\footnote{Throughout this section and the remainder of the paper, we use the overhead tilde (e.g. $\widetilde \mu_t$, $\widetilde \sigma_t$, and $\widetilde C_t$) to emphasize that these quantities can change over time. For example, \cref{fig:cs_wavy} explicitly displays means and \cs{}s with sinusoidal behaviors resembling a tilde.}

Given the additional complexity introduced by considering time-varying conditional distributions, we will first explicitly spell out the conditions required to achieve a time-varying analogue of Theorem~\ref{theorem:acs}.
Suppose $(Y_t)_{t=1}^\infty$ is a sequence of random variables with conditional means and variances given by $\mu_t := \EE (Y_t \mid Y_1^{t-1})$ and $\sigma_t^2 := \Var(Y_t \mid Y_1^{t-1})$, respectively where we use the shorthand $Y_1^{t-1}$ for $\{Y_1, \dots, Y_{t-1}\}$. First, we require that the average conditional variance $\widetilde \sigma_t^2 := \frac{1}{t} \sum_{i=1}^t \sigma_i^2$ either does not vanish, or does so superlinearly; equivalently, we require that the cumulative conditional variance diverges almost surely.
\begin{conditionlind}[Cumulative variance diverges almost surely]
  \label{assumption:lindeberg_variance_does_not_vanish}
  For each $t \geq 1$, let $\sigma_t^2 := \Var(Y_t \mid Y_1^{t-1})$ be the conditional variance of $Y_t$. Then,
  \begin{equation}
    \label{eq:lindeberg_variance_vanish_superlinear}
    V_t := \sum_{i=1}^t \sigma_i^2 \to \infty~~\text{almost surely}.
  \end{equation}
\end{conditionlind}

\cref{eq:lindeberg_variance_vanish_superlinear} can also be interpreted as saying that the average conditional variance $\widetilde \sigma_t^2 := \frac{1}{t}\sum_{i=1}^t \sigma_i^2$ does not vanish faster than $1/t$ (if at all), meaning $\widetilde \sigma_t^2 = \omega_{\as}(1/t)$. For example,~\cref{assumption:lindeberg_variance_does_not_vanish} would hold if $\widetilde \sigma_t^2 \xrightarrow{\as} \sigma_\star^2$ for some $\sigma_\star^2 > 0$ or in the \iid{} case where $\sigma^2_1 = \sigma^2_2 = \cdots = \sigma_\star^2$.
Second, we require a Lindeberg-type uniform integrability condition on the tail behavior of $\infseqt{Y_t}$.
\begin{conditionlind}[Lindeberg-type uniform integrability]
  \label{assumption:lindeberg-type_condition}
  For $t \geq 1$, let $\sigma_t^2 := \Var(Y_t \mid Y_1^{t-1})$ be the conditional variance of $Y_t$. Then there exists some $0 < \kappa < 1$ such that
  \begin{equation}
    \label{eq:lindeberg_type_condition}
    \sum_{t=1}^\infty \frac{\EE \left [ (Y_t - \mu_t)^2 \1 \left ((Y_t-\mu_t)^2 > V_t^{\kappa} \right ) \mid Y_1^{t-1} \right ] }{V_t^{\kappa}} < \infty ~~\text{almost surely.}
  \end{equation}
\end{conditionlind}
Notice that~\cref{eq:lindeberg_type_condition} is satisfied if all conditional $q^\mathrm{th}$ moments are almost surely uniformly bounded for some $q > 2$, meaning $1/K \leq \EE (|Y_t - \mu_t|^q \mid Y_1^{t-1}) < K$ a.s.~for all $t \geq 1$ and for some constant $K > 0$, or more generally under a Lyapunov-type condition that states $\sum_{t=1}^\infty \left [\EE(|Y_t - \mu_t|^{2+\delta} \mid Y_1^{t-1}) / \sqrt{V_t}^{2+\delta} \right ] < \infty$ a.s.~for some $\delta  > 0$.\footnote{We show that the Lyapunov-type condition implies \cref{assumption:lindeberg-type_condition} in \cref{section:lyapunov-implies-lindeberg}.}
Third and finally, we require a consistent variance estimator.
\begin{conditionlind}[Consistent variance estimation]
  \label{assumption:lindeberg_variance_estimator_consistent}
  Let $\widehat \sigma_t^2$ be an estimator of $\widetilde \sigma_t^2$ constructed using $Y_1, \dots, Y_t$ such that
  \begin{equation}
    \label{eq:lyapunov_variance_o1}
    \widehat \sigma_t^2 / \widetilde \sigma_t^2 \xrightarrow{\mathrm{a.s.}} 1.
  \end{equation}
\end{conditionlind}
Note that in the \iid{} case,~\eqref{eq:lyapunov_variance_o1} would hold using the sample variance by the SLLN. More generally for independent but non-identically distributed data, \cref{assumption:lindeberg_variance_estimator_consistent} holds as long as the variation in means vanishes --- i.e.~$\frac{1}{t}\sum_{i=1}^t (\mu_i - \widetilde \mu_t)^2 = o(1)$ --- but we will expand on this later in \cref{corollary:lyapunov-asympcs}. Given \namecref{assumption:lindeberg_variance_does_not_vanish}s~\ref{assumption:lindeberg_variance_does_not_vanish}, \ref{assumption:lindeberg-type_condition}, and \ref{assumption:lindeberg_variance_estimator_consistent}, we have the following \asympcs{} for the time-varying conditional mean $\widetilde \mu_t := \frac{1}{t} \sum_{i=1}^t \mu_i$.

\begin{proposition}[Lindeberg-type Gaussian mixture martingale \asympcs{}]\label{theorem:lindeberg-martingale-asympcs}
  Let $\infseqt{Y_t}$ be a sequence of random variables with conditional mean $\mu_t := \EE(Y_t \mid Y_1^{t-1})$ and conditional variance $\sigma_t^2 := \Var(Y_t\mid Y_1^{t-1})$. Then under Assumptions~\ref{assumption:lindeberg_variance_does_not_vanish},~\ref{assumption:lindeberg-type_condition},~and~\ref{assumption:lindeberg_variance_estimator_consistent}, we have that
  \begin{equation}
    \label{eq:lyapunov-style-asympcs}
    \widetilde C_t \equiv (\widehat \mu_t \pm \widetilde \boundary_t) := \left (\widehat \mu_t \pm \sqrt{\frac{2(t \widehat \sigma_t^2 \rho^2 + 1)}{t^2 \rho^2} \log \left ( \frac{\sqrt{t \widehat \sigma_t^2 \rho^2 + 1}}{\alpha} \right )} \right )
  \end{equation}
  forms a $(1-\alpha)$-\asympcs{} for the running average conditional mean $\widetilde \mu_t := \frac{1}{t}\sum_{i=1}^t \mu_i$.
\end{proposition}

At a high level, the proof of \cref{theorem:lindeberg-martingale-asympcs} (found in \cref{proof:lyapunov-type-asympcs}) follows from the general \asympcs{} procedure of \cref{theorem:general-asympcs} by using \cref{assumption:lindeberg-type_condition} and Strassen's \citeyear{strassen1967and} strong approximation \citep{strassen1967and} (not to be confused with his \citeyear{strassen1964invariance} result that we used in \cref{theorem:acs}) to satisfy Conditions~\ref{condition:strong-approx} and~\ref{condition:strong-approx-rate}, and a variant of Robbins' mixture martingale for non-\iid{} random variables along with Conditions~\ref{assumption:lindeberg_variance_does_not_vanish} and~\ref{assumption:lindeberg_variance_estimator_consistent} to satisfy Conditions~\ref{condition:boundary-for-approximating-process} and~\ref{condition:boundary-approximation}.

Notice that if the data happen to be \iid{}, then $\infseqt{\widetilde C_t}$ is asymptotically equivalent to $\infseqt{\widebar C_t^\Gcal}$ given in \cref{theorem:acs} (here, ``asymptotic equivalence'' simply means that the ratio of the two boundaries converges a.s.~to 1). In other words, \cref{theorem:lindeberg-martingale-asympcs} is valid in a more general (non-\iid{}) setting, but will essentially recover \cref{theorem:acs} in the \iid{} case. \cref{fig:cs_wavy} illustrates what $\widetilde C_t$ may look like in practice. Note that when $\infseqt{Y_t}$ are independent with $\mu_1 = \mu_2 = \cdots = \mu_\star$, and $\sigma_1^2 = \sigma_2^2 = \cdots = \sigma_\star^2$, it is nevertheless the case that $\widetilde C_t$ forms a $(1-\alpha)$-\asympcs{} for $\mu_\star$ under the same assumptions as Theorem~\ref{theorem:acs}. In this sense, we can view $(\widetilde C_t)_{t=1}^\infty$ as ``robust'' to deviations from independence and stationarity.\footnote{Here, the term ``robust'' should not be interpreted in the same spirit as ``doubly robust'', where the latter is specific to the discussions surrounding functional estimation and causal inference in \cref{section:csate}.} A one-sided analogue of~\cref{theorem:lindeberg-martingale-asympcs} is presented in~\cref{proposition:martingale-asympcs-one-sided} within~\cref{section:one-sided}.

As a near-immediate corollary of~\cref{theorem:lindeberg-martingale-asympcs}, we have the following Lyapunov-type \asympcs{} under independent but non-identically distributed random variables.

\begin{corollary}[Lyapunov-type \asympcs{}]\label{corollary:lyapunov-asympcs}
  Suppose $\infseqt{Y_t}$ is a sequence of independent random variables with individual means and variances given by $\mu_t := \EE(Y_t)$ and $\sigma_t^2 := \Var(Y_t)$, respectively. Suppose that in addition to \cref{assumption:lindeberg_variance_does_not_vanish} and the Lyapunov-type condition $\sum_{i=1}^\infty \left [ \EE|Y_i - \mu_i|^{2+\delta}/\sqrt{V_i}^{2+\delta} \right ] < \infty$, we have the following regularity conditions: 
  \begin{equation}\label{eq:lyapunov-asympcs-conditions}
    \sum_{i=1}^\infty \frac{\EE |Y_i^2 - \EE Y_i^2|^{1+\beta}}{V_i^{1+\beta}} < \infty, \quad\widetilde \mu_t^2 = o(V_t)~~\mathrm{a.s.},\quad \text{and}\quad\frac{1}{t}\sum_{i=1}^t (\mu_i - \widetilde \mu_t)^2 = o(\widetilde \sigma_t^2)
  \end{equation}
  for some $\beta \in (0, 1)$.
  In other words, the higher moments of $\infseqt{Y_t}$, the running mean $\widetilde \mu_t$, and the cumulative ``variation in means'' $\sum_{i=1}^t (\mu_i - \widetilde \mu_t)^2$ all cannot diverge too quickly relative to $\infseqt{V_t}$.
Then using the sample variance for $\widehat \sigma_t^2$, $\widetilde C_t$ forms a $(1-\alpha)$-\asympcs{} for the running average mean $\widetilde \mu_t := \frac{1}{t}\sum_{i=1}^t \mu_i$.
\end{corollary}
Clearly, the conditions of \eqref{eq:lyapunov-asympcs-conditions} are trivially satisfied if the means and the $(2+2\beta)^\text{th}$ absolute central moments are uniformly bounded over time. Since \namecref{assumption:lindeberg_variance_does_not_vanish}s~\ref{assumption:lindeberg_variance_does_not_vanish} and \ref{assumption:lindeberg-type_condition} hold by the assumptions of \cref{corollary:lyapunov-asympcs}, the proof in \cref{proof:lyapunov-asympcs} simply shows how the conditions in \eqref{eq:lyapunov-asympcs-conditions} imply \cref{assumption:lindeberg_variance_estimator_consistent}.

As suggested by \cref{section:general-framework}, we can combine \cref{theorem:general-asympcs} with essentially any other Gaussian boundary, and indeed there are others that can yield Lindeberg- and Lyapunov-type \asympcs{}s but we do not enumerate any more here, though we do mention one inspired by \citet[Eq. (20)]{robbins1970statistical} in passing in \cref{section:delayed-start}. The next section discusses how all of the aforementioned \asympcs{}s satisfy a certain formal asymptotic coverage guarantee.

\subsection{Asymptotic coverage and type-I error control}\label{section:type-I-error}
While the \asympcs{}s derived thus far serve as sequential analogues of CLT-based \ci{}s, it is not immediately obvious whether the bounds introduced in the previous section enjoy similar asymptotic coverage (equivalently, type-I error) guarantees. We will now give a positive answer to this question by showing that after appropriate tuning, our \asympcs{}s have asymptotic $(1-\alpha)$-coverage uniformly for all $t \geq m$ as $m \to \infty$ (to be formalized in \cref{definition:asymptotic-coverage}).

Recall that the coverage of CLT-based \ci{}s is at least $(1-\alpha)$ in the limit:
\begin{equation}\label{eq:asymptotic-fixed-time-coverage}
  \liminf_{n \to \infty}\PP(\mu \in \dot C_n) \geq 1-\alpha,
\end{equation}
but what is the right time-uniform analogue of \eqref{eq:asymptotic-fixed-time-coverage}? Since any single \asympcs{} will simply have \emph{some} coverage, we provide the following definition as a time-uniform analogue of \eqref{eq:asymptotic-fixed-time-coverage} for \emph{sequences} of \asympcs{}s that start later and later.

\begin{definition}[Asymptotic time-uniform coverage]\label{definition:asymptotic-coverage}
  For each $m \in \NN$, let $(C_t(m))_{t=m}^\infty$ be a sequence of sets, and let $\alpha \in (0, 1)$ be the desired miscoverage level. We say that $(C_t(m))_{t=m}^\infty$ has \emph{asymptotic time-uniform $(1-\alpha)$-coverage} for $(\mu_t)_{t=1}^\infty$ if
  \begin{equation}\label{eq:asymptotic-time-uniform-coverage}
    \liminf_{m \to \infty}\PP(\forall t \geq m ,\ \mu_t \in C_t(m) ) \geq 1-\alpha,
  \end{equation}
  and we say that this coverage is \emph{sharp} if the above inequality holds with equality.
\end{definition}
To the best of our knowledge, the existing literature lacks a concrete definition of asymptotic time-uniform coverage (or type-I error control) like \cref{definition:asymptotic-coverage}, but sequences of \asympcs{}s satisfying \eqref{eq:asymptotic-time-uniform-coverage} have been implicit in \citet{robbins1970statistical} and \citet{robbins1970boundary}, and the followup work of \citet{bibaut2022near}. 
In what follows, we provide (sharp) coverage guarantees for our \asympcs{}s. Furthermore, in \cref{section:delayed-start} we strictly improve on aforementioned bounds by \citet{robbins1970statistical} and \citet{robbins1970boundary}. Furthermore, we note that a bound~\citet{bibaut2022near} is in a certain sense equivalent to one that we provide here.

In order to obtain asymptotic time-uniform coverage, we  need a stronger variant of \cref{assumption:lindeberg_variance_estimator_consistent} so that variances are estimated at polynomial rates (rather than at arbitrary rates).
\begin{conditionalt}{\ref*{assumption:lindeberg_variance_estimator_consistent}-$\eta$}[Polynomial rate variance estimation]
  \label{assumption:lindeberg-variance-polynomial-rate}
  There exists some $0<\eta<1$ such that
  \begin{equation}
    \label{eq:lindeberg-variance-polynomial-rate}
  \widehat \sigma_t^2 - \widetilde \sigma_t^2 = o\left ( \frac{(t\widetilde \sigma_t^2)^\eta}{t} \right )\quad\text{almost surely}.
  \end{equation}
\end{conditionalt}
Note that while \cref{assumption:lindeberg-variance-polynomial-rate} is stronger than \cref{assumption:lindeberg_variance_estimator_consistent}, it is still quite mild. For instance, if $\widetilde \sigma_t^2$ is uniformly bounded, then \eqref{eq:lindeberg-variance-polynomial-rate} simply requires that $\widehat \sigma_t^2 - \widetilde \sigma_t^2 = o(t^{\eta-1})$ (i.e.~strong consistency at \emph{any} polynomial rate, potentially much slower than $t^{-1/2}$). Moreover, in the \iid{} case with at least $(2+\delta)$ finite absolute moments, \cref{assumption:lindeberg-variance-polynomial-rate} always holds by the SLLNs of \citet{marcinkiewicz1937fonctions}.

Our goal now is to show that sequences of \asympcs{}s given in \cref{theorem:lindeberg-martingale-asympcs} have asymptotic time-uniform coverage, and we will achieve this by effectively tuning them for later and later start times. Recall that \cref{section:optimizingMixture} allows us to choose the parameter $\rho > 0$ so that the \asympcs{} is tightest at some particular time --- we will now choose $\rho_m$ based on the first peeking time $m$ as
  \begin{equation}\label{eq:rho-m}
   \rho_m := \rho(\widehat \sigma_m^2m \log (m \lor e)) \equiv \sqrt{\frac{-2 \log\alpha + \log(-2\log \alpha) + 1}{\widehat \sigma_m^2 m \log(m \lor e) }}.
 \end{equation}
Then, let $\infseqtm{\widetilde C_t(m)}$ be the Gaussian mixture \asympcs{} with $\rho_m$ plugged into the expression of the boundary for all $t \geq m$:
  \begin{equation}\label{eq:asympcs-after-m}
   \widetilde C_t(m) := \left (\widehat \mu_t \pm \sqrt{\frac{2(t \widehat \sigma_t^2 \rho_m^2 + 1)}{t^2 \rho_m^2} \log \left ( \frac{\sqrt{t\widehat \sigma_t^2 \rho_m^2 + 1}}{\alpha} \right )} \right ).
  \end{equation}
  In other words, $\infseqtm{\widetilde C_t(m)}$ should be thought of as an \asympcs{} that only starts after time $m$, and is vacuous beforehand. The following theorem formalizes the coverage guarantees satisfied by this \emph{sequence} of \asympcs{}s as $m \to \infty$.

\begin{theorem}[Asymptotic $(1-\alpha)$-coverage for Gaussian mixture \asympcs{}s]\label{theorem:type-I-error}
  Given the same setup as \cref{theorem:lindeberg-martingale-asympcs} and \namecref{assumption:lindeberg_variance_does_not_vanish}s~\ref{assumption:lindeberg_variance_does_not_vanish},~\ref{assumption:lindeberg-type_condition}, and ~\ref{assumption:lindeberg-variance-polynomial-rate}, the \asympcs{}s $\infseqtm{\widetilde C_t(m)}$ given in \eqref{eq:asympcs-after-m} have sharp asymptotic $(1-\alpha)$-coverage for $\widetilde \mu_t := \frac{1}{t}\sum_{i=1}^t \mu_i$ as $m \to \infty$, meaning
  \begin{equation}\label{eq:type-I-error}
    \lim_{m \to \infty} \PP \left ( \forall t \geq m ,\ \widetilde \mu_t \in \widetilde C_t(m) \right ) = 1- \alpha.
  \end{equation}
\end{theorem}
The proof can be found in \cref{proof:type-I-error}. Notice that in the \iid{} setting, as $m \to \infty$, $\infseqtm{\widetilde C_t(m)}$ is asymptotically equivalent to $\infseqtm{\widebar C_t^\Gcal(m)}$ given by
\begin{equation}\label{eq:type-I-error-iid}
  \widebar C_t^\Gcal(m) :=  \left ( \widehat \mu_t \pm \widehat \sigma_t\sqrt{\frac{2(t\widebar \rho_m^2 + 1)}{t^2\widebar \rho_m^2} \log \left ( \frac{\sqrt{t \widebar \rho_m^2 + 1}}{\alpha} \right )} \right ),
\end{equation}
where $\widebar \rho_m := \rho(m \log (m\lor e))$. A quick inspection of the proof will reveal that \eqref{eq:type-I-error-iid} also satisfies the coverage guarantee provided in \cref{theorem:type-I-error} under the condition that $\widetilde \sigma_t^2 \to \sigma_\star^2 > 0$ almost surely. In summary, \eqref{eq:type-I-error-iid} can be thought of as an analogue of \eqref{eq:asympcs-after-m} for the \asympcs{}s that were derived in \cref{theorem:acs}.


\subsection{Asymptotic confidence sequences using Robbins' delayed start}\label{section:delayed-start}
As is clear from \cref{theorem:general-asympcs}, virtually any boundary for Gaussian observations can be used to derive an \asympcs{} as long as an appropriate strong invariance principle can be applied under the given assumptions --- indeed, \cref{theorem:acs}, \cref{proposition:acsLIL}, \cref{theorem:lindeberg-martingale-asympcs}, and \cref{corollary:lyapunov-asympcs} are all instantiations of the general phenomenon outlined in \cref{theorem:general-asympcs}.

Another \asympcs{} that may be of interest to practitioners is one that leverages Robbins' \cs{} for means of Gaussian random variables with a delayed start time \citep[Eq.~(20)]{robbins1970statistical}.
 In a nutshell, Robbins calculated a lower bound on the probability that a centered Gaussian random walk would remain within a particular two-sided boundary for all times $t \geq m$ given some starting time $m \geq 1$. That is, he showed that for \iid{} Gaussians $\infseqt{Z_t}$ with mean zero and variance $\sigma^2$, letting $G_t := \sum_{i=1}^t Z_i / \sigma$, and for any $a > 0$,
\begin{equation}\label{eq:robbins-delayed-start}
  \PP \left ( \forall t \geq m : \left \lvert G_t \right \rvert  < \sqrt{t(a^2 + \log (t / m))} \right ) \geq 1 -\Psi(a) ,
\end{equation}
where $\Psi(a) := 2(1-\Phi(a) + a\phi(a))$ and $\Phi$ and $\phi$ are the CDF and PDF of a standard Gaussian, respectively. In particular, setting $a \in (0, \infty)$ so that $\Psi(a) = \alpha$ yields a two-sided $(1-\alpha)$-boundary for the Gaussian random walk $\infseqt{G_t}$, and indeed, a solution to $\Psi(a) = \alpha$ always exists and is trivial to compute due to the fact that $\Psi$ is strictly decreasing, starting at $\Psi(0)=1$ and $\lim_{a \to \infty}\Psi(a) = 0$.
In a followup paper, \citet{robbins1970boundary} extended the ideas of \citet{robbins1970statistical} to a large class of boundaries for Wiener processes so that the probabilistic inequality in \eqref{eq:robbins-delayed-start} can be shown to be an equality when $|G_t|$ is replaced by the absolute value of a Wiener process (which would imply the inequality in \eqref{eq:robbins-delayed-start} for \iid{} standard Gaussians as a corollary).
Using this fact within the general framework of \cref{theorem:general-asympcs} combined with the strong invariance principle of~\citet{strassen1967and} and the techniques found in the proof of \cref{theorem:type-I-error} yields the following result.
\begin{proposition}[Delayed-start \asympcs{}]\label{proposition:asympcs-delayed-start}
  Consider the same setup as \cref{theorem:type-I-error} so that $\infseqt{Y_t}$ have conditional means and variances given by $\mu_t := \EE(Y_t \mid Y_1^{t-1})$ and $\sigma_t^2 := \Var(Y_t \mid Y_1^{t-1})$. Then under \namecref{assumption:lindeberg_variance_does_not_vanish}s~\ref{assumption:lindeberg_variance_does_not_vanish}, \ref{assumption:lindeberg-type_condition}, and \ref{assumption:lindeberg_variance_estimator_consistent}, we have that for any $m \geq 1$,
  \begin{equation}\label{eq:asympcs-delayed-start}
    \widetilde C_t^\mathrm{DS}(m) := 
                    \left ( \widehat \mu_t \pm \widehat \sigma_t \sqrt{ \frac{1}{t}\left [a^2 + \log\left ( \frac{t \widehat \sigma_t^2}{m \widehat \sigma_m^2}\right )\right ]} \right ) \quad \text{if}~t \geq m
   \end{equation}
   (and all of $\RR$ otherwise) forms a $(1-\alpha)$-\asympcs{} for $\widetilde \mu_t$, where $a$ is chosen so that $\Psi(a) = \alpha$. Furthermore, under~\cref{assumption:lindeberg-variance-polynomial-rate}, $\widetilde C_t^\mathrm{DS}(m)$ has sharp asymptotic time-uniform $(1-\alpha)$-coverage in the sense of \cref{definition:asymptotic-coverage}.
 \end{proposition}
The proof is provided in \cref{proof:asympcs-delayed-start}.
Similar to the relationship between \cref{theorem:type-I-error} and \eqref{eq:type-I-error-iid},
we have that if variances happen to converge $\widetilde \sigma_t^2 \to \sigma_\star^2 > 0$ almost surely, then as a corollary of \cref{proposition:asympcs-delayed-start}, the sequence $\infseqt{\widetilde C_t^{\mathrm{DS}\star}(m)}$ given by
  \begin{equation}\label{eq:delayed-start-corollary}
    \widetilde C_t^{\mathrm{DS}\star}(m) := \left ( \widehat \mu_t \pm \widehat \sigma_t \sqrt{ \frac{1}{t} \left ( a^2 + \log \left ( \frac{t}{m} \right ) \right )} \right ) \quad \text{for $t \geq m$}
  \end{equation}
  has asymptotic time-uniform coverage as in \cref{definition:asymptotic-coverage}. This can be seen as a generalization and improvement of results implied by \citet{robbins1970statistical} and \citet{robbins1970boundary}, and has connections to \citet{bibaut2022near}. We elaborate on these below.

\subsubsection{Relationship to Robbins and Siegmund}\label{remark:relationship-to-robbins-siegmund}
  Informally, \citet[Theorem 2(i)]{robbins1970boundary} show that for independent and identically distributed random variables $\infseqt{Y_t} \sim \PP$ with mean zero and unit variance (without loss of generality), the probability of their scaled partial sums $S_t := \sum_{i=1}^t Y_i$ exceeding a particular boundary behaves like a rescaled Wiener process exceeding that boundary. Consequently, for \iid{} data with known variance, their result combined with Robbins' delayed start \citep[Eq.~(20)]{robbins1970statistical} implies an asymptotic coverage guarantee (in the sense of \cref{definition:asymptotic-coverage}) for $\widetilde C_t^{\DS\star}$ given in \eqref{eq:delayed-start-corollary} but with $\widehat \sigma_t$ replaced by the true $\sigma$. As such, \eqref{eq:delayed-start-corollary} should be thought of as a generalization of their boundary when variances are unknown under martingale dependence. Nevertheless, \cref{proposition:asympcs-delayed-start} is strictly more general, allowing for variances to never converge.

\subsubsection{Relationship to \citet{bibaut2022near}}\label{remark:relationship-to-bibaut}
  In version 1 of \citet{bibaut2022near}, the authors derive a particular \asympcs{} and show that sequences thereof satisfy an asymptotic coverage guarantee in the same sense as \cref{definition:asymptotic-coverage}. That bound resembles --- but is always looser than --- a corollary of \cref{proposition:asympcs-delayed-start} in \eqref{eq:delayed-start-corollary} for a fixed $m \geq 1$; see \cref{claim:bv1}. Nevertheless, it was their asymptotic type-I error results that inspired us to show that the same guarantees hold for the bounds in \eqref{eq:asympcs-after-m}, \eqref{eq:asympcs-delayed-start}, and \eqref{eq:delayed-start-corollary}, and without the assumption that $\widetilde \sigma_t^2$ converges almost surely for the former two. After our advances, the second version of their paper also weakened their almost-sure convergent variance assumption and introduced a bound called the ``running maximum likelihood SPRT'' (rmlSPRT) which is identical to the corollary of \cref{proposition:asympcs-delayed-start} found in \eqref{eq:delayed-start-corollary} (modulo subtle differences in variance estimation techniques); see \cref{claim:bv2}. Since the exact connections may not be obvious to the reader, we derive them explicitly in \cref{section:explicit-connections-bibaut}. Despite their rmlSPRT being identical to \eqref{eq:delayed-start-corollary}, we remark that their paper focuses on testing guarantees and contains several additional interesting investigations including a sophisticated analysis of the expected rejection time, enriching the landscape of asymptotic anytime-valid methodology.


\section{Illustration:~Causal effects and semiparametric estimation}\label{section:csate}

Given the groundwork laid in \cref{section:time-uniform-clt}, we now demonstrate the use of \asympcs{}s for conducting anytime-valid causal inference. Since it is an important and thoroughly studied functional, we place a particular emphasis on the average treatment effect (ATE) for illustrative purposes but we discuss how these techniques apply to semiparametric functional estimation more generally in \cref{section:generalFunctionals}. The literature on semiparametric functional inference often falls within the asymptotic regime and hence \asympcs{}s form a natural time-uniform extension thereof.

It is important to note that obtaining \asympcs{}s for the ATE is not as simple as directly applying the theorems of \cref{section:asymptoticCS} to some appropriately chosen augmented inverse-probability-weighted (AIPW) influence functions (otherwise the illustration of this section would have been trivial). Indeed, satisfying the conditions of the aforementioned theorems --- \cref{theorem:general-asympcs} in particular --- in the presence of infinite-dimensional nuisance parameters is nontrivial and the analysis proceeds rather differently from the fixed-$n$ setting. Nevertheless, after introducing and carefully analyzing sequential sample splitting and cross-fitting (\cref{section:seqSampleSplit}), we will see that asymptotic time-uniform inference for the ATE is possible.

To solidify the notation and problem setup, suppose that we observe a (potentially infinite) sequence of \iid{} variables $Z_1, Z_2, \dots$ from a distribution $\PP$
where $Z_t := (X_t, A_t, Y_t)$ denotes the $t^\text{th}$ subject's triplet and $X_t \in \mathbb R^d$ are their measured baseline covariates, $A_t \in \{0, 1\}$ is the treatment that they receive, and $Y_t \in \mathbb R$ is their measured outcome after treatment.
Our target estimand is the average treatment effect (ATE) $\psi$ defined as
\begin{equation}
   \psi := \EE(Y^1 - Y^0),
\end{equation}
where $Y^a$ is the counterfactual outcome for a randomly selected subject had they received treatment $a \in \{0, 1\}$. The ATE $\psi$ can be interpreted as the average population outcome if everyone were treated $\EE(Y^1)$ versus if no one were treated $\EE(Y^0)$. Under standard causal identification assumptions --- typically referred to as consistency, positivity, and exchangeability (see e.g.~\citet[\S 2.2]{kennedy2016semiparametric}) --- we have that $\psi$ can be written as a (non-counterfactual) functional of the distribution $\PP$:
\begin{equation}
 \psi \equiv \psi(\PP) = \EE\{ \EE(Y \mid X, A = 1) - \EE(Y \mid X, A= 0) \}.
\end{equation}
Throughout the remainder of this section, we will operate under these identification assumptions and aim to derive efficient \asympcs{}s for $\psi$ using tools from semiparametric theory.
At a high level, we will construct \asympcs{}s for $\psi$ by combining the results of \cref{section:time-uniform-clt} with sample averages of \emph{influence functions} for $\psi$ and in the ideal case, these influence functions will be \emph{efficient} (in the semiparametric sense).

\subsection{Sequential sample splitting and cross fitting}\label{section:seqSampleSplit}
Following \citet{robins2008higher}, \citet{zheng2010asymptotic}, and \citet{chernozhukov2017double}, we employ sample splitting to derive an estimate $\widehat f$ of the influence function $f$ on a ``training'' sample, and evaluate $\widehat f$ on values of $Z_t$ in an independent ``evaluation'' sample.
Sample splitting sidesteps complications introduced from ``double-dipping'' (i.e. using $Z_t$ to both construct $\widehat f$ and evaluate $\widehat f(Z_t)$) and greatly simplifies the analysis of the downstream estimator.
Since the aforementioned authors employed sample splitting in the \emph{batch} (non-sequential) regime while we are concerned with settings where data are continually observed in an online stream over time, we modify the sample splitting procedure as follows. We will denote $\splitSetInfty^\train$ and $\splitSetInfty^\eval$ as the ``training'' and ``evaluation'' sets, respectively. At time $t$, we assign $Z_t$ to either group with equal probability:
\[ Z_t \in
\begin{cases}
\splitSetInfty^\train & \text{with probability } 1/2,\\
\splitSetInfty^\eval & \text{otherwise.}\end{cases}\]
Note that at time $t+1$, $Z_t$ is \emph{not} re-randomized into either split --- once $Z_t$ is randomly assigned to one of $\splitSetInfty^\train$ or $\splitSetInfty^\eval$, they remain in that split for the remainder of the study. In this way, we can write $\splitSetInfty^\train = (Z_1^\train, Z_2^\train, \dots)$ and $\splitSetInfty^\eval = (Z_1^\eval, Z_2^\eval, \dots )$ and think of these as independent, sequential observations from a common distribution $\PP$. To keep track of how many subjects have been randomized to $\splitSetInfty^\train$ and $\splitSetInfty^\eval$ at time $t$, define
\begin{equation}
\label{eq:tau_and_tauprime}
\Teval := | \splitSetInfty^\eval | ~~~\text{and}~~~ \Ttrain := | \splitSetInfty^\train  | \equiv t - \Teval,
\end{equation}
where we have left the dependence on $t$ implicit.
\begin{remark}
    Strictly speaking, under the \iid{} assumption, we do not need to randomly assign subjects to training and evaluation groups for the forthcoming results to hold
    (e.g. we could simply assign even-numbered subjects to $\splitSetInfty^\train$ and odd-numbered subjects to $\splitSetInfty^\eval$).
    However, the analysis is not further complicated by this randomization, and it can be used to combat bias in treatment assignments when the \iid{} assumption is violated \citep{efron1971forcing}.
\end{remark}
\subsubsection{The sequential sample-split estimators $(\widehat \psi_t^\samplesplit)_{t=1}^\infty$} After employing sequential sample splitting, the sequence of sample-split estimators $(\widehat \psi_{t}^\samplesplit)_{t=1}^\infty$ for $\psi$ are given by
\begin{equation}
\label{eq:dr_estimator}
     \widehat \psi_t^\samplesplit := \frac{1}\Teval \sum_{i=1}^\Teval \widehat f_{\Ttrain}(Z_i^\eval),
\end{equation}
where $\widehat f_{\Ttrain}$ is given by the so-called \emph{efficient influence function} (a brief review of semiparametric efficient estimators can be found in \cref{section:efficient-estimators-review}),
\begin{equation}\label{eq:efficientInfluenceFunction}
    f(z) \equiv f(x, a, y) := \left \{ \mu^1(x) - \mu^0(x) \right \} + \left ( \frac{a}{\pi(x)} - \frac{1-a}{1-\pi(x)} \right )\left \{ y - \mu^a(x) \right \},
\end{equation}
with $\eta \equiv (\mu^1, \mu^0, \pi)$ replaced by $\widehat \eta_{\Ttrain} \equiv (\widehat \mu^1_{\Ttrain}, \widehat \mu^0_{\Ttrain}, \widebar \pi_{\Ttrain}) $ --- where $\widebar \pi_{\Ttrain}$ may be an estimator $\widehat \pi_{\Ttrain}$ of the propensity score $\pi$, or the propensity score itself, depending on whether one is considering an observational study or randomized experiment --- so that $\widehat \eta_\Ttrain$ is built solely from $\splitSetInfty^\train$. The sample splitting procedure for constructing $\widehat \psi_t^\samplesplit$ is summarized pictorially in \cref{fig:seqSampleSplitting}. In the batch setting for a fixed sample size,~\eqref{eq:dr_estimator} is often referred to as the \emph{augmented inverse probability weighted} (AIPW) estimator \citep{rotnitzky1998semiparametric, scharfstein1999adjusting} (an instantiation of so-called ``one-step correction'' in the semiparametrics literature) and we adopt similar nomenclature here.
\begin{figure}[!htbp]
  \centering
\begin{tikzpicture}[>=latex]
  \node[font=\fontsize{12}{14}\selectfont] (varZ) {$Z_t$};
  \node[draw, above right=-0.295cm and 2cm of varZ, minimum width=1.5cm, minimum height=2cm, font=\fontsize{12}{14}\selectfont] (box1) {$\splitSetInfty^\train$};
  \node[draw, below right=-0.295cm and 2cm of varZ, minimum width=1.5cm, minimum height=2cm, font=\fontsize{12}{14}\selectfont] (box2) {$\splitSetInfty^\eval$};
  \node[right of=box1, node distance=3cm, font=\fontsize{12}{14}\selectfont] (varY) {$\left (\widehat \mu_{\Ttrain}^1, \widehat \mu_\Ttrain^0, \widehat \pi_\Ttrain \right )$};
  \node[right of=box2, node distance=3cm, font=\fontsize{12}{14}\selectfont] (varZeval) {$\{ Z_1^\eval, \dots, Z_T^\eval \}$};
  \node[right of=varZ, node distance=9cm, font=\fontsize{12}{14}\selectfont] (varW) {$\frac{1}{T}\sum_{i=1}^\Teval \widehat f_\Ttrain(Z_i^\eval) =: \widehat \psi_t^\samplesplit $ };

  \draw[->] (varZ) -- node[above, rotate=25, anchor=south] {w.p. 1/2} (box1);
  \draw[->] (varZ) -- node[above, rotate=-25, anchor=north] {w.p. 1/2} (box2);

  \draw[->] (box1) -- (varY);
  \draw[->] (box2) -- (varZeval);

  \draw[->] ($(varY.east)!0.5!(varY.west)+(1.3cm, -0.2cm)$) -- ($(varW.west)!0.3!(varW.east)+(0,0.3cm)$);
  \draw[->] ($(varZeval.east)!0.5!(varZeval.west)+(1.6cm, 0cm)$) -- ($(varW.west)!0.3!(varW.east)+(0.9,-0.3cm)$);
\end{tikzpicture}
  \caption{A schematic illustrating sequential sample splitting. At each time step $t$, the new observation $Z_t$ is randomly assigned to $\splitSetInfty^\train$ or $\splitSetInfty^\eval$ with equal probability (1/2). Nuisance function estimators $(\widehat \mu_{\Ttrain}^1, \widehat \mu_{\Ttrain}^0, \widehat \pi_{\Ttrain})$ are constructed using $\splitSetInfty^\train$ which then yield $\widehat f_{\Ttrain}$. The sample-split estimator $\widehat \psi_t^\samplesplit$ is defined as the sample average $\frac{1}{\Teval} \sum_{i=1}^\Teval \widehat f_{\Ttrain}(Z_i^\eval)$ where each $Z_i^\eval \in \splitSetInfty^\eval$.}
    \label{fig:seqSampleSplitting}
\end{figure}
\subsubsection{The sequential cross-fit estimators $(\widehat \psi_t^\times)_{t=1}^\infty$}
A commonly cited downside of sample splitting is the loss in efficiency by using $\Teval \approx t/2$ subjects instead of $t$ when evaluating the sample mean $\frac{1}\Teval \sum_{i=1}^\Teval \widehat f_{\Ttrain} (Z_i^\eval)$. An easy fix is to \emph{cross-fit}: swap the two samples, using the evaluation set $\splitSetInfty^\eval$ for training and the training set $\splitSetInfty^\train$ for evaluation to recover the full sample size of $t \equiv \Teval + \Ttrain$ \citep{robins2008higher,zheng2010asymptotic,chernozhukov2017double}. That is, construct $\widehat f_\Teval$ solely from $\splitSetInfty^\eval$ and define the cross-fit estimator $\widehat \psi_t^\times$ as
\begin{equation}
\label{eq:cross_fit_estimator}
     \widehat \psi_t^\times := \frac{\sum_{i=1}^\Teval \widehat f_{\Ttrain}(Z_i^\eval) + \sum_{i=1}^{\Ttrain} \widehat f_{\Teval}(Z_i^\train)}{t},
\end{equation}
and the associated cross-fit variance estimate
\begin{equation}
  \label{eq:cross_fit_variance}
  \widehat \Var_t(\widehat f) := \frac{\widehat \Var_\Teval(\widehat f_{\Ttrain}) + \widehat \Var_{\Ttrain}(\widehat f_{\Teval})}{2},
\end{equation}
where $\widehat \Var_\Teval(\widehat f_\Ttrain)$ is the $\splitSetInfty^\eval$-sample variance of the pseudo-outcomes $(\widehat f_{\Ttrain}(Z_i^{\eval}))_{i=1}^{\Teval}$ and similarly for $\widehat \Var_\Ttrain$ (we deliberately omit the subscript on $\widehat f$ in the left-hand side of \eqref{eq:cross_fit_variance}). All of the results that follow are stated in terms of the cross-fit estimators $(\widehat \psi_t^\times)_{t=1}^\infty$ but they also hold using $(\widehat \psi_t^\samplesplit)_{t=1}^\infty$ instead.
With the setup of \cref{section:efficient-estimators-review} and \cref{section:seqSampleSplit} in mind, we are ready to derive \asympcs{}s for $\psi$, first in randomized experiments.

\subsection{Asymptotic confidence sequences in randomized experiments}
\label{section:randomized}
Consider a sequential randomized experiment so that a subject with covariates $x$ has a known propensity score  
\begin{equation}
    \pi(x) := \PP(A = 1 \mid X = x).
\end{equation}
Consider the cross-fit AIPW estimator $\widehat \psi_t^\times$ as given in \eqref{eq:cross_fit_estimator} but with estimated propensity scores --- $\widehat \pi_{\Ttrain}(x)$ and $\widehat \pi_{\Teval}(x)$ --- replaced by their true values $\pi(x)$, and with $\widehat \mu^a_{\Ttrain}$ and $\widehat \mu^a_\Teval$ being possibly misspecified estimators for $\mu^a$. We will assume that $\widehat \mu^a_t$ converges to some function $\widebar \mu^a$, which need not coincide with $\mu^a$. In what follows, when we use $\widehat \mu_t^a$ or $\widehat f_t$ in writing $\|\widehat \mu_t^a - \widebar \mu^a\|_\LP$ or $\|\widehat f_t - \widebar f\|_\LP$, we are referring to large-sample properties of the estimator (and hence $\widehat f_t$ could be replaced by $\widehat f_\Ttrain$ or $\widehat f_\Teval$ without loss of generality).
\begin{theorem}[Confidence sequences for the ATE in randomized experiments]
\label{theorem:csate_randomized}
    Let $\widehat \psi_t^\times$ be the cross-fit AIPW estimator as in \eqref{eq:cross_fit_estimator}. Suppose $\| \widehat \mu^a_t(X) - \widebar \mu^a(X) \|_{L_2(\PP)} = o(1)$ for each $a \in \{0, 1\}$ where $\widebar \mu^a$ is some function (but need not be $\mu^a$), and hence $\|\widehat f_{t} - \widebar f\|_{L_2(\PP)} = o(1)$ for some influence function $\widebar f$. Suppose that propensity scores are bounded away from 0 and 1, i.e. $\pi(X) \in [\delta, 1-\delta]$ almost surely for some $\delta > 0$, and suppose that $\EE |\widebar f(Z)|^{2+\eps} < \infty$ for some $\eps > 0$. Then for any constant $\rho > 0$,
\begin{equation}
  \widehat \psi_t^\times \pm \sqrt{\widehat \Var_t(\widehat f)} \cdot\sqrt{\frac{2(t \rho^2 + 1)}{t^2 \rho^2} \log \left ( \frac{\sqrt{t\rho^2 + 1}}{\alpha} \right ) }
\end{equation}
forms a $(1-\alpha)$-\asympcs{} for $\psi$.
\end{theorem}
The proof in Appendix~\ref{proof:csate_randomized} combines an analysis of the almost-sure convergence of $(\widehat \psi_t^\times - \psi)$ with the \asympcs{} of Theorem~\ref{theorem:acs}.
Notice that since $\widehat \mu_{t}^a$ is consistent for a function $\widebar \mu^a$, we have that $\widehat f_{t}$ is converging to some influence function $\widebar f$ of the form
\begin{equation}
\label{eq:randomizedExperimentInfluenceFunction}
     \widebar f(z) \equiv \widebar f(x, a, y) := \left \{ \widebar \mu^1(x) - \widebar \mu^0(x) \right \} + \left ( \frac{a}{\pi(x)} - \frac{1-a}{1-\pi(x)} \right )\left \{ y - \widebar \mu^a(x) \right \}.
\end{equation}
In practice, however, one must choose $\widehat \mu_{t}^a$. As alluded to at the beginning of \cref{section:csate}, the best possible influence function is the EIF $f(z)$ defined in \eqref{eq:efficientInfluenceFunction}, and thus it is natural to attempt to construct $\widehat \mu_{t}^a$ so that $\| \widehat f_{t} - f \|_{L_2(\PP)} = o(1)$. The resulting confidence sequences would inherit such optimality properties, a point which we discuss further in \cref{section:unimprovability}.

\begin{figure}[!htbp]
    \centering
    \includegraphics[width=\columnwidth]{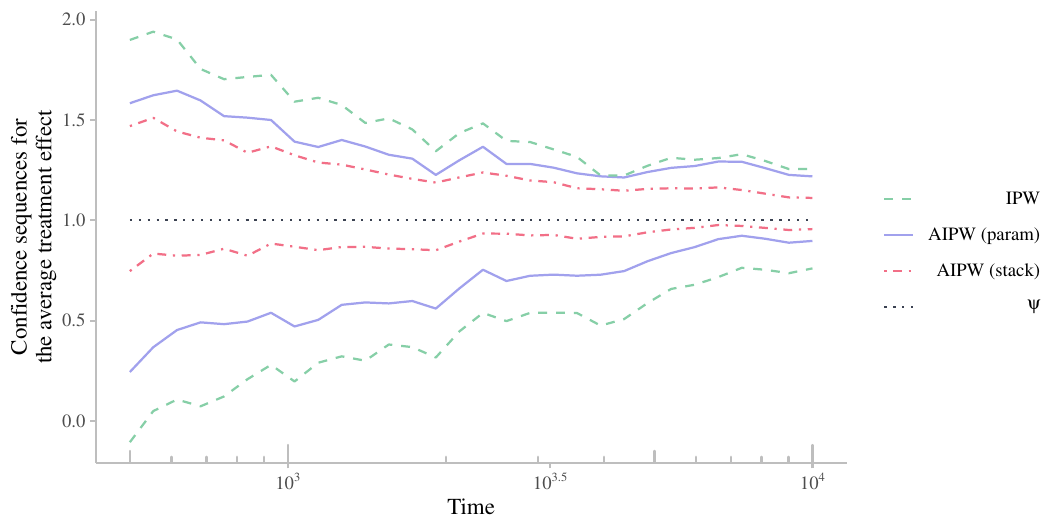}
    \caption{Three 90\%-\asympcs{}s for the average treatment effect in a simulated randomized experiment using different regression estimators. Notice that all three confidence sequences uniformly capture the average treatment effect $\psi$, but more sophisticated models do so more efficiently, with AIPW+stacking greatly outperforming IPW\@.}
    \label{fig:CS_randomized}
\end{figure}
Since $\mu^a$ is simply a conditional mean function, we can use virtually any regression techniques to estimate it. Here we will consider the general approach of \emph{stacking} introduced by \citet{breiman1996stacked} and further studied by \citet{tsybakov2003optimal} and \citet{van2007super} (see also \citep{polley2010super,polley2011super}) under the names of ``aggregation'' and ``Super Learning'' respectively. In short, stacking uses cross-validation to choose a weighted combination of $K$ candidate predictors where the weights are chosen based on data in held-out samples. Importantly (and under certain conditions), the stacked predictor will have a mean squared error that scales with that of the best of the $K$ predictors up to an additive $\log K$ term \citep{tsybakov2003optimal,van2006oracle}. This advantage can be seen empirically in \cref{fig:CS_randomized} where the true regression functions $\mu^0$ and $\mu^1$ are nonsmooth and nonlinear in $x$. Such advantages via stacking are not new --- we are only highlighting the observation that similar phenomena carry over to \asympcs{}s. 

So far, the use of flexible regression techniques like stacking were used only for the purposes of deriving sharper \asympcs{}s in sequential randomized experiments. In observational studies, however, consistent estimation of nuisance functions at fast rates is essential to the construction of \emph{valid} fixed-$n$ \ci{}s, and indeed the same is true for \asympcs{}s.

\subsection{Asymptotic confidence sequences in observational studies}
\label{section:observational}

Consider now a sequential observational study (e.g.~we are able to continuously monitor $\infseqt{X_t, A_t, Y_t}$ but do not know $\pi(x)$ exactly, or we are in a sequentially randomized experiment with noncompliance, etc.). The only difference in this setting with respect to setup is the fact that $\pi(x)$ is no longer known and must be estimated. As in the fixed-$n$ setting, this complicates estimation and inference.
The following theorem provides the conditions under which we can construct \asympcs{}s for $\psi$ using the cross-fit AIPW estimator in observational studies.
\begin{figure}[!htbp]
    \centering
    \includegraphics[width=\columnwidth]{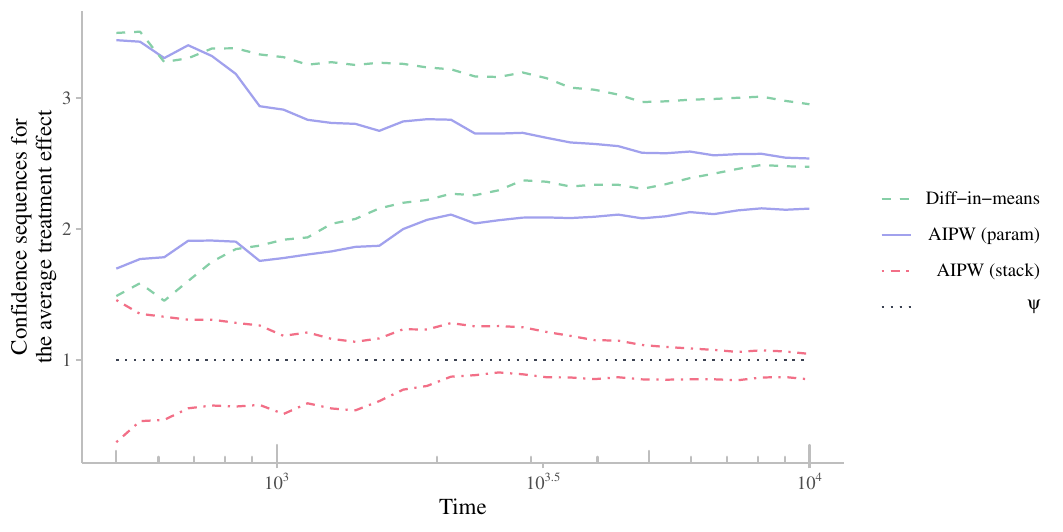}
    \caption{Three 90\%-\asympcs{}s for the ATE in an observational study using three different estimators --- a difference-in-means estimator, AIPW with parametric models, and AIPW with an ensemble of predictors combined via stacking. Unlike the randomized setup, only the stacking ensemble is consistent, since the other two are misspecified. Not only is the stacking-based \asympcs{} converging to $\psi$, but it is also the tightest of the three models at each time step.}
    \label{fig:CS_observational}
\end{figure}
\begin{theorem}[Confidence sequence for the ATE in observational studies]
  \label{theorem:csate_observational}
    Consider the same setup as Theorem~\ref{theorem:csate_randomized} but with $\pi(x)$ no longer known. Suppose that regression functions and propensity scores are consistently estimated in $L_2(\PP)$ at a product rate of $o(\sqrt{\log t / t})$, meaning that \[\| \widehat \pi_{t} - \pi \|_{L_2(\PP)} \sum_{a = 0}^1 \| \widehat \mu_{t}^a - \mu^a\|_{L_2(\PP)} = o \left ( \sqrt{\log t / t} \right ).\] Moreover, suppose that $\|\widehat f_{t} - f\|_{L_2(\PP)} = o(1)$ where $f$ is the efficient influence function \eqref{eq:efficientInfluenceFunction} and that $\EE |f(Z)|^{2+\eps} < \infty$ for some $\eps > 0$. Then for any constant $\rho > 0$,
    \[ \widehat \psi_t^\times \pm \sqrt{\widehat \Var_t(\widehat f)} \cdot \sqrt{\frac{2(t\rho^2 + 1)}{t^2 \rho^2} \log \left ( \frac{\sqrt{t\rho^2 + 1}}{\alpha} \right ) } \]
    forms a $(1-\alpha)$-\asympcs{} for $\psi$.
\end{theorem}
The proof in Appendix~\ref{proof:csate_observational} proceeds similarly to the proof of Theorem~\ref{theorem:csate_randomized} by combining Theorem~\ref{theorem:acs} with an analysis of the almost-sure behavior of $(\widehat \psi_t^\times - \psi)$.
Notice that the nuisance estimation rate of $\sqrt{\log t / t}$ is slower than $1/\sqrt{t}$ which is usually required in the fixed-$n$ regime, but we do require almost-sure convergence rather than convergence in probability.

Unlike the experimental setting of \cref{section:randomized}, Theorem~\ref{theorem:csate_observational} requires that $\widehat \mu_{t}^a$ and $\widehat \pi_{t}$ consistently estimate $\mu^a$ and $\pi$, respectively. As a consequence, the stacking-based AIPW \asympcs{} is both the tightest of the three \emph{and} is uniquely consistent for $\psi$ (see~\cref{fig:CS_observational}).

\subsection{The running average of individual treatment effects}
The results in Sections~\ref{section:randomized} and~\ref{section:observational} considered the classical regime where the ATE $\psi$ is a fixed functional that does not change over time. Consider a strict generalization where distributions --- and hence individual treatment effects in particular --- may change over time. In other words,
\begin{equation}
  \label{eq:ate_t}
  \psi_t := \EE\left \{ Y^1_t - Y^0_t \right \} \overset{(\star)}{=} \EE \left \{ \EE( Y_t \mid X_t, A_t = 1 ) - \EE(Y_t \mid X_t, A_t = 0) \right \},
\end{equation}
where the equality $(\star)$ holds under the familiar causal identification assumptions discussed earlier. Despite the non-stationary and non-\iid{} structure, it is nevertheless possible to derive \asympcs{}s for the \emph{running average} of individual treatment effects $\widetilde \psi_t := \frac{1}{t} \sum_{i=1}^t \psi_i$ --- or simply, the running average treatment effect --- using the Lyapunov-type bounds of \cref{corollary:lyapunov-asympcs}. However, given this more general and complex setup, the assumptions required are more subtle (but no more restrictive) than those for \namecref{theorem:csate_randomized}s \ref{theorem:csate_randomized} and \ref{theorem:csate_observational}; as such, we explicitly describe their details here but handle the randomized and observational settings simultaneously for brevity.

\begin{conditiontvate}[Regression estimator is uniformly well-behaved in $\LP$]
  \label{assumption:sup-regression-estimators-converge}
  We assume that regression estimators $\mu_t^a(X_i)$ converge in $\LP$ to any function $\widebar \mu^a(X_i)$ uniformly for $i \in \{1, 2, \dots\}$ i.e.
  \begin{equation}
    \sup_{1 \leq i < \infty }\| \widehat \mu_{t}^a(X_i) - \widebar \mu^a(X_i) \|_\LP = o(1)
  \end{equation}
for each $a \in \{0, 1\}$.
\end{conditiontvate}
\cref{assumption:sup-regression-estimators-converge} simply requires that the regression estimator $\widehat \mu_t^a$ must converge to some function $\widebar \mu^a$, which need not coincide with true regression function $\mu^a$. In the \iid{} setting where $X_1,X_2, \dots$ all have the same distribution, we would simply drop the $\sup_{1 \leq i \leq \infty}$, recovering the conditions for Theorems~\ref{theorem:csate_randomized} and~\ref{theorem:csate_observational}.


\begin{conditiontvate}[Convergence of average nuisance errors]
  \label{assumption:doubly-robust-average-bias-convergence}
  Let $\widehat \mu_t^a$ be an estimator of the regression function $\mu^a$, $a \in \{0, 1\}$ and $\widehat \pi_t$ an estimator of the propensity score $\pi$. We assume that the average bias shrinks at a $\sqrt{\log t / t}$ rate, i.e.
\begin{equation}
  \label{eq:doubly-robust-average-bias-convergence}
  \frac{1}{t}\sum_{i=1}^t \left \{ \|\widehat \pi_t(X_i) - \pi(X_i) \|_{L_2(\PP)} \sum_{a=0}^1  \| \widehat \mu_{t}^a(X_i) - \mu^a(X_i) \|_{L_2(\PP)} \right \} = o\left (\sqrt{\frac{\log t}{t}} \right).
\end{equation}
\end{conditiontvate}
Note that \cref{assumption:doubly-robust-average-bias-convergence} would hold in two familiar scenarios. Firstly, in a randomized experiment (Theorem~\ref{theorem:lyapunov_ate_randomized}) where $\widehat \pi_t = \pi$ is known by design, we have that $\eqref{eq:doubly-robust-average-bias-convergence}$ is always zero, satisfying \cref{assumption:doubly-robust-average-bias-convergence} trivially. Second, in an observational study where the product of errors $\|\widehat \pi_t(X_i) - \pi(X_i)\|_\LP \|\widehat \mu_t^a(X_i) - \mu^a(X_i)\|_\LP$ vanishes at a rate faster than $\sqrt{\log t / t}$, for each $i$ and for both $a \in \{0, 1\}$, we also have that their average product errors vanish at the same rate \eqref{eq:doubly-robust-average-bias-convergence}. With these assumptions in mind, let us summarize how running average treatment effects can be captured in randomized experiments.

\begin{theorem}[\asympcs{}s for the running average treatment effect]
  \label{theorem:lyapunov_ate_randomized}
  Suppose $Z_1, Z_2, \dots$ are independent triples $Z_t := (X_t, A_t, Y_t)$ and that \namecref{assumption:doubly-robust-average-bias-convergence}s~\ref{assumption:sup-regression-estimators-converge} and~\ref{assumption:doubly-robust-average-bias-convergence} hold. Finally, suppose that the conditions of \cref{corollary:lyapunov-asympcs} hold, but with $\infseqt{Y_t}$ replaced by the influence functions $\infseqt{\widebar f(Z_t)}$. Then,
  \begin{equation}
    \label{eq:lyapunov_ate_randomized_cs}
    \widehat \psi_t^\times \pm \sqrt{\frac{2(t \rho^2 \widehat \Var_t(\widebar f) + 1)}{t^2 \rho^2} \log \left ( \frac{\sqrt{t \rho^2\widehat \Var_t(\widebar f) + 1}}{\alpha} \right )}
 \end{equation}
 forms a $(1-\alpha)$-\asympcs{} for the running average treatment effect $\widetilde \psi_t := \frac{1}{t} \sum_{i=1}^t \psi_i$.
\end{theorem}
\begin{figure}[!htbp]
  \centering
  \includegraphics[width=0.9\columnwidth]{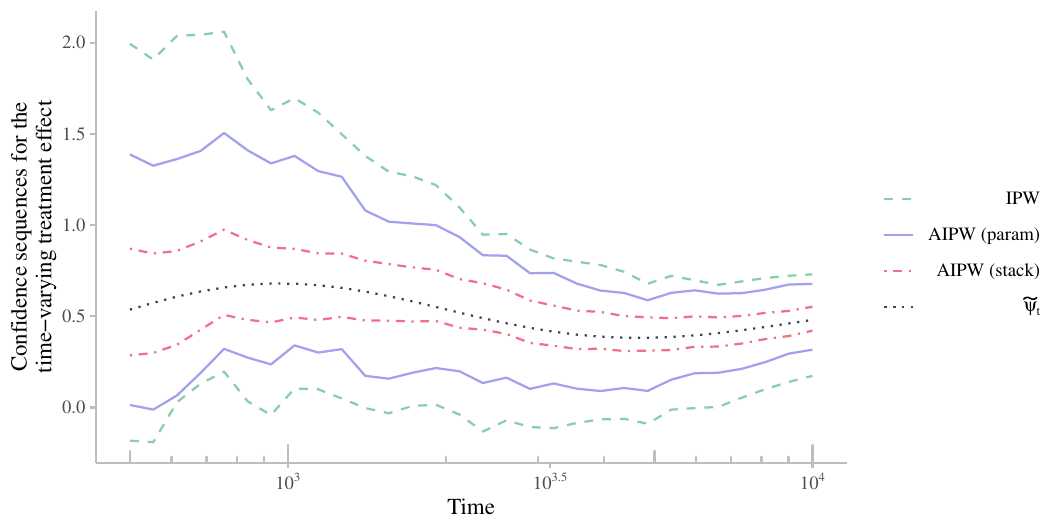}
  \caption{Three 90\% \asympcs{}s for $\widetilde \psi_t$ constructed using various estimators via Theorem~\ref{theorem:lyapunov_ate_randomized}. Since this is a randomized experiment, all three \cs{}s capture $\widetilde \psi_t$ uniformly over time with high probability. Similar to \cref{fig:CS_randomized}, however, the stacking-based AIPW estimator greatly outperforms those based on parametric models or IPW.}
  \label{fig:time-varying-ate-randomized}
\end{figure}
The proof can be found in \cref{proof:lyapunov_ate}, and it is not hard to see that both \namecref{theorem:csate_randomized}s~\ref{theorem:csate_randomized} and~\ref{theorem:csate_observational} are particular instantiations of \cref{theorem:lyapunov_ate_randomized}. The important takeaway from \cref{theorem:lyapunov_ate_randomized} is that under some rather mild conditions on the moments of $(\widebar f(Z_t))_{t=1}^\infty$, it is possible to derive an \asympcs{} for a running average treatment effect $\widetilde \psi_t$ (see \cref{fig:time-varying-ate-randomized} for what these look like in practice). Nevertheless, under the commonly considered regime where the treatment effect is constant $\psi_1 = \psi_2 = \cdots = \psi$, we have that \eqref{eq:lyapunov_ate_randomized_cs} forms a $(1-\alpha)$-\asympcs{} for $\psi$. Note that unlike \namecref{theorem:csate_randomized}s~\ref{theorem:csate_randomized} and~\ref{theorem:csate_observational}, \cref{theorem:lyapunov_ate_randomized} actually does require the use of the cross-fit AIPW estimator $\widehat \psi_t^\times$ and would not capture $\widetilde \psi_t$ if the sample-split version were used in its place.



\begin{remark}[Avoiding sample splitting via martingale \asympcs{}s]\label{remark:predictable_nuisance_estimation}
The reader may wonder whether it is possible to simply plug in a \emph{predictable} estimate of $\widehat \mu_t^a$ in a randomized experiment --- i.e.~so that $\widehat \mu_t^a$ only depends on $Z_1^{t-1}$ --- and employ the Lindeberg-type martingale \asympcs{} of \cref{theorem:lindeberg-martingale-asympcs} in place of \cref{corollary:lyapunov-asympcs}, thereby sidestepping the need for sequential sample splitting and cross fitting altogether. Indeed, such an analogue of \cref{theorem:lyapunov_ate_randomized} is possible to derive, but the conditions required are less transparent than those we have provided above so we defer it to \cref{section:martingale-time-varying-trt-effects}.
\end{remark}


\subsection{Extensions to general semiparametric estimation and the delta method}
\label{section:generalFunctionals}
    
The discussion thus far has been focused on deriving confidence sequences for the ATE in the context of causal inference. However, the tools presented in this paper are more generally applicable to any pathwise differentiable functional with positive (and finite) semiparametric information bound. Some prominent examples in causal inference include modified interventions, complier-average effects, time-varying effects, and controlled mediation effects, among others. Examples outside causal inference include the expected density, entropy, the expected conditional variance, and the expected conditional covariance, to list a few.

All of the aforementioned problems, including estimation of the ATE in \cref{section:csate} can be written in the following general form. Suppose $Z_1, Z_2, \dots \simiid \QQ$ and let $\theta(\QQ)$ be some functional (such as those listed above) of the distribution $\QQ$. In the case of a finite sample size $n$, $\widehat \theta_n$ is said to be an asymptotically linear estimator \citep{tsiatis2007semiparametric} for $\theta$ if
\begin{equation}
 \widehat \theta_n - \theta = \frac{1}{n} \sum_{i=1}^n \phi(Z_i) + o_{\QQ} \left (\frac{1}{\sqrt{n}} \right), 
\end{equation}
where $\phi$ is the influence function of $\widehat \theta_n$. When the sample size is not fixed in advance, we may analogously say that $\widehat \theta_t$ is an \emph{asymptotically linear time-uniform estimator} if instead,
\begin{equation}
\label{eq:asymp-linear-time-uniform-estimator}
\widehat \theta_t - \theta \eqas \frac{1}{t} \sum_{i=1}^t \phi(Z_i) + o\left ( \sqrt{\log t /t} \right ),
\end{equation}
with $\phi$ being the same influence function as before. For example, in the case of the ATE with $(Z_t)_{t=1}^\infty \simiid \PP$, we presented an efficient estimator $\widehat \psi_t$ which took the form,
\begin{equation}
  \widehat \psi_t - \psi \eqas \frac{1}{t} \sum_{i=1}^t (f(Z_i) - \psi) + o\left (\sqrt{\log t / t} \right ),
\end{equation}
where $f$ is the uncentered efficient influence function (EIF) defined in  \eqref{eq:efficientInfluenceFunction}. In order to justify that the remainder term is indeed $o( \sqrt{ \log t / t} )$, we used sequential sample splitting and additional analysis in the randomized and observational settings (see the proofs in Sections~\ref{proof:csate_randomized} and \ref{proof:csate_observational} for more details). 
In general, as long as an estimator $\widehat \theta_t$ for $\theta$ has the form \eqref{eq:asymp-linear-time-uniform-estimator}, we may derive \asympcs{}s for $\theta$ as a simple corollary of Theorem~\ref{theorem:acs}.
\begin{corollary}[\asympcs{}s for general functional estimation]\label{corollary:general-functionals}
    Suppose $\widehat \theta_t$ is an asymptotically linear time-uniform estimator of $\theta$ with influence function $\phi$, that is, satisfying \eqref{eq:asymp-linear-time-uniform-estimator}. Additionally, suppose that $\Var (\phi) < \infty$. Then,
\begin{equation}
  \widehat \theta_t ~\pm~ \sqrt{ \Var(\phi) } \cdot \sqrt{\frac{2(t \rho^2 + 1)}{t^2\rho^2} \log \left ( \frac{\sqrt{t \rho^2 + 1}}{\alpha} \right )}
\end{equation}
    forms a $(1-\alpha)$-\asympcs{} for $\theta$. Moreover, if $\rho$ is replaced by $\rho_m := \rho(d_m m)$ as in \cref{eq:type-I-error-iid} to obtain
  \begin{equation}
    C_t(m) := \left ( \widehat \theta_t \pm \sqrt{\Var(\phi)} \sqrt{\frac{2(t \rho_m^2 + 1)}{t^2\rho_m^2} \log \left ( \frac{\sqrt{t \rho_m^2 + 1}}{\alpha} \right )} \right ),
  \end{equation}
  then $\infseqt{C_t(m)}$ has asymptotic time-uniform coverage as $m \to \infty$, meaning
  \begin{equation}
    \liminf_{m \to \infty}\PP(\forall t \geq m,\ \theta \in C_t(m)) = 1-\alpha.
  \end{equation}
\end{corollary}
Clearly, the boundaries of \namecref{section:lil}s~\ref{section:lil} or~\ref{section:delayed-start} could be used here in place of Theorem~\ref{theorem:acs}, though for the boundary in \cref{proposition:acsLIL}, we would need to strengthen the $o(\sqrt{\log t / t})$ rate in \eqref{eq:asymp-linear-time-uniform-estimator} to $o(\sqrt{\log \log t / t})$. 
If computing $\widehat \theta_t$ additionally involves the estimation of a nuisance parameter $\eta$ such as in Theorems~\ref{theorem:csate_randomized} and \ref{theorem:csate_observational}, this must be handled carefully on a case-by-case basis where sequential sample splitting and cross fitting (\cref{section:seqSampleSplit}) may be helpful, and higher moments on $\phi(Z_i)$ may be needed.
We now derive an analogue of the delta method for asymptotically linear time-uniform estimators.

\begin{proposition}[The delta method for \asympcs{}s]\label{proposition:delta-method}
  Consider the same setup as \cref{corollary:general-functionals} and let $g: \RR \to \RR$ be a continuously differentiable function with first derivative $g'$. Then,
  $g(\widehat \theta_t)$ is an asymptotically linear time-uniform estimator for $g(\theta)$ with influence function given by $g'(\theta)\phi(\cdot)$, i.e.
  \begin{equation}
g(\widehat \theta_t) - g(\theta) \eqas \frac{1}{t} \sum_{i=1}^t g'(\theta) \phi(Z_i) + o\left ( \sqrt{\log t / t} \right ).
  \end{equation}
\end{proposition}
In particular, \cref{proposition:delta-method} can be combined with \cref{corollary:general-functionals} to obtain delta method-like \asympcs{}s with asymptotic time-uniform coverage guarantees.
The short proof in \cref{proof:delta-method} is similar to the proof of the classical delta method but with the almost-sure continuous mapping theorem used in place of the in-probability one, and with the law of the iterated logarithm used in place of the central limit theorem.


\section{Simulation studies: Widths and empirical coverage}\label{section:simulations-width-coverage}

We now perform simulations focusing on the setting where $\infseqt{Y_t}$ are bounded random variables (with known bounds) and the parameters of interest may include means or treatment effects. We consider this setting as it is well-studied, allowing us to draw on a rich literature containing several nonasymptotic \cs{}s to which we will compare \asympcs{}s (though it is important to keep in mind that there are many unbounded problems for which nonasymptotic \cs{}s do not exist, and we discuss these at the end of the section). In particular, we will compare the \asympcs{}s of \namecref{theorem:acs}s~\ref{theorem:acs} and~\ref{theorem:csate_randomized} to Robbins' sub-Gaussian mixture \cs{} \citep{robbins1970statistical} (see also \citep[\S 3.2]{howard2018uniform}), the empirical Bernstein \cs{}s of \citet[Thm 4 and Cor 2]{howard2018uniform}, and CLT-based \ci{}s.
Of course, CLT-based \ci{}s are not time-uniform and are only included for reference. \cref{fig:widths_miscoverage} considers three cases of parameter estimation for bounded random variables:

\begin{figure}[!htbp]
  \centering
  \textbf{(a) Estimating the mean of bounded random variables}
  \includegraphics[width=\columnwidth]{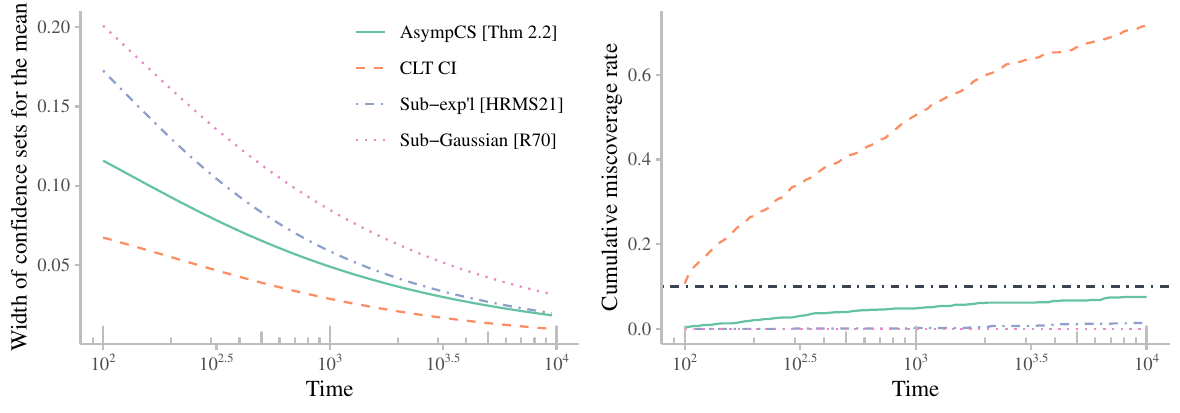}
  \textbf{(b) Estimating the ATE in a completely randomized Bernoulli experiment}
  \includegraphics[width=\columnwidth]{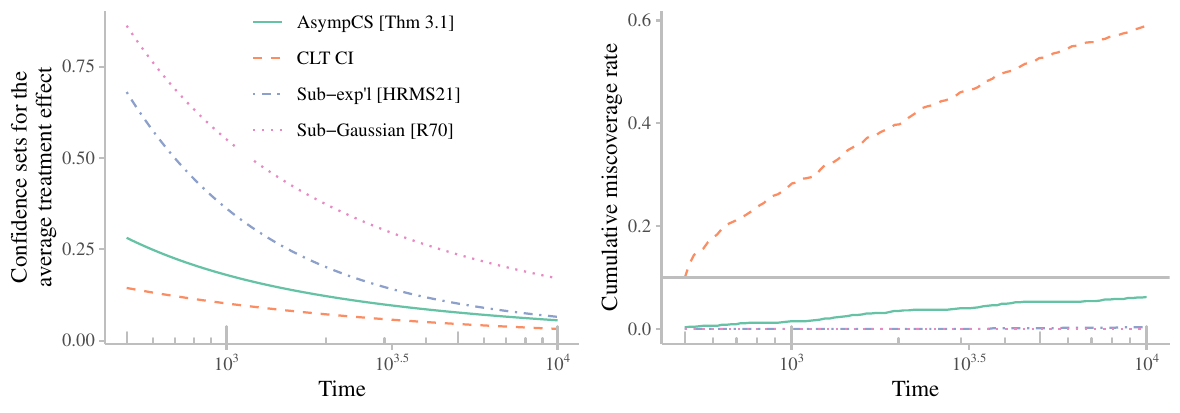}
  \textbf{(c) Estimating the ATE in an experiment with personalized randomization}
  \includegraphics[width=\columnwidth]{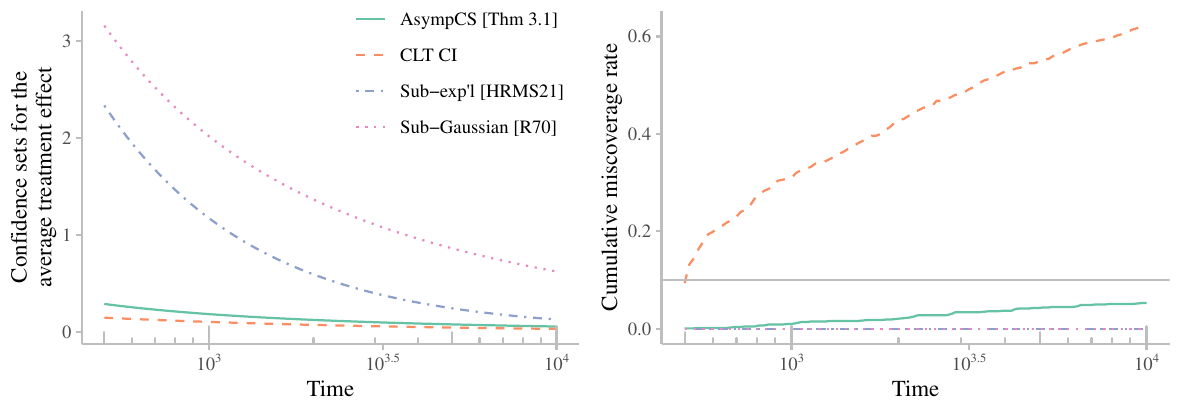}
  \caption{A comparison of $(1-\alpha) \equiv 90$\% confidence sets for parameters in three simulation scenarios: (a) mean estimation of bounded random variables, (b) ATE estimation from a completely randomized Bernoulli experiment, and (c) ATE estimation from an experiment with covariate-dependent (``personalized'') randomization. Empirical widths and miscoverage rates were computed with $1000$ replications beginning at time $100$ for mean estimation (Simulation (a)) and time $500$ for ATE estimation (Simulations (b) and (c)). Notice that in all three scenarios, only \cs{}s have miscoverage rates below $\alpha$, but \asympcs{}s are the only ones that appear to sharply approach this level. The tuning parameter $\widebar \rho_m$ was chosen for a start time of $m$ as $\widebar \rho_m := \rho(m \log (m\lor e))$ following~\eqref{eq:type-I-error-iid}.}
  \label{fig:widths_miscoverage}
\end{figure}

\begin{enumerate}[label = (\alph*)]

\item The first simulation focuses on estimating the mean $\mu$ of \iid{} Uniform(0, 1) random variables only using knowledge of $[0, 1]$-boundedness. For this setting, boundedness allows the \cs{}s of \citet{robbins1970statistical} and \citet{howard2018uniform} to be immediately applicable.

\item The second considers average treatment effect estimation in a randomized experiment with $\{0, 1\}$-valued outcomes where everyone is randomly assigned to treatment or control with probability $1/2$. Since all propensity scores are equal to $1/2$, we have that estimates of the influence functions in \eqref{eq:randomizedExperimentInfluenceFunction} from \cref{section:randomized} are bounded in $[-2, 2]$, and hence the techniques of \citet{robbins1970statistical} and \citet{howard2018uniform} are applicable (as suggested by \citet[\S 4.2]{howard2018uniform}).

\item The third and final simulation considers a similar setup to (b) with the key difference that propensity scores are now covariate-dependent (``personalized''). In this case, as suggested by \citet[\S 4.2]{howard2018uniform}, note that estimates of the influence functions \eqref{eq:randomizedExperimentInfluenceFunction} lie in $[-1/p_\mathrm{min} - 1, 1/p_\mathrm{min} + 1]$ where $p_\mathrm{min} := \min\{\essinf_{x} \pi(x), \essinf_{x} [1-\pi(x)]\}$, permitting the application of \citet{robbins1970statistical} and \citet{howard2018uniform} as before.\footnote{Note that \citet[\S 4.2]{howard2018uniform} use a more conservative bound of $[-2/p_\mathrm{min}, 2/p_\mathrm{min}]$ but it is not hard to see that this can be improved to $[-1/p_\mathrm{min} - 1, 1/p_\mathrm{min} + 1]$. Consequently, we are ultimately comparing our \asympcs{}s to a \emph{strictly tighter} nonasymptotic \cs{} than the one proposed by \citet[\S 4.2]{howard2018uniform}.}
\end{enumerate}

\begin{figure}[!htbp]
  \centering
  \includegraphics[width=\textwidth]{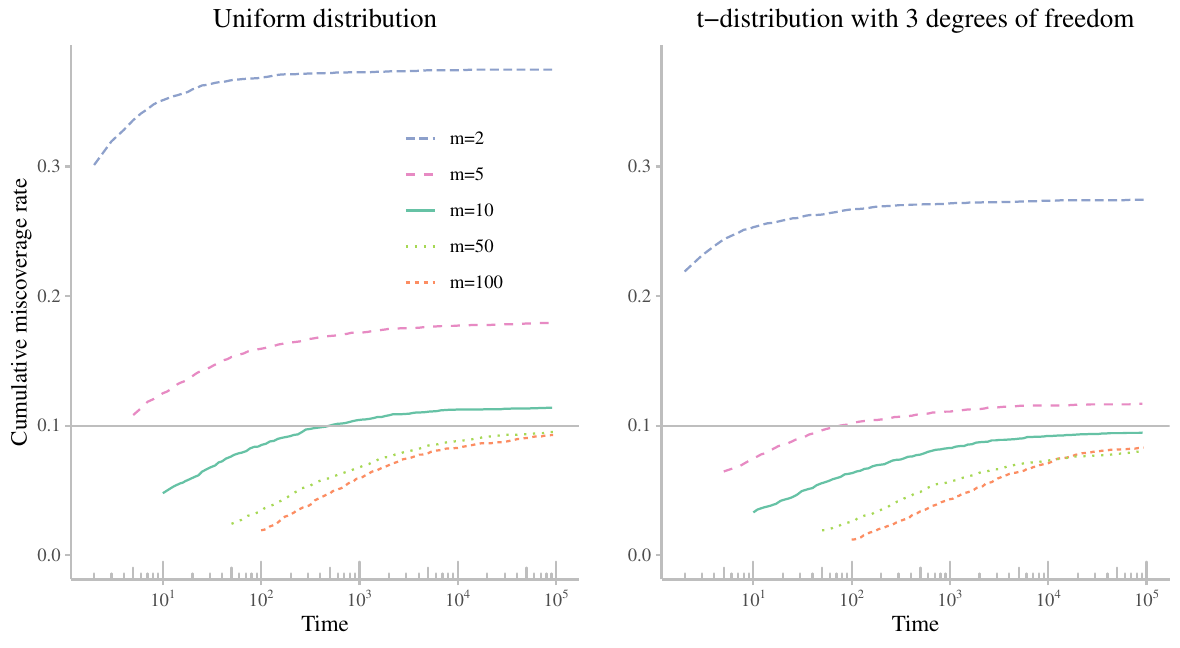}
  \caption{Cumulative miscoverage rates using \eqref{eq:type-I-error-iid} at level $\alpha = 0.1$ to estimate the mean of \iid{} Uniform[0, 1] and $t$-distributed random variables in the left-hand and right-hand side plots, respectively. Notice that in both cases including the heavy-tailed setting of a $t$-distribution with 3 degrees of freedom (so that the variance is finite but third and higher absolute moments are all infinite), cumulative miscoverage rates do not exceed $\alpha = 0.1$ even after $10^5$ observations as long as the first peeking time $m$ is at least 50. It is worth remarking that asymptotic approximations appear to ``kick in'' earlier for the heavy-tailed $t$-distribution.}
  \label{fig:first-peeking-time}
\end{figure}

Notice that in all three scenarios, CLT-based \ci{}s have cumulative miscoverage rates that quickly diverge beyond $\alpha = 0.1$ while those of \cs{}s --- both asymptotic and nonasymptotic --- never exceed $\alpha$ before time $10^4$. (Note that a longer time horizon of $10^5$ is considered only for \asympcs{}s and CLT-based \ci{}s in \cref{fig:CSvsCI}, but is omitted here due to the computational expense of \citet[Thm 2]{howard2018uniform} at large $t$.) Moreover, notice that nonasymptotic \cs{}s appear to be conservative, while our \asympcs{}s are much tighter and have miscoverage rates approaching $\alpha$ (as expected in light of \cref{theorem:type-I-error}). Of course, \cref{theorem:type-I-error} states that miscoverage will only be close to $\alpha$ for large values of the first peeking time $m$. In \cref{fig:first-peeking-time}, we illustrate this phenomenon by using the bound in \eqref{eq:type-I-error-iid} to estimate the means of two distributions (uniform and Student $t$-distributed in the left-hand and right-hand side plots, respectively) given several initial peeking times $m \in \{2, 5, 10, 50, 100\}$.

These particular simulation scenarios illustrate situations where \asympcs{}s are increasingly beneficial over nonasymptotic bounds. First, estimating means of bounded random variables (\cref{fig:widths_miscoverage}(a)) is a problem for which several nonasymptotic \ci{}s and \cs{}s exist, and indeed both the bounds of \citet{robbins1970statistical} and especially \citet{howard2018uniform} fare well in this setting. The relatively small variance of Uniform(0, 1) random variables means that asymptotic methods can quickly tighten in this setting while the empirical Bernstein \cs{}s of \citet{howard2018uniform} take a while to adapt to the variance (while those of \citet{robbins1970statistical} never will).

\asympcs{}s are particularly well-suited to the settings of ATE estimation in Figs~(b) and (c) where almost-sure bounds on observed random variables can be quite large in comparison to their variances. \cref{fig:widths_miscoverage}(b) considers the setting where all propensity scores are equal to $1/2$ which is the ``easiest'' regime for nonasymptotic methods, meaning that a.s.~bounds on the influence functions are as tight as possible. Even here, \asympcs{}s are tighter than the nonasymptotic bounds. On the other hand, \cref{fig:widths_miscoverage}(c) considers a setting where propensity scores are highly covariate-dependent so that bounds on propensity scores $\pi(x) \in [p_\mathrm{min}, 1-p_\mathrm{min}]$ must hold for almost all $x$ (in this simulation, $p_\mathrm{min} = 0.2$). Here we see that \asympcs{}s drastically outperform nonasymptotic \cs{}s without inflating miscoverage rates. Taking this to a logical extreme, it is possible to construct scenarios where $p_\mathrm{min}$ is closer and closer to 0, but influence function variances remain bounded. In other words, it is possible to consider scenarios where \asympcs{}s are arbitrarily tighter than nonasymptotic \cs{}s, without inflating miscoverage rates.

Finally, we remark that while \asympcs{}s demonstrate substantial benefits over nonasymptotic \cs{}s in terms of \emph{tightness}, we wish to also highlight their benefits of \emph{versatility}, and in particular, note that there are many settings for which no simulations could be run since \asympcs{}s provide the first (asymptotically) time-uniform solution in the literature. For example, as a consequence of \citet{bahadur1956nonexistence}, it is impossible to derive nonasymptotic \cs{}s (or \ci{}s) for the mean of random variables without prior bounds on their moments. By contrast, \asympcs{}s can handle mean estimation under finite (but unknown) moment assumptions. It is also impossible to derive nonasymptotic \cs{} (and \ci{}s) for the ATE from observational studies (without substantial and unrealistic knowledge of nuisance function estimation errors) but \cref{section:observational} outlines an asymptotically time-uniform solution. In both settings, we do not run simulations akin to \cref{fig:widths_miscoverage} since there are no prior \cs{}s to compare to.



\section{Real data application: effects of IV fluid caps in sepsis patients}
\label{section:illustration}

Let us now illustrate the use of \cref{theorem:csate_observational} by sequentially estimating the effect of fluid-restrictive strategies on mortality in an observational study of real sepsis patients. We will use data from the Medical Information Mart for Intensive Care III (MIMIC-III), a freely available database consisting of health records associated with more than 45,000 critical care patients at the Beth Israel Deaconess Medical Center \citep{mimiciii, mimiciiidata}. The data are rich, containing demographics, vital signs, medications, and mortality, among other information collected over the span of 11 years.

Following \citet{shahn2020fluid}, we aim to estimate the effect of restricting intravenous (IV) fluids within 24 hours of intensive care unit (ICU) admission on 30-day mortality in sepsis patients. In particular, we considered patients at least 16 years of age satisfying the Sepsis-3 definition --- i.e. those with a suspected infection and a Sequential Organ Failure Assessment (SOFA) score of at least 2 \citep{singer2016third}. Sepsis-3 patients can be obtained from MIMIC-III using SQL scripts provided by \citet{alistair_johnson_2018_1256723}, but we provide detailed instructions for reproducing our data collection and analysis process on GitHub.\footnote{\reproduceDataLink{}} This resulted in a total of 5231 sepsis patients, each of whom received out-of-hospital followup of at least 90 days.

\begin{figure}[!htbp]
    \centering
    \includegraphics[width=\columnwidth]{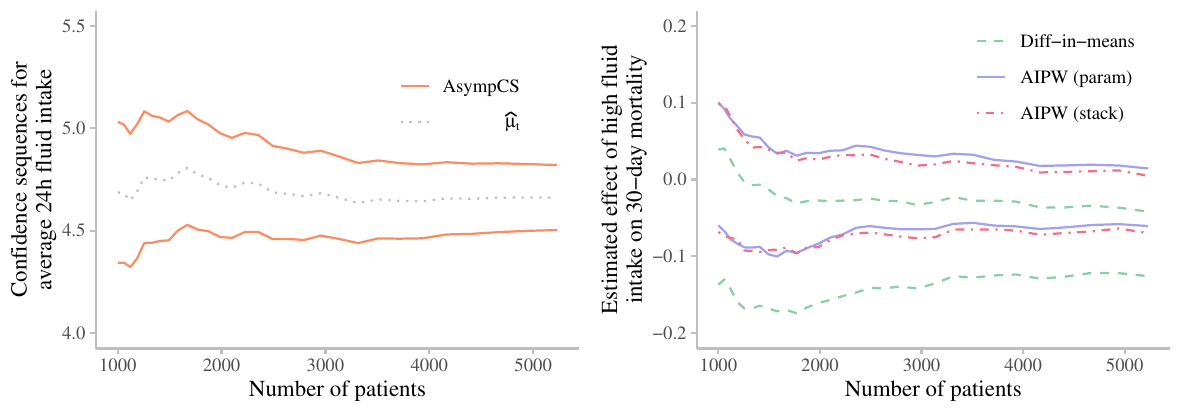}
    \caption{Left-hand side: a $90\%$ \asympcs{} used to track the average 24h fluid intake over time. Right-hand side: Three 90\%-\asympcs{}s for the causal effect of capped IV fluid intake (defined as $\leq 6$ litres) on 30-day mortality using the same three estimators as those outlined in~\cref{fig:CS_observational}. Notice that an analysis using a difference-in-means estimator would conclude that the treatment effect is negative after observing fewer than 1500 patients.}
    \label{fig:sepsis}
\end{figure}

Consider IV fluid intake within 24 hours of ICU admission $\fluid$. To construct a binary treatment $A \in \{0, 1\}$, we dichotomize $\fluid$ so that $A_i = \Ind(\fluid_i \leq 6\mathrm{L})$. 30-day mortality $Y$ is defined as 1 if the patient died within 30 days of hospital admission, and 0 otherwise. We will consider baseline covariates $X$ including a patient's age and sex, whether they are diabetic, modified Elixhauser scores \citep{van2009modification}, and SOFA scores. We are interested in the causal estimand
\begin{equation}
  \psi := \PP\left (Y^{\fluid \leq 6L} = 1 \right) - \PP\left (Y^{\fluid > 6L} = 1 \right),
\end{equation}
i.e.~the difference in average 30-day mortality that would be observed if all sepsis patients were randomly assigned an IV fluid level according to the lower truncated distribution $\PP(\fluid \leq l \mid \fluid \leq 6\mathrm{L})$ versus the upper truncated distribution $\PP(\fluid \leq l \mid \fluid > 6\mathrm{L})$ \citep{munoz2012population}. While this is technically a stochastic intervention effect, we have that under the same causal identification assumptions discussed in \cref{section:csate}, $\psi$ is identified as
\begin{equation}
  \psi = \EE \left \{ \EE(Y \mid X, A=1) - \EE(Y \mid X, A = 0) \right \},
\end{equation}
the same functional considered in the previous sections. Therefore, we can estimate $\psi$ under the same assumptions and with the same techniques as \cref{section:observational}. \cref{fig:sepsis} contains \asympcs{}s for $\psi$ using difference-in-means, parametric AIPW, and stacking-based AIPW estimators to demonstrate the impacts of different modeling choices on \asympcs{} width. 
\begin{remark}
These simple binary treatment and outcome variables were used for simplicity so that the methods outlined in \cref{section:observational} are immediately applicable, but  \cref{section:generalFunctionals} points out that our \asympcs{}s may be used to sequentially estimate other causal functionals.
\end{remark}

The stacking-based AIPW \asympcs{}s cover the null treatment effect of 0 from the $1000^\text{th}$ to the $5231^\text{st}$ observed patient, and thus we cannot conclude whether 6L IV fluid caps have an effect on 30-day mortality in sepsis patients. Interestingly, the width and center of \asympcs{}s based on both stacking and parametric estimators are roughly equal, which is in contrast to the simulations provided in the previous sections. This could be due to the fact that the true regression or propensity score functions lie in the parametric models considered by ``AIPW (param)'' or because they are good approximations. Since this analysis is not a simulation, we cannot know with certainty one way or another.

Note that these stacking-based \asympcs{}s nearly drop below 0 after observing the $5231^\text{st}$ patient's outcome. If we were using fixed-time confidence intervals, the analyst would need to resist the temptation to resume data collection (e.g. to see whether the null hypothesis $H_0: \psi = 0$ could be rejected with a slightly larger sample size) as this would inflate type-I error rates (as seen in \cref{fig:CSvsCI}). On the other hand, \asympcs{}s flexibly permit precisely this form of continued sampling.


\section{Conclusion}
This paper introduced the notion of an ``asymptotic confidence sequence'' as the time-uniform analogue of an asymptotic confidence interval based on the central limit theorem. We derived an explicit universal asymptotic confidence sequence for the mean from \iid{} observations under weak moment assumptions by appealing to the strong invariance principle of \citet{strassen1964invariance}. These results were extended to the setting where observations' distributions (including means and variances) can vary over time under martingale dependence, such that our confidence sequences capture a moving parameter --- the running average of the conditional means so far.
We then applied the aforementioned results to the problem of doubly robust sequential inference for the average treatment effect in both randomized experiments and observational studies under \iid{} sampling. Finally, we showed how these causal applications remain valid in the non-\iid{} setting where distributions change over time, in which case our confidence sequences capture a running average of individual treatment effects. The aforementioned results will enable researchers to continuously monitor sequential experiments --- such as clinical trials and online A/B tests --- as well as sequential observational studies even if treatment effects do not remain stationary over time.



\subsection*{Acknowledgements}
IW-S thanks Tudor Manole for early conversations regarding strong approximation theory. AR acknowledges funding from NSF Grant DMS1916320, NSF Grant DMS2053804, an Adobe faculty research award, and an NSF DMS (CAREER) 1945266. EK gratefully acknowledges support from NSF Grant DMS1810979. The authors thank Ziyu (Neil) Xu, Dae Woong Ham, Abhinandan Dalal, Eric J. Tchetgen Tchetgen, and the CMU causal inference working group for helpful discussions.


\bibliographystyle{plainnat}
\bibliography{references}

\newpage
\appendix

\section{Proofs of the main results}

\subsection{Proof of Theorem~\ref{theorem:acs}}
\label{proof:acs}

We first introduce a lemma that will later be used in the main proof.

\begin{lemma}[Strong Gaussian approximation of the sample average]
\label{lemma:strongGaussianApproxSampleAvg}
    Let $(Y_t)_{t=1}^\infty$ be an \iid{} sequence of random variables with mean $\mu$, variance $\sigma^2$, and let $\widehat \sigma_t^2$ be the sample variance. Then (after sufficiently enriching the probability space), there exist \iid{} Gaussian random variables $(G_t)_{t=1}^\infty$ such that
    \[ \frac{1}{t} \sum_{i=1}^t (Y_i - \mu) \eqas \frac{\widehat \sigma_t}{t} \sum_{i=1}^t G_i + \varepsilon_t \]
    where $\varepsilon_t = o(\sqrt{\log \log t / t})$. One could also replace $\sum_{i=1}^t G_i$ with a Wiener process $(W_t)_{t\geq 0}$.
\end{lemma}
\begin{proof}
    By the strong approximation theorem of \citet{strassen1964invariance}, we have that (after sufficiently enriching the probability space) there exist \iid{} Gaussian random variables $\infseq Gt1$ such that
    \begin{equation}
    \label{eq:strassen}
         \frac{1}{t} \sum_{i=1}^t (Y_i - \mu) = \frac{\sigma}{t} \sum_{i=1}^t G_i + \kappa_t 
    \end{equation}
    with $\kappa_t = o(\sqrt{\log \log t / t})$. Since $(G_t)_{t=1}^\infty$ have mean zero and unit variance, we have by the law of the iterated logarithm (LIL) that
    \begin{equation}
    \label{eq:GaussianLIL}
         \frac1t \sum_{i=1}^t G_i = O\left ( \sqrt{\frac{\log \log t}{t}} \right ). 
    \end{equation}
    Now, by the strong law of large numbers, $\widehat \sigma_t \xrightarrow{a.s.} \sigma$. Combining this fact with with \eqref{eq:strassen}, we have
    \begin{align}
        \frac{1}{t} \sum_{i=1}^t (Y_i - \mu) &= \frac{\widehat \sigma_t + o(1)}{t} \sum_{i=1}^t G_i + \kappa_t \\
        &= \frac{\widehat \sigma_t}{t} \sum_{i=1}^t G_i + \kappa_t + o\left ( \sqrt{\frac{\log \log t}{t}} \right )\\
        &= \frac{\widehat \sigma_t}{t} \sum_{i=1}^t G_i + o\left ( \sqrt{\frac{\log \log t}{t}} \right )
    \end{align}
    almost surely, which completes the proof.
\end{proof}

\begin{proof}[Proof of \cref{theorem:acs}]
The proof proceeds in 3 steps. First, we use the fact that for any martingale $M_t(\lambda)$, we have that the mixture $\int_{\RR} M_t(\lambda) dF(\lambda)$ is also a martingale where $F$ is any probability distribution on $\RR$ \citep{ howard_exponential_2018, howard2018uniform}. We apply this fact to an exponential Gaussian martingale and use a Gaussian density $f(\lambda; 0, \rho^2)$ as the mixing distribution. Second, we apply Ville's inequality \citep{ville1939etude} to this mixture exponential Gaussian martingale to obtain Robbins' normal mixture confidence sequence \citep{robbins1970statistical}. Third, we use Lemma~\ref{lemma:strongGaussianApproxSampleAvg} to approximate $\sum_{i=1}^t (Y_i - \mu)$ by a cumulative sum of Gaussian random variables and apply the results from steps 1 and 2.

\paragraph*{Step 1} Let $(W_t)_{t\geq 0}$ be a Wiener process and write the exponential process for any $\lambda \in \RR$,
\[ M_t(\lambda) := \exp\left \{ \lambda W_t - t \lambda^2 / 2  \right \}. \]
It is well-known that $M_t(\lambda)$ is a \emph{nonnegative martingale starting at $M_0 \equiv 1$} with respect to the Brownian filtration. By Fubini's theorem, for any probability distribution $F(\lambda)$ on $\RR$, we also have that the mixture,
\[ \int_{\lambda \in \RR} M_t(\lambda) dF(\lambda) \]
is a nonnegative martingale with initial value one with respect to the same filtration \citep{robbins1970statistical}. In particular, consider the Gaussian probability distribution function $f(\lambda; 0, \rho^2)$ with mean zero and variance $\rho^2 > 0$ as the mixing distribution. The resulting martingale can be written as
\begin{align}
  M_t &:= \int_{\lambda \in \mathbb R} \exp\left \{ \lambda W_t - \frac{t \lambda^2}{2} \right \}f(\lambda; 0, \rho^2) d\lambda \\
      &= \frac{1}{\sqrt{2\pi \rho^2}} \int_\lambda \exp \left \{ \lambda W_t - \frac{t\lambda^2}{2} \right \}\exp \left \{ \frac{-\lambda^2}{2\rho^2} \right \} d\lambda \\
      &= \frac{1}{\sqrt{2\pi \rho^2}} \int_\lambda \exp \left \{ \lambda W_t - \frac{\lambda^2 (t\rho^2 + 1)}{2 \rho^2}\right \} d\lambda \\
      &= \frac{1}{\sqrt{2\pi \rho^2}} \int_\lambda \exp \left \{ \frac{-\lambda^2 (t\rho^2 + 1) + 2\lambda \rho^2 W_t }{2\rho^2} \right \} d\lambda \\
      &= \frac{1}{\sqrt{2\pi \rho^2}} \int_\lambda \exp \left \{ \frac{-a(\lambda^2 - \frac{b}{a} 2\lambda) }{2\rho^2} \right \} d\lambda
\end{align}
by setting $a:= t\rho^2 + 1$ and $b := \rho^2 W_t$. Focusing on the integrand and completing the square, we have
\begin{align}
    \exp \left \{ \frac{-\lambda^2 + 2\lambda \frac{b}{a} + \left ( \frac{b}{a} \right )^2 - \left ( \frac{b}{a} \right )^2 }{2 \rho^2 /a} \right \} &= \exp \left \{ \frac{-(\lambda - b/a)^2}{2\rho^2/a} + \frac{a \left ( b/a \right )^2}{2\rho^2} \right \} \\
    &= \underbrace{\exp \left \{ \frac{-(\lambda - b/a)^2}{2\rho^2/a} \right \}}_{\propto f(\lambda, b/a, \rho^2/a)} \exp \left \{  \frac{b^2}{2a\rho^2} \right \}. 
\end{align}
Plugging this back into the integral and multiplying the entire quantity by $\frac{a^{-1/2}}{a^{-1/2}}$, we finally get the closed-form expression of the mixture exponential Wiener process,
\begin{align}
    M_t &:= \underbrace{\frac{1}{\sqrt{2\pi \rho^2/a}} \int_{\lambda \in \mathbb R} \exp\left \{ \frac{-(\lambda - b/a)^2}{2 \rho^2/a} \right \} d\lambda}_{=1} \frac{\exp \left \{ \frac{b^2}{2a\rho^2} \right \} }{\sqrt{a}} \nonumber \\ 
    &=\frac{\exp \left \{ \frac{\rho^2 W_t^2}{2(t\rho^2 + 1)} \right \} }{\sqrt{t\rho^2 + 1}}.
    \label{eq:mixed-exponential-wiener}
\end{align}

\paragraph*{Step 2} 
Since $M_t$ is a nonnegative martingale with initial value one, we have by Ville's inequality \citep{ville1939etude} that
\[ \PP( \forall t \geq 1,\ M_t < 1/\alpha) \geq 1-\alpha. \]
Writing this out explicitly for $M_t$ and solving for $W_t$ algebraically, we have that 
\begin{align}
    &\PP\left ( \forall t \geq 1,\ \frac{\rho^2 W_t^2}{2(t\rho^2 + 1)} < \log (1/\alpha) + \log \left (\sqrt{t\rho^2 + 1} \right) \right )\\
    =\ &\PP \left ( \forall t \geq 1,\ \left | W_t/t \right | < \sqrt {\frac{2(t\rho^2 + 1)}{t^2 \rho^2} \log \left ( \frac{\sqrt{t\rho^2 + 1}}{ \alpha } \right )} \right ) \geq 1 - \alpha.
\end{align}

\paragraph*{Step 3} First, note that by the triangle inequality,
\[ \left |\frac{1}{t} \sum_{i=1}^t Y_i - \mu \right | \leq \left | \frac{1}{t} \sum_{i=1}^t (Y_i - \mu) - \frac{\widehat \sigma_t}{t}W_t \right | + \frac{\widehat \sigma_t}{t} \left | W_t \right |, \]
and thus by Lemma~\ref{lemma:strongGaussianApproxSampleAvg} and Step 2, we have with probability at least $(1-\alpha)$, 
\begin{align}
    \forall t \geq 1,\ \left |\frac{1}{t} \sum_{i=1}^t Y_i - \mu \right | < \underbrace{\widehat \sigma_t\sqrt {\frac{2(t\rho^2 + 1)}{t^2\rho^2} \log \left ( \frac{\sqrt{t\rho^2 + 1}}{ \alpha } \right )}}_{\widebar \boundary_t^\Gcal} + |\varepsilon_t|
\end{align} 
where $\varepsilon_t$ is defined as in Lemma~\ref{lemma:strongGaussianApproxSampleAvg}. Finally, note that
\begin{equation}
  \frac{\widebar \boundary_t^\Gcal + |\eps_t|}{\widebar \boundary_t^\Gcal} \xrightarrow{\as} 1,
\end{equation}
completing the proof.
\end{proof}

\subsection{Proof of \cref{theorem:general-asympcs}}\label{proof:general-asympcs}
\begin{proof}
  Using \namecref{condition:strong-approx}s~\ref{condition:strong-approx} and~\ref{condition:boundary-for-approximating-process}, take $L_t^\star := \widehat L_t + (\widehat \theta_t - \theta_t + Z_t)$ and $U_t^\star := \widehat U_t - (\widehat \theta_t - \theta_t + Z_t)$. Then we have that
  \begin{align}
    &\PP \left ( \forall t \in \Tcal,\ \theta_t \in \left [\widehat \theta_t - L_t^\star, \widehat \theta_t + U_t^\star \right]  \right )\\
    =\ &\PP \left ( \forall t \in \Tcal,\ \theta_t \in \left [\widehat \theta_t - \widehat L_t - (\widehat \theta_t - \theta_t + Z_t), \widehat \theta_t + \widehat U_t - (\widehat \theta_t - \theta_t + Z_t) \right]  \right )\\
    =\ &\PP \left ( \forall t \in \Tcal,\ Z_t \in [-\widehat L_t, \widehat U_t ]  \right ) \geq 1-\alpha.
  \end{align}
  It remains to show that 
  \begin{equation}
    L_t^\star/ L_t \xrightarrow{\as} 1\quad\text{and}\quad U_t^\star/ U_t \xrightarrow{\as} 1.
  \end{equation}
  Indeed, we have that
  \begin{align}
    L_t^\star / L_t &= (L_t^\star / \widehat L_t) \cdot (\widehat L_t / L_t)\\
                    &= \left ( \frac{\widehat L_t + (\widehat \theta_t - \theta_t + Z_t)}{\widehat L_t} \right ) \cdot \left ( \widehat L_t / L_t \right )\\
                    &= \left ( 1 + o(1) \right ) \cdot \left ( 1 + o(1) \right ) \\
    &= 1 + o(1)
 \end{align} 
 almost surely, where the second-last line follows from the combination of \namecref{condition:strong-approx-rate}s~\ref{condition:strong-approx-rate} and~\ref{condition:boundary-approximation}. A similar calculation goes through for $U_t^\star / U_t$, completing the proof.
\end{proof}

\subsection{Proof of \cref{theorem:lindeberg-martingale-asympcs}}
\label{proof:lyapunov-type-asympcs}

First, we present a lemma that is implicit in the proof of \citet[Theorem 4.4]{strassen1967and} and which will be central to the proof of \cref{theorem:lindeberg-martingale-asympcs}.
\begin{lemma}[Strong approximation under martingale dependence]\label{lemma:strong-approx-martingale}
  Let $\infseqt{Y_t}$ be a sequence of random variables with conditional means and variances given by $\mu_t = \EE(Y_t \mid Y_1^{t-1})$ and $\sigma_t^2 = \Var(Y_t \mid Y_1^{t-1})$, respectively. Let $V_t = \sum_{i=1}^t \sigma_i^2$ be the cumulative conditional variance process, and suppose that $V_t \to \infty$ almost surely. Furthermore, assume that the Lindeberg-type condition given in \cref{assumption:lindeberg-type_condition} holds. Then, on a potentially enriched probability space, there exist \iid{} standard Gaussians $\infseqt{G_t}$ such that
  \begin{equation}
    \frac{1}{t}\sum_{i=1}^t (Y_i - \mu_i) - \frac{1}{t}\sum_{i=1}^t \sigma_i G_i = o \left ( \frac{V_t^{3/8} \log V_t}{t} \right ).
  \end{equation}
\end{lemma}

\begin{proof}[Proof of \cref{lemma:strong-approx-martingale}]
  The proof centrally relies on \citet[Eq. 159]{strassen1967and} which states that on a potentially enriched probability space,
  \begin{equation}
    \sum_{i=1}^t (Y_i - \mu_i) = \xi(V_t) + o(h(V_t))
  \end{equation}
  where $\xi$ is a standard Brownian motion, and $h(v) = (v f(v))^{1/4} \log v$ is any function so that $f(v)$ is increasing in $v$, but $f(v)/v$ is decreasing. For the purposes of this proof, we set $f(v) = v^{1/2}$, and hence $h(v) = v^{3/8} \log v$. Since a standard Brownian motion evaluated at $V_t = \sum_{i=1}^t \sigma_i^2$ is equal in distribution to the discrete time process $\sum_{i=1}^t \sigma_i G_i$ at each $t$, we have that
  \begin{equation}
    \sum_{i=1}^t (Y_i - \mu_i) = \sum_{i=1}^t \sigma_i G_i + o \left ( V_t^{3/8} \log V_t \right ),
  \end{equation}
  which completes the proof after dividing both sides by $t$.
\end{proof}
With \cref{lemma:strong-approx-martingale} in mind, we can now prove the main result (\cref{theorem:lindeberg-martingale-asympcs}).

\begin{proof}[Proof of \cref{theorem:lindeberg-martingale-asympcs}]
  The proof proceeds in four steps, each step being dedicated to satisfying one of \cref{theorem:general-asympcs}'s conditions (\ref{condition:strong-approx}--\ref{condition:boundary-approximation}). Step 1 follows quickly from \cref{lemma:strong-approx-martingale}, while Step 2 requires us to derive a (sub)-Gaussian boundary for non-\iid{} observations under martingale dependence. Step 3 follows from the arguments in Steps 1 and 2, and Step 4 follows from a simple argument that uses \cref{assumption:lindeberg_variance_estimator_consistent}.

  \paragraph*{Step 1: Satisfying \cref{condition:strong-approx} via \citet{strassen1967and} and \namecref{assumption:lindeberg_variance_does_not_vanish}s~\ref{assumption:lindeberg_variance_does_not_vanish} and~\ref{assumption:lindeberg-type_condition}}
  Notice that the assumptions of \cref{lemma:strong-approx-martingale} are satisfied by \namecref{assumption:lindeberg-type_condition}s~\ref{assumption:lindeberg_variance_does_not_vanish} and~\ref{assumption:lindeberg-type_condition}, and hence \cref{condition:strong-approx} follows using the approximating process formed by $\frac{1}{t}\sum_{i=1}^t \sigma_i G_i$ as in \cref{lemma:strong-approx-martingale} with a rate of $r_t := t^{-1} V_t^{3/8} \log V_t$. That is,
  \begin{equation}
    \left ( \widehat \mu_t - \widetilde \mu_t \right ) - \frac{1}{t} \sum_{i=1}^t \sigma_i G_i = o \left ( \frac{V_t^{3/8} \log V_t}{t} \right ).
  \end{equation}
  In Step 2, we provide a nonasymptotic boundary for the process given by $\frac{1}{t}\sum_{i=1}^t \sigma_i G_i$.
  \paragraph*{Step 2: Satisfying \cref{condition:boundary-for-approximating-process} using a nonasymptotic sub-Gaussian boundary}
Let $\sigma_i G_i$ be as in \cref{lemma:strong-approx-martingale}. Define the conditional variance $\sigma_t^2 := \Var(Y_t \mid Y_1^{t-1})$ and note that
\[ \widetilde M_t(\lambda) := \exp \left \{ \sum_{i=1}^t (\lambda \sigma_i G_i - \lambda^2 \sigma_i^2 / 2) \right \} \]
is a nonnegative martingale starting at one (with respect to the filtration generated by $\infseqt{G_t, \sigma_t}$). Mixing over $\lambda$ with the probability density $dF(\lambda)$ of a mean-zero Gaussian with variance $\rho^2$ as in the proof of Theorem~\ref{theorem:acs}, we have that
\[ \widetilde M_t := \int_{\lambda \in \RR} \widetilde M_t(\lambda) dF(\lambda) =\exp \left \{ \frac{\rho^2 (\sum_{i=1}^t \sigma_i G_i)^2}{2(V_t \rho^2 + 1)} \right \}\cdot \left(V_t \rho^2 + 1\right)^{-1/2}\]
is also a martingale. By Ville's inequality for nonnegative (super)martingales, we have that $\PP(\exists t : \widetilde M_t \geq 1/\alpha) \leq \alpha$ and hence with probability at least $(1-\alpha)$,
\begin{equation}
\label{eq:nonid_bound}
    \forall t \geq 1, \ \left |\frac{1}{t} \sum_{i=1}^t \sigma_i G_i \right | < \sqrt{\frac{2(V_t \rho^2 + 1)}{t^2 \rho^2} \log \left ( \frac{\sqrt{V_t \rho^2 + 1}}{\alpha} \right )}.
\end{equation}
In particular, combined with Step 1, we have that
\begin{equation}
  \label{eq:asympcs-widetilde-sigma}
  \left ( \widehat \mu_t \pm \sqrt{\frac{2(V_t \rho^2 + 1)}{t^2 \rho^2} \log \left ( \frac{\sqrt{V_t \rho^2 + 1}}{\alpha} \right )}\right )
\end{equation}
forms a $(1-\alpha)$-\asympcs{} for $\widetilde \mu_t$.

\paragraph*{Step 3: Satisfying \cref{condition:strong-approx-rate} as a consequence of \namecref{condition:boundary-for-approximating-process}s~\ref{condition:strong-approx} and~\ref{condition:boundary-for-approximating-process}}
Inspecting the boundary in \eqref{eq:nonid_bound}, we notice that $V_t^{3/8} \log V_t / t = o(\sqrt{V_t \log V_t} / t)$ and hence \cref{condition:strong-approx-rate} is satisfied.

  \paragraph*{Step 4: Satisfying \cref{condition:boundary-approximation} via \cref{assumption:lindeberg_variance_estimator_consistent}}

Writing out the margin of \eqref{eq:asympcs-widetilde-sigma} combined with \cref{assumption:lindeberg_variance_estimator_consistent}: $\widehat \sigma_t^2 - \widetilde \sigma_t^2 = o(\widetilde \sigma_t^2)$ and recalling that $V_t := t\widetilde \sigma_t$, we have
\begin{align}
  \sqrt{\frac{2(V_t \rho^2 + 1)}{t^2 \rho^2} \log\left ( \frac{\sqrt{ V_t \rho^2 + 1}}{\alpha} \right )} & = \sqrt{\frac{2(t(\widehat \sigma_t^2 + o(\widetilde \sigma_t^2)) \rho^2 + 1)}{t^2 \rho^2} \log\left ( \frac{\sqrt{t(\widehat \sigma_t^2 + o(\widetilde \sigma_t^2)\rho^2 + 1}}{\alpha} \right )} \nonumber\\
                                                                                                                                             & = \sqrt{\frac{t(\widehat \sigma_t^2 + o(\widetilde \sigma_t^2) \rho^2 + 1}{t^2 \rho^2} \log\left ( \frac{t(\widehat \sigma_t^2 + o(\widetilde \sigma_t^2)\rho^2 + 1}{\alpha^2} \right )}         \nonumber\\
                                                                                                                                             & = \sqrt{\frac{t\widehat \sigma_t^2\rho^2 + o(t\widetilde \sigma_t^2) + 1}{t^2 \rho^2} \log\left ( \frac{t\widehat \sigma_t^2\rho^2 + o(t\widetilde \sigma_t^2) + 1}{\alpha^2} \right )} \nonumber\\
                                      & = \sqrt{\left (\frac{t\widehat \sigma_t^2\rho^2 + 1}{t^2 \rho^2} + o(\widetilde \sigma_t^2/t) \right ) \log\left ( \frac{t\widehat \sigma_t^2\rho^2 + o(t\widetilde \sigma_t^2) + 1}{\alpha^2} \right )}.\label{eq:margin-asymptotic}
\end{align}
Focusing on the logarithmic factor, we have
\begin{align}
  \log\left ( \frac{t\widehat \sigma_t^2\rho^2 + o(t\widetilde \sigma_t^2) + 1}{\alpha^2} \right ) &= \log\left ( \frac{ 1+t\widehat \sigma^2_t \rho^2}{\alpha^2} + o(t\widetilde \sigma_t^2) \right ) \nonumber\\
                                                                              &= \log\left ( \frac{ 1+t\widehat \sigma^2_t \rho^2}{\alpha^2}\left [ 1 + o\left ( 1 \right ) \right ] \right ) \nonumber\\
  &= \log \left ( \frac{1 + t\widehat \sigma_t^2 \rho^2}{\alpha^2} \right ) + \log \left ( 1 + o\left ( 1 \right )\right )\nonumber \\
                                                                              &= \log \left ( \frac{1 + t\widehat \sigma_t^2 \rho^2}{\alpha^2} \right ) + o\left ( 1 \right ) \label{eq:log-approx}
\end{align}
where the last line follows from the Taylor expansion $\log(1 + x) = x + o(1)$ for $|x| < 1$. Combining \eqref{eq:margin-asymptotic} and \eqref{eq:log-approx}, we have that the margin of \eqref{eq:asympcs-widetilde-sigma} can be written as
\begin{align}
  \sqrt{\frac{2(V_t \rho^2 + 1)}{t^2 \rho^2} \log\left ( \frac{\sqrt{V_t \rho^2 + 1}}{\alpha} \right )} & = \sqrt{\left (\frac{t\widehat \sigma_t^2\rho^2 + 1}{t^2 \rho^2} + o(V_t /t^2) \right ) \left [\log \left ( \frac{1 + t\widehat \sigma_t^2 \rho^2}{\alpha^2} \right )  + o(1)\right]} \nonumber\\
                                                                                                                                             &= \sqrt{\frac{t \widehat \sigma_t^2 \rho^2 + 1}{t^2 \rho^2} \log \left ( \frac{1 + t\widehat \sigma_t^2\rho^2}{\alpha^2} \right ) + o( V_t/t^2) + o(V_t \log V_t / t^2) + o(V_t/t^2)} \nonumber\\
                                                                                                                                             &= \sqrt{\frac{2(t \widehat \sigma_t^2 \rho^2 + 1)}{t^2 \rho^2} \log \left ( \frac{\sqrt{1 + t\widehat \sigma_t^2\rho^2}}{\alpha} \right ) + o(V_t \log V_t / t^2)} \nonumber\\
  &\leq \sqrt{\frac{2(t \widehat \sigma_t^2 \rho^2 + 1)}{t^2 \rho^2} \log \left ( \frac{\sqrt{1 + t\widehat \sigma_t^2\rho^2}}{\alpha} \right )} + o \left ( \frac{\sqrt{V_t \log V_t}}{t} \right ),
\end{align}
where the last inequality follows from $\sqrt{a+ b} \leq \sqrt{a} + \sqrt{b}$ for $a, b \geq 0$. In particular, letting
\begin{align}
  [L_t, U_t] &:= \left [ \widehat \mu_t \pm \sqrt{\frac{2(V_t \rho^2 + 1)}{t^2 \rho^2} \log\left ( \frac{\sqrt{V_t \rho^2 + 1}}{\alpha} \right )} \right ]\\
  \text{and}\quad [\widehat L_t, \widehat U_t] &:= \left [ \widehat \mu_t \pm \sqrt{\frac{2(t \widehat \sigma_t^2 \rho^2 + 1)}{t^2 \rho^2} \log\left ( \frac{\sqrt{t \widehat \sigma_t^2 \rho^2 + 1}}{\alpha} \right )} \right ],
\end{align}
we have that $\widehat L_t / L_t \xrightarrow{\as} 1$ and $\widehat U_t / U_t \xrightarrow{\as} 1$, satisfying \cref{condition:boundary-approximation} and completing the proof of \cref{theorem:lindeberg-martingale-asympcs}.

\end{proof}

\subsection{Proof of \cref{corollary:lyapunov-asympcs}}\label{proof:lyapunov-asympcs}
\begin{proof}
  
  As alluded to in the paragraph following \cref{corollary:lyapunov-asympcs}, it suffices to show that if the following regularity conditions hold, i.e.
  \begin{equation}\label{eq:lyapunov-asympcs-proof-regularity}
    \sum_{i=1}^\infty \frac{\EE |Y_i^2 - \EE Y_i^2|^{1+\beta}}{V_i^{1+\beta}} < \infty, \quad\widetilde \mu_t^2 = o(V_t),\quad \text{and}\quad\frac{1}{t}\sum_{i=1}^t (\mu_i - \widetilde \mu_t)^2 = o(\widetilde \sigma_t^2),
  \end{equation}
  then \cref{assumption:lindeberg_variance_estimator_consistent} is satisfied, meaning $\widehat \sigma_t^2 / \widetilde \sigma_t^2 \xrightarrow{\as} 1$ where $\widehat \sigma_t^2 := \frac{1}{t}\sum_{i=1}^t Y_i^2 - \left ( \frac{1}{t}\sum_{i=1}^t Y_i \right )^2$.
  To this end, notice that by the independent SLLN (see \citet[\S IX, Theorem 12]{petrov2022sums}), we have that the left-most inequality of \eqref{eq:lyapunov-asympcs-proof-regularity} implies that
  \begin{equation}\label{eq:lyapunov-proof-second-moment-estimation}
    \frac{1}{t}\sum_{i=1}^t \left [Y_i^2 - \EE (Y_i^2) \right ] = o(\widetilde \sigma_t^2).
  \end{equation}
  Writing out the difference $\widehat \sigma_t^2 - \widetilde \sigma_t^2$, we see that
  \begin{align}
    \widehat \sigma_t^2 - \widetilde \sigma_t^2 &\equiv \frac{1}{t}\sum_{i=1}^t Y_i^2 - \left ( \frac{1}{t}\sum_{i=1}^t Y_i \right )^2 - \frac{1}{t}\sum_{i=1}^t \left [ \EE(Y_i^2) -  (\EE Y_i)^2 \right ]  \\
                                                &= \underbrace{\frac{1}{t}\sum_{i=1}^t \left [ Y_i^2 - \EE(Y_i^2) \right ]}_{=o(\widetilde \sigma_t^2) \text{ by } \eqref{eq:lyapunov-proof-second-moment-estimation}} - \underbrace{\left [  \left ( \frac{1}{t}\sum_{i=1}^t Y_i \right )^2 -  \frac{1}{t}\sum_{i=1}^t (\EE Y_i)^2  \right ]}_{(\star)}.
  \end{align}
  Focusing now only $(\star)$, we use the Lyapunov-type condition assumption $\sum_{i=1}^\infty \frac{\EE|Y_i - \mu_i|^{2+\delta}}{\sqrt{V_i}^{2+\delta}} < \infty$ combined with \citet[\S IX, Theorem 12]{petrov2022sums} to note that $\frac{1}{t}\sum_{i=1}^t (Y_i - \mu_i) = o(\sqrt{V_t} / t)$ and hence
  \begin{align}
    (\star) &\equiv \left ( \widetilde \mu_t + o(\sqrt{V_t} / t) \right )^2 - \frac{1}{t} \sum_{i=1}^t (\EE Y_i)^2 \\
&= \widetilde \mu_t^2 + o( \widetilde \mu_t \sqrt{V_t} / t) + o(V_t / t^2) - \frac{1}{t} \sum_{i=1}^t (\EE Y_i)^2\\
            &= \left ( \frac{1}{t}\sum_{i=1}^t \mu_i  \right )^2 + o( \widetilde \mu_t \sqrt{V_t} / t) + o(V_t / t^2) - \frac{1}{t} \sum_{i=1}^t (\EE Y_i)^2 \\
    &= - \frac{1}{t}\sum_{i=1}^t \left ( \mu_i - \widetilde \mu_t \right )^2 + o( \widetilde \mu_t \sqrt{V_t} / t) + o(V_t / t^2) \\
    &= o(\widetilde \sigma_t^2), 
  \end{align}
  where the final line follows from the assumptions of \eqref{eq:lyapunov-asympcs-proof-regularity}. This completes the proof of \cref{corollary:lyapunov-asympcs}.
  
\end{proof}

\subsection{Proof of Theorems~\ref{theorem:csate_randomized} and~\ref{theorem:csate_observational}}
\label{proof:csate_randomized}

In the proofs that follow, we will make extensive use of some convenient notation, namely the sample average operator $\PP_tf(Z) \equiv \frac{1}{t} \sum_{i=1}^t f(Z_i)$ and the conditional expectation operator $\PP \widehat f(Z) \equiv \PP(\widehat f(Z_i) \mid Z_1', \dots, Z_n')$ where $Z_1', \dots, Z_n'$ are the data used to construct $\widehat f$.

First, let us analyze the almost-sure behavior of the AIPW estimator $\widehat \psi_t$ for the average treatment effect $\psi$.

\begin{lemma}[Decomposition of $\widehat \psi_t - \psi$]
\label{lemma:estimatorDecomp}
Let $\widehat \psi_t := \PP_\Teval(\widehat f_{\Ttrain}) = \frac{1}{\Teval} \sum_{i=1}^\Teval \widehat f_{\Ttrain}(Z_i^\eval)$ be a (possibly misspecified) estimator of $\psi := \PP(f) = \EE\{ f(Z^\eval) \}$ based on $(Z^\eval_1, \dots, Z^\eval_\Teval)$ where $\widehat f_{\Ttrain}$ can be any estimator built from $(Z^\train_1, \dots, Z^\train_{\Ttrain})$ and $f : \mathcal Z \rightarrow \mathbb R$ any function. Furthermore, assume that there exists $\widebar f$ such that $\| \widehat f_{\Ttrain} - \widebar f\|_{L_2(\PP)} \rightarrow 0$. In other words, $\widehat f_{\Ttrain}$ is an estimator of $f$ but may instead converge to $\widebar f$. Then we have the decomposition, 
\[ \widehat \psi_t - \psi = \SA_t + \EP_t + \Bias_t \]
where 
\begin{alignat}{3}
    &\SA_t& :=\ & (\PP_\Teval - \PP)\widebar f \quad &\text{is the centered sample average term,}& \\
    &\EP_t& :=\ & (\PP_\Teval - \PP)(\widehat f - \widebar f) \quad &\text{is the empirical process term, and}&\\
    &\Bias_t& :=\ & \PP(\widehat f - f) \quad &\text{is the bias term.}&
\end{alignat}
\end{lemma}
\begin{proof}
By definition of the quantities involved, we decompose
    \begin{align}
        \widehat \psi_t - \psi &= \PP_\Teval (\widehat f_{\Ttrain}) - \PP(f) \\
        &= (\PP_\Teval - \PP) (\widehat f_{\Ttrain}) + \PP(\widehat f_{\Ttrain} - f ) \\
        &= \underbrace{(\PP_\Teval - \PP) (\widehat f_{\Ttrain} - \widebar f)}_{\EP_t} + \underbrace{(\PP_\Teval - \PP)\widebar f}_{\SA_t} + \underbrace{\PP(\widehat f_{\Ttrain} - f )}_{\Bias_t},
    \end{align}
    which completes the proof.
\end{proof}

Now, let us analyze the almost-sure behaviour of the empirical process term $\EP_t$ and the bias term $\Bias_t$ to show that they vanish asymptotically at sufficiently fast rates. First, let us examine $\EP_t$.
\begin{lemma}[Almost sure convergence of $\EP_t$]
\label{lemma:as_converge_EP}
Let $\PP_\Teval$ denote the empirical measure over $\mathbf{Z}^\eval_\Teval := (Z^\eval_1, \dots, Z_\Teval^\eval)$ and let $\widehat f_{\Ttrain}(z)$ be any function estimated from a sample $\splitSet^\train_{\Ttrain} = (Z^\train_1, Z^\train_2, \dots, Z^\train_{\Ttrain})$ which is independent of $\splitSet^\eval_\Teval$. If $\widehat \pi_t \in [\delta, 1-\delta]$ almost surely, then,
\[ \EP_t := (\PP_\Teval - \PP) (\widehat f_{\Ttrain} - \widebar f) = O \left ( \left \{  \sum_{a=0}^1\| \widehat \mu^a_t - \widebar \mu^a \|_{L_2(\PP)} \right \} \sqrt{ \frac{\log \log t}{t}}\right ).\]
In particular, if $\|\widehat \mu^a_t - \widebar \mu^a \|_{L_2(\PP)} = o(1)$ for each $a$, then we have that $\EP_t$ almost-surely converges to 0 at a rate of $o(\sqrt{\log \log t / t})$, but possibly faster.
\end{lemma}

The proof proceeds in two steps. First, we use an argument from \citet{kennedy2020sharp} and the law of the iterated logarithm to show $\EP_t = O \left ( \| \widehat f_t - \widebar f \| \sqrt{\log \log t / t} \right )$. Second and finally, we upper bound $\|\widehat f_t - \widebar f\|$ by $O\left ( \sum_{a = 0}^1 \| \widehat \mu_t^a - \widebar \mu \| \right )$.
\begin{proof}
\textbf{Step 1.}
Following the proof of \citet[Lemma 2]{kennedy2020sharp}, we have that conditional on $\splitSet^\train_\infty := (Z_t^\train)_{t=1}^\infty$ and $\whichSplitSet^\train_\infty := (\Ind(Z_t \in \splitSet_\infty^\train))_{t=1}^\infty$ the term of interest has mean zero:
\[ \EE \left \{ \PP_\Teval(\widehat f_{\Ttrain} - \widebar f) \bigm | \splitSet^\train_\infty, \whichSplitSet_\infty^\train \right \} = \EE (\widehat f_{\Ttrain} - \widebar f \mid \splitSet^\train_{\infty}, \whichSplitSet^\train_\infty) = \PP(\widehat f_{\Ttrain} - \widebar f).\] 
Now, we upper bound the conditional variance of a single summand,
\begin{align}
    \Var\left \{ (1-\PP)(\widehat f_{\Ttrain} - \widebar f) \bigm | \splitSet^\train_\infty, \whichSplitSet_\infty^\train \right \} &= \Var\left \{ (\widehat f_{\Ttrain} - \widebar f) \bigm | \splitSet^\train_\infty, \whichSplitSet_\infty^\train \right \}\\
    &\leq \| \widehat f_{\Ttrain} - \widebar f \|^2.
\end{align}
In particular, this means that
\[ \left (\frac{\Teval (\PP_\Teval - \PP)(\widehat f_{\Ttrain} - \widebar f)}{\| \widehat f_{\Ttrain} - \widebar f \|} \Bigm | \splitSet_\infty^\train,  \whichSplitSet_\infty^\train \right ) \]
is a sum of \iid{} random variables with conditional mean zero and conditional variance at most $1$, and thus by the law of the iterated logarithm,
\[ \PP \left ( \limsup_{t \rightarrow \infty} \frac{\pm \sqrt \Teval  (\PP_\Teval - \PP)(\widehat f_{\Ttrain} - \widebar f)}{\sqrt{2 \log \log \Teval}\| \widehat f_{\Ttrain} - \widebar f \|} \leq 1 \Bigm | \splitSet_\infty^\train, \whichSplitSet_\infty^\train \right ) = 1. \]
Therefore, we have that 
\begin{align}
    &\PP\left ( \frac{(\PP_\Teval - \PP)(\widehat f_{\Ttrain} - \widebar f)}{\|\widehat f_t - \widebar f \| \sqrt{ \log \log t / t}} = O(1) \Bigm | \splitSet_\infty^\train, \whichSplitSet_\infty^\train \right )\\
    =\ &\PP\left ( \frac{\| \widehat f_{\Ttrain} - \widebar f\|}{\|\widehat f_t - \widebar f \|} \cdot \frac{\sqrt t(\PP_\Teval - \PP)(\widehat f_{\Ttrain} - \widebar f)}{ \sqrt{ \log \log t} \| \widehat f_{\Ttrain} - \widebar f \|} = O(1) \Bigm | \splitSet_\infty^\train, \whichSplitSet_\infty^\train \right )\\
    =\ &\PP\left ( \lim_{t \rightarrow \infty} \underbrace{\frac{\| \widehat f_{\Ttrain} - \widebar f \|}{\|\widehat f_t - \widebar f \|}}_{O(1)} \cdot \underbrace{\frac{\sqrt{t \log \log \Teval} }{\sqrt{\Teval \log \log t}}}_{O(1)} \cdot \underbrace{\frac{\sqrt \Teval (\PP_\Teval - \PP)(\widehat f_{\Ttrain} - \widebar f)}{\sqrt{ \log \log \Teval}\| \widehat f_{\Ttrain} - \widebar f \|}}_{O(1)} = O(1) \Biggm | \splitSet_\infty^\train, \whichSplitSet_\infty^\train \right ) = 1.\\
\end{align}
Finally, by iterated expectation,
\begin{align} 
\PP \left ( \frac{(\PP_\Teval - \PP)(\widehat f_{\Ttrain} - \widebar f)}{\|\widehat f_t - \widebar f \| \sqrt{ \log \log t / t}} = O(1) \right ) &= \EE \left [\PP \left ( \frac{(\PP_\Teval - \PP)(\widehat f_{\Ttrain} - \widebar f)}{\|\widehat f_t - \widebar f \| \sqrt{ \log \log t / t}} = O(1) \Biggm | \splitSet_\infty^\train, \whichSplitSet_\infty^\train \right ) \right ]\\
&= \EE 1 = 1,
\end{align}
which completes Step 1.
\\

\noindent \textbf{Step 2.} Now, let us upper bound $\| \widehat f_t - \widebar f \|$ by $O \left ( \sum_{a=0}^1 \| \widehat \mu_t^a - \widebar \mu^a \| \right )$. To simplify the calculations which follow, define
\[\widehat f^1(Z_i) := \widehat \mu^1(X_i) +  \frac{A_i}{\widehat \pi(X_i)} \left \{ Y_i - \widehat \mu^{A_i}(X_i) \right \} ~~ \text{and} ~~ \widebar f^1(Z_i) := \widebar \mu^1(X_i) +  \frac{A_i}{\pi(X_i)} \left \{ Y_i - \widebar \mu^1 (X_i)\right \}. \]
Analogously define $\widehat f^0$ and $\widebar f^0$ so that $\widehat f = \widehat f^1 - \widehat f^0$ and $\widebar f = \widebar f^1 - \widebar f^0$. Writing out $\| \widehat f^1_t - \widebar f^1\|$,
\begin{align}
  \| \widehat f_t^1 - \widebar f^1\| &= \left \| \widehat \mu^1_t + \frac{A}{\widehat \pi_t} \left \{ Y - \widehat \mu^1_t \right \} - \widebar \mu^1_t - \frac{A}{\widehat \pi_t} \left \{ Y - \widebar \mu^1 \right \} \right \| \\
                                     &= \left \| \left \{ \widehat \mu^1_t - \widebar \mu^1 \right \} \left \{ 1 - \frac{A}{\widehat \pi_t} \right \} \right \| \\
                                     &\leq \frac{1}{\delta}\cdot \left \| \widehat \mu^1_t - \widebar \mu^1 \right \| = O\left ( \| \widehat \mu_t^1 - \widebar \mu^1 \| \right ),
\end{align}
where the last inequality follows from the assumed bounds on $\widehat \pi(X)$. A similar story holds for $\| \widehat f^0_t - \widebar f^0 \|$, and hence by the triangle inequality,
\[ \|\widehat f_t - \widebar f\| = O \left ( \sum_{a=0}^1 \| \widehat \mu_t^a - \widebar \mu^a \| \right ), \]
which completes the proof.
\end{proof}
Now, we examine the asymptotic almost-sure behaviour of the bias term, $\Bias_t$ by upper-bounding this term by a product of $L_2(\PP)$ estimation errors of nuisance functions.
\begin{lemma}[Almost-surely bounding $\Bias_\Teval$ by $L_2(\PP)$ errors of nuisance functions]
\label{lemma:boundingBias}
Suppose $\widehat \pi_t \in [\delta, 1-\delta]$ almost surely for some $\delta > 0$. Then
\[ \Bias_\Teval = O\left ( \|\widehat \pi_{t} - \pi \|_{L_2(\PP)} \Big \{ \| \widehat \mu_{t}^1 - \mu^1 \|_{L_2(\PP)}  + \| \widehat \mu_{t}^0 - \mu^0 \|_{L_2(\PP)}\Big \} \right )\]
\end{lemma}
This is an immediate consequence of the usual proof for $O_{\PP}$ combined with the fact that expectations are real numbers, and thus stochastic boundedness is equivalent to almost-sure boundedness. For completeness, we recall this proof here as it is short and illustrative.
\begin{proof}
To simplify the calculations which follow, define
\[\widehat f^1(Z_i) := \widehat \mu^1(X_i) +  \frac{A_i}{\widehat \pi(X_i)} \left \{ Y_i - \widehat \mu^{A_i}(X_i) \right \} ~~ \text{and} ~~ f^1(Z_i) := \mu^1(X_i) +  \frac{A_i}{\pi(X_i)} \left \{ Y_i - \mu^1 (X_i)\right \}. \]
Analogously define $\widehat f^0$ and $f^0$ so that $\widehat f = \widehat f^1 - \widehat f^0$ and $f = f^1 - f^0$.
Therefore,
\begin{align}
    \PP\left ( \widehat f^1 - f^1 \right ) &\overset{(i)}{=} \PP\left ( \frac{A}{\widehat \pi}(Y - \widehat \mu^A) + \widehat \mu^1 - \mu^1 \right ) \\
    &\overset{(ii)}{=} \PP\left [ \left ( \frac{\pi}{\widehat \pi} - 1 \right )\left ( \widehat \mu^1 - \mu^1 \right ) \right ] \\
    &\overset{(iii)}{\leq} \frac{1}{\delta}\PP\left ( \left | (\widehat \pi - \pi) (\widehat \mu^1 - \mu^1) \right | \right ) \\
    &\overset{(iv)}{\leq} \frac{1}{\delta} \|\widehat \pi - \pi\|_{L_2(\PP)} \|\widehat \mu^1 - \mu^1\|_{L_2(\PP)},
\end{align}
where $(i)$ and $(ii)$ follow by iterated expectation, $(iii)$ follows from the assumed bounds on $\widehat \pi$, and $(iv)$ by H\"older's inequality. Similarly, we have
\[ \PP(\widehat f^0 - f^0) \leq \frac{1}{1-\delta} \|\widehat \pi - \pi\| \|\widehat \mu^0 - \mu^0 \|. \]
Finally by the triangle inequality, 
\[ \PP(\widehat f - f) = O\left ( \left \|\widehat \pi - \pi \right \| \sum_{a = 0}^1  \left \|\widehat \mu^a - \mu^a \right \| \right ), \]
which completes the proof.
\end{proof}

\begin{lemma}[Almost-sure consistency of the influence function variance estimator]
\label{lemma:InfluenceFunctionVarianceEstimator}
    Suppose that $\| \widehat f_{\Ttrain} - \widebar f \|_2 = o(1)$ and that $\widebar f(Z)$ has a finite second moment. Then, \[\widehat \Var_{\Teval} (\widehat f_{\Ttrain}) = \Var(\widebar f) + o(1).\]
\end{lemma}

\begin{proof}
By direct calculation, we see that
  \begin{align}
    \PP_\Teval \left (\widehat f_\Ttrain - \PP_\Teval \widehat f_\Ttrain \right )^2 &= \underbrace{\PP_\Teval \left ( \widehat f_\Ttrain - \PP \widebar f \right )^2}_{(i)} + \underbrace{\left [ \PP_\Teval (\widehat f_\Ttrain - \PP \widebar f) \right ]^2}_{(ii)},
  \end{align}
  and we have by the preceding lemmas that $(ii) \xrightarrow{\as} 0$ since $a_t \xrightarrow{\as} 0 \implies a_t^2 \xrightarrow{\as}0 $ and thus it suffices to show that $(i) \xrightarrow{\as} \Var(\widebar f)$. Indeed, we use the inequality $2ab \leq |ab| \leq a^2 + b^2$ and notice that
  \begin{align}
    (i) &= \PP_\Teval \left ( \widehat f_\Ttrain - \widebar f + \widebar f - \PP \widebar f \right )^2\\
    &= \PP_\Teval \left ( \widehat f_\Ttrain - \widebar f \right )^2 + 2\PP_\Teval \left ( \widehat f_\Ttrain - \widebar f \right ) \left ( \widebar f - \PP \widebar f \right ) + \PP_\Teval \left ( \widebar f - \PP \widebar f \right )^2\\
    &\leq \PP_\Teval \left ( \widehat f_\Ttrain - \widebar f \right )^2 + \PP_\Teval \left ( \widehat f_\Ttrain - \widebar f \right )^2 + \PP_\Teval \left ( \widebar f - \PP \widebar f \right )^2 + \PP_\Teval \left ( \widebar f - \PP \widebar f \right )^2\\
    &= 2 \Bigg [\underbrace{\PP_\Teval \left ( \widehat f_\Ttrain - \widebar f \right )^2}_{(i.i)} + \underbrace{\PP_\Teval \left ( \widebar f - \PP \widebar f \right )^2}_{(i.ii)} \Bigg ]
  \end{align}
  Let us now study $(i.i)$ and $(i.ii)$ separately. Note that $\left \{ (\widehat f_\Ttrain(Z_i) - \widebar f(Z_i))^2 \right \}_{i \in \NN}$ are conditionally \iid{} and hence $\PP_\Teval \left ( \widehat f_\Ttrain - \widebar f \right )^2$ is a (conditional) reverse submartingale. By \citet[Theorem 2]{manole2023martingale}, we have
  \begin{equation}
    \PP \left ( \exists k \geq \Ttrain : \PP_\Teval (\widehat f_\Ttrain - \widebar f)^2 \geq \eps \right ) \leq \PP (\widehat f_\Ttrain - \widebar f)^2 = \| \widehat f_\Ttrain - \widebar f\|^2 \xrightarrow \as 0,
  \end{equation}
  so that $(i.i) \to 0$ $\PP$-almost surely. Furthermore, we have by the finite second moment assumption $\PP |\widebar f - \PP \widebar f|^2 < \infty$ that $(i.ii) \to \Var(\widebar f)$ $\PP$-almost surely by the strong law of large numbers. Putting these together, we have that $(i) \to \Var(\widebar f)$ $\PP$-almost surely, completing the proof.  
\end{proof}

\begin{proposition}[General \asympcs{}s under sequential cross fitting]
  \label{proposition:general-asympcs-cross-fit}
    Consider the cross-fit estimator as defined in \eqref{eq:cross_fit_estimator}:
    \begin{equation}
      \widehat \psi_t^\times := \frac{\sum_{i=1}^\Teval f_{\Ttrain}(Z_i^\eval) + \sum_{i=1}^{\Ttrain} f_{\Teval}(Z_i^\train)}{t},
    \end{equation}
    and the cross-fit variance estimator as defined in \eqref{eq:cross_fit_variance}:
    \begin{equation}
      \label{eq:csate_randomized_variance}
      \widehat \Var_t(f) := \frac{\widehat \Var_\Teval(\widehat f_{\Ttrain}) + \widehat \Var_{\Ttrain}(\widehat f_{\Teval})}{2}.
    \end{equation}
    Suppose that $\Bias_t$ and $\EP_t$ are both $o(\sqrt{\log \log t / t})$. Then,
    \[ \widehat \psi_t^\times \pm \sqrt{\widehat{\Var}_t(\widehat f) }\sqrt{\frac{2(t \rho^2 + 1)}{t^2 \rho^2} \log \left ( \frac{\sqrt{t \rho^2 + 1}}{\alpha} \right )} \]
    forms a $(1-\alpha)$-\asympcs{} for $\psi$.
\end{proposition}
\begin{proof}
  Writing out the centered cross-fit estimator $\widehat \psi_t^\times - \psi$ using the decomposition of Lemma~\ref{lemma:estimatorDecomp}, we have
  \begin{align}
    \widehat \psi_t^\times - \psi &= \frac{\sum_{i=1}^\Teval \widehat f_\Ttrain(Z_i^\eval) + \sum_{i=1}^\Ttrain \widehat f_\Teval(Z_i^\train) - t \psi}{t} \nonumber \\
                                  &= \frac{\sum_{i=1}^\Teval(f_\Ttrain(Z_i^\eval) - \psi) + \sum_{i=1}^\Ttrain (\widehat f_\Teval(Z_i^\train) - \psi)}{t} \nonumber \\
                                  &= \frac{( \Teval \SA_{t, \eval} + \Ttrain \SA_{t,\train}) + \Teval \EP_{t, \eval} + \Ttrain \EP_{t, \train} +  \Teval \Bias_{t, \eval} + \Ttrain \Bias_{t, \train}  }{t}\nonumber \\
                                  &= (\PP_t - \PP)\widebar f(Z) + \frac{\Teval \EP_{t, \eval} + \Ttrain \EP_{t, \train} +  \Teval \Bias_{t, \eval} + \Ttrain \Bias_{t, \train} }{t}\nonumber \\
                                  &= (\PP_t - \PP)\widebar f(Z) + \underbrace{O(\Bias_t + \EP_{t} )}_{o(\sqrt{\log \log t / t})} \label{eq:cross-fit-decomp}
  \end{align}
  where $\EP_{t, \eval} := \frac{1}{\Teval}\sum_{i=1}^\Teval \widehat f_\Ttrain(Z_i^\eval)$ and $\EP_{t,\train} := \frac{1}{\Ttrain}\sum_{i=1}^\Ttrain \widehat f_\Teval(Z_i^\train)$, and similarly for $\SA_{t, \eval}, \SA_{t,\train}, \Bias_{t,\eval}$, and $\Bias_{t,\train}$. Applying the proof of Theorem~\ref{theorem:acs} (but with variance consistency $\widehat \Var_t(\widehat f) \xrightarrow{\as} \Var(\widebar f)$ obtained via Lemma~\ref{lemma:InfluenceFunctionVarianceEstimator}), we have that
  \[ \widehat \psi_t^\times \pm \sqrt{\widehat{\Var}_t(\widehat f) }\sqrt{\frac{2(t \rho^2 + 1)}{t^2 \rho^2} \log \left ( \frac{\sqrt{t \rho^2 + 1}}{\alpha} \right )} + o(\sqrt{\log \log t / t}) \]
  forms a nonasymptotic $(1-\alpha)$-\cs{} for $\psi$. Consequently,
  \[ \widehat \psi_t^\times \pm \sqrt{\widehat{\Var}_t(\widehat f) }\sqrt{\frac{2(t \rho^2 + 1)}{t^2 \rho^2} \log \left ( \frac{\sqrt{t \rho^2 + 1}}{\alpha} \right )}\]
  forms a $(1-\alpha)$-\asympcs{} for $\psi$ with rate $o(\sqrt{\log \log t / t})$ which completes the proof.
\end{proof}

\subsubsection{Proof of Theorem~\ref{theorem:csate_randomized}}
\begin{proof}

When propensity scores are known, we have that $\Bias_t \eqas 0$ by Lemma~\ref{lemma:boundingBias}. By assumption, $\EE \| \widehat \mu_{\Ttrain}^a(X) - \widebar \mu^a(X) \|_2 = o(1)$, and thus by Lemma~\ref{lemma:as_converge_EP}, $\EP_t = o\left ( \sqrt{ \log \log t / t} \right )$. Combining these conditions on $\Bias_t$ and $\EP_t$ with Proposition~\ref{proposition:general-asympcs-cross-fit}, we obtain the desired result. This completes the proof of Theorem~\ref{theorem:csate_randomized}.
\end{proof}

\subsubsection{Proof of Theorem~\ref{theorem:csate_observational}}
\label{proof:csate_observational}
\begin{proof}
  By Lemmas~\ref{lemma:as_converge_EP} and \ref{lemma:boundingBias} we have that both $\Bias_t$ and $\EP_t$ are $o(\sqrt{\log \log t / t})$. Applying Proposition~\ref{proposition:general-asympcs-cross-fit}, we obtain the desired result, completing the proof of Theorem~\ref{theorem:csate_observational}.
\end{proof}

\subsection{Proof of \cref{theorem:lyapunov_ate_randomized}}
\label{proof:lyapunov_ate}

\begin{lemma}[Decomposition of $t \widehat \psi_t^\times - t\widetilde \psi_t$]
\label{lemma:estimatorDecomp-time-varying}
Let $\widehat \psi_t^\times$ be as in \eqref{eq:cross_fit_estimator}. Furthermore, assume that there exists $\widebar f$ such that $\| \widehat f_{t} - \widebar f\|_{L_2(\PP)} \rightarrow 0$. In other words, $\widehat f_{t}$ is an estimator of $f$ but may instead converge to $\widebar f$. Then we have the decomposition,
\begin{equation}
  \label{eq:decomposition-time-varying}
  t\widehat \psi_t^\times - t\widetilde \psi_t = \SAtilde_t + \underbrace{\EPtilde_{t, \eval} + \EPtilde_{t,\train}}_{\EPtilde_t} + \underbrace{\Biastilde_{t, \eval} + \Biastilde_{t, \train}}_{\Biastilde_t}
\end{equation}
where
\begin{alignat}{3}
    &\SAtilde_{t}& :=\ & \sum_{i=1}^t \left [\widebar f(Z_i) - \PP(\widebar f(Z_i)) \right ], \\
    &\EPtilde_{t, \eval}& :=\ &  \sum_{i=1}^{\Teval}\left \{\left [ \widehat f_{\Ttrain}(Z_i^\eval) - \PP (\widehat f_{\Ttrain}(Z_i^\eval)) \right ] - \left [ \widebar f(Z_i^\eval) - \PP(\widebar f(Z_i^\eval)) \right ] \right \}, \\
    & \EPtilde_{t, \train} & :=\ &  \sum_{i=1}^{\Ttrain}\left \{\left [ \widehat f_{\Teval}(Z_i^\train) - \PP( \widehat f_{\Teval}(Z_i^\train) )\right ] - \left [ \widebar f(Z_i^\train) - \PP(\widebar f(Z_i^\train)) \right ] \right \}, \\
    &\Biastilde_{t,\eval}& :=\ & \sum_{i=1}^\Teval \PP(\widehat f_{\Ttrain}(Z_i^\eval) - f(Z_i^\eval)), ~~~\text{and} \\
    &\Biastilde_{t, \train}& :=\ & \sum_{i=1}^{\Ttrain} \PP(\widehat f_{\Teval}(Z_i^\train) - f(Z_i^\train)).
\end{alignat}
\end{lemma}
\begin{proof}
  First, note that $t\widehat \psi_t^\times - t\widetilde \psi_t$ can be written as
  \begin{align}
    t\widehat \psi_t^\times - t\widetilde \psi_t &= \sum_{i=1}^{\Teval} \widehat f_{\Ttrain}(Z_i^\eval) + \sum_{i=1}^{\Ttrain} \widehat f_{\Teval}(Z_i^\train) - \sum_{i=1}^{\Teval} \PP f(Z_i^\eval) - \sum_{i=1}^{\Ttrain} \PP f(Z_i^\train)\\
    &= \underbrace{\sum_{i=1}^{\Teval} \left [ \widehat f_{\Ttrain}(Z_i^\eval) - \PP f(Z_i^\eval) \right]}_{(i)}+ \underbrace{\sum_{i=1}^{\Ttrain} \left [ \widehat f_{\Teval}(Z_i^\train) - \PP f(Z_i^\train) \right ]}_{(ii)}.
  \end{align}
  We will handle each sum separately and then combine them to arrive at the final decomposition \eqref{eq:decomposition-time-varying}. Taking a closer look at $(i)$ first, we have
 \begin{alignat}{2}
   (i) &= &\sum_{i=1}^{\Teval} \left \{ \left [ \widehat f_{\Ttrain}(Z_i^\eval) - \PP ( \widehat f_{\Ttrain}(Z_i^\eval) )\right  ] - \left [ \widebar f(Z_i^\eval) - \PP(\widebar f(Z_i^\eval))\right ] \right \}\\
   & &+ \sum_{i=1}^{\Teval} \PP\left  \{ \widehat f_{\Ttrain}(Z_i^\eval) - f(Z_i^\eval) \right \} + \sum_{i=1}^{\Teval} \left \{\widebar f(Z_i^\eval) - \PP(\widebar f (Z_i^\eval)) \right \}.
 \end{alignat}
 Similarly for $(ii)$, we have
 \begin{alignat}{2}
   (ii) &= &\sum_{i=1}^{\Ttrain} \left \{ \left [ \widehat f_{\Teval}(Z_i^\train) - \PP ( \widehat f_{\Teval}(Z_i^\train) )\right  ] - \left [ \widebar f(Z_i^\train) - \PP(\widebar f(Z_i^\train))\right ] \right \}\\
   & &+ \sum_{i=1}^{\Ttrain} \PP\left  \{ \widehat f_{\Teval}(Z_i^\train) - f(Z_i^\train) \right \} + \sum_{i=1}^{\Ttrain} \left \{\widebar f(Z_i^\train) - \PP(\widebar f (Z_i^\train)) \right \}.
 \end{alignat}
 Putting $(i)$ and $(ii)$ together, we have
 \begin{align}
   t\widehat \psi_t^\times - t \widetilde \psi_t &= \sum_{i=1}^{\Ttrain} \left \{ \widebar f(Z_i^\eval) - \PP(\widebar f(Z_i^\eval)) \right \} + \sum_{i=1}^\Teval \left \{ \widebar f(Z_i^\train) - \PP(\widebar f(Z_i^\train)) \right \} + \EPtilde_t + \Biastilde_t \\
   &= \underbrace{\sum_{i=1}^t \left \{ \widebar f(Z_i) - \PP(\widebar f(Z_i)) \right \}}_{\SAtilde_t} + \EPtilde_t + \Biastilde_t,
 \end{align}
    which completes the proof.
\end{proof}

\begin{lemma}[Almost sure behavior of $\EPtilde_t$]
\label{lemma:as_converge_EPtilde}
Suppose that there exists $\delta > 0$ such that $\widehat \pi_t \in [\delta, 1-\delta]$ almost surely for all $t$. Then,
\begin{equation}
  \label{eq:EP-time-varying-bigO}
  \EPtilde_t = o \left ( \left \{ \sum_{a=0}^1 \sup_{i} \| \widehat \mu_t^a(X_i) - \widebar \mu^a(X_i) \|_\LP \right \} \sqrt{t \log \log t}\right ).
\end{equation}

\end{lemma}

\begin{proof}
We will show that the result \eqref{eq:EP-time-varying-bigO} holds for each of $\EPtilde_{t, \eval}$ and $\EPtilde_{t, \train}$, thereby yielding the same result for their sum $\EPtilde_t$. The proof proceeds in two steps. First, we use an argument from \citet{kennedy2020sharp} and the law of the iterated logarithm to bound $\SAtilde_t$ in terms of $\sup_i \|\widehat f_{\Teval}(Z_i) - \widebar f(Z_i) \|$. Second and finally, we upper bound $\sup_i \|\widehat f_t(Z_i) - \widebar f(Z_i)\|$ by $O\left ( \sum_{a = 0}^1 \sup_i \| \widehat \mu_t^a(Z_i) - \widebar \mu(Z_i) \| \right )$.

\paragraph*{Step 1}
Let us first consider $\EPtilde_{t, \eval}$. Following the proof of \citet[Lemma 2]{kennedy2020sharp} and of Lemma~\ref{lemma:as_converge_EP}, note that conditional on $\splitSet^\train_\infty := (Z_t^\train)_{t=1}^\infty$ and $\whichSplitSet^\train_\infty := (\Ind(Z_t \in \splitSet_\infty^\train))_{t=1}^\infty$ the summands of $\EPtilde_{t, \eval}$ have mean zero:
\[ \PP \left \{\left [\widehat f_{\Ttrain}(Z_i^\eval) - \widebar f(Z_i^\eval) \right ] - \left [ \widebar f(Z_i^\eval) - \PP(\widebar f(Z_i^\eval)) \right ] \mid \splitSetInfty^\train, \whichSplitSet_\infty^\train \right \} = 0. \]
Similar to the proof of Lemma~\ref{lemma:as_converge_EP}, we upper bound the conditional variance of a single summand,
\begin{align}
    \Var\left \{ (1-\PP)(\widehat f_{\Ttrain}(Z_i^\eval) - \widebar f(Z_i^\eval)) \bigm | \splitSet^\train_\infty, \whichSplitSet_\infty^\train \right \} &= \Var\left \{ (\widehat f_{\Ttrain}(Z_i^\eval) - \widebar f(Z_i^\eval)) \bigm | \splitSet^\train_\infty, \whichSplitSet_\infty^\train \right \}\\
    &\leq \| \widehat f_{\Ttrain}(Z_i^\eval) - \widebar f(Z_i^\eval) \|^2_\LP.
\end{align}
Denote the following process $\nu_{t, \eval}$ as the supremum of the above with respect to $i \in \{1,2,\dots\}$:
\[ \nu_{t, \eval} := \sup_{1\leq i \leq \infty} \|\widehat f_\Ttrain (Z_i^\eval) - \widebar f(Z_i^\eval) \|_\LP. \]
Then we can lower bound the following conditional probability
\begin{align}
  &\PP \left ( \limsup_{t \to \infty} \frac{\pm \EPtilde_{t, \eval}}{\nu_{t, \eval}\sqrt{ 2 t \log \log t}} \leq 1 \Bigm \vert \splitSet_\infty^\train, \whichSplitSet_\infty^\train \right ) \nonumber \\
  =\ & \PP\left ( \limsup_{t\to \infty} \sum_{i=1}^{\Teval}\frac{ \pm \left \{\left [ \widehat f_{\Ttrain}(Z_i^\eval) - \PP (\widehat f_{\Ttrain}(Z_i^\eval)) \right ] - \left [ \widebar f(Z_i^\eval) - \PP(\widebar f(Z_i^\eval)) \right ] \right \}}{ \nu_{t, \eval}\sqrt{2 t \log \log t} } \leq 1 \Bigm \vert \splitSet_\infty^\train, \whichSplitSet_\infty^\train\right )\nonumber\\
  \geq\ & \PP\left ( \limsup_{t\to \infty} \sum_{i=1}^{\Teval}\frac{ \pm \left \{\left [ \widehat f_{\Ttrain}(Z_i^\eval) - \PP (\widehat f_{\Ttrain}(Z_i^\eval)) \right ] - \left [ \widebar f(Z_i^\eval) - \PP(\widebar f(Z_i^\eval)) \right ] \right \}}{\| \widehat f_\Ttrain(Z_i^\eval) - \widebar f (Z_i^\eval) \|_{L_2(\PP)} \sqrt{2 t \log \log t} } \leq 1 \Bigm \vert \splitSet_\infty^\train, \whichSplitSet_\infty^\train\right )\nonumber\\
  =\ & \PP\left ( \limsup_{t\to \infty} \sum_{i=1}^{\Teval}\frac{ \pm \zeta_i }{\sqrt{2 t \log \log t} } \leq 1 \Bigm \vert \splitSet_\infty^\train, \whichSplitSet_\infty^\train\right ), \label{eq:EP-time-varying-LIL-zeta}
\end{align}
where $\zeta_i$ are independent mean-zero random variables with variance at most one (conditional on $\splitSet_\infty^\train, \whichSplitSet_\infty^\train$). By the law of the iterated logarithm, we have that $\eqref{eq:EP-time-varying-LIL-zeta} = 1$. In particular, since this event happens with probability one conditionally, it also happens with probability one marginally. It follows that
\begin{equation}
  \EPtilde_{t, \eval} = O\left ( \sup_{i} \|\widehat f_t(Z_i) - \widebar f(Z_i) \|_{\LP} \sqrt{t \log \log t} \right ).
\end{equation}
Applying the same technique to $\EPtilde_{t,\train}$, we have that $\EPtilde_{t, \eval} = O\left ( \sup_{i} \|\widehat f_t(Z_i) - \widebar f(Z_i) \|_{\LP} \sqrt{t \log \log t} \right )$, and hence
\begin{equation}
  \label{eq:EPtilde-upper-bound-f}
  \EPtilde_t = O\left ( \sup_{i} \|\widehat f_t(Z_i) - \widebar f(Z_i) \|_{\LP} \sqrt{t \log \log t} \right ).
\end{equation}

\paragraph*{Step 2} Now, following the same technique as Step 2 in the proof of Lemma~\ref{lemma:as_converge_EP}, we have that
\begin{equation}
  \label{eq:EPtilde-proof-step2}\|\widehat f_{\Ttrain}(Z_i) - \widebar f(Z_i)\| = O \left ( \sum_{a=0}^1 \| \widehat \mu_{t}^a(X_i) - \widebar \mu^a_i(X_i) \| \right ).
\end{equation}
Combining \eqref{eq:EPtilde-upper-bound-f} and \eqref{eq:EPtilde-proof-step2}, we have the desired result,
\begin{equation}
  \EPtilde_t = O \left ( \left \{ \sum_{a=0}^1  \sup_{1\leq i \leq \infty} \| \widehat \mu_{t}^a(X_i) - \widebar \mu^a_i(X_i) \| \right \} \sqrt{t \log \log t} \right ),
\end{equation}
which completes the proof.
\end{proof}
Now, we examine the asymptotic almost-sure behaviour of the bias term, $\Bias_t$ by upper-bounding this term by a product of $L_2(\PP)$ estimation errors of nuisance functions.
\begin{lemma}[Almost-sure behavior of $\Biastilde_t$]
\label{lemma:boundingBiastilde}
Suppose $\widehat \pi_t \in [\delta, 1-\delta]$ for every $t$ almost surely for some $\delta > 0$. Then,
\begin{align}
  \Biastilde_t = O \left ( \sum_{i=1}^t \|\widehat \pi_t(X_i) - \pi(X_i) \|_{L_2(\PP)} \sum_{a=0}^1  \| \widehat \mu_{t}^a(X_i) - \mu^a(X_i) \|_{L_2(\PP)}  \right )
\end{align}

\end{lemma}
The proof proceeds similarly to that of Lemma~\ref{lemma:boundingBias} but with additional care given to the fact that observations are no longer identically distributed.

\begin{proof}
Similar to the proof of Lemma~\ref{lemma:as_converge_EPtilde}, we will first prove the result for $\Biastilde_{t, \eval}$, and the proof proceeds similarly for $\Biastilde_{t, \train}$, thereby yielding the desired result for $\Biastilde_t \equiv \Biastilde_{t,\eval} + \Biastilde_{t, \train}$.
Following the same technique as Lemma~\ref{lemma:boundingBias}, we have that
\[ \PP(\widehat f_\Ttrain(Z_i^\eval) - f(Z_i^\eval)) = O\left ( \left \|\widehat \pi_\Ttrain(X_i^\eval) - \pi(X_i^\eval) \right \| \sum_{a = 0}^1  \left \|\widehat \mu_\Ttrain^a(X_i^\eval) - \mu^a(X_i^\eval) \right \| \right ). \]
Putting the above term back into the sum $\Biastilde_{t, \eval}$, we have
\begin{align}
  \Biastilde_{t, \eval} &:= \sum_{i=1}^T \PP(\widehat f_\Ttrain(Z_i^\eval) - f(Z_i^\eval))\\
  &= O\left ( \sum_{i=1}^T \left \|\widehat \pi_\Ttrain(X_i^\eval) - \pi(X_i^\eval) \right \| \sum_{a = 0}^1  \left \|\widehat \mu_\Ttrain^a(X_i^\eval) - \mu^a(X_i^\eval) \right \| \right ).
\end{align}
Using a similar argument to bound $\Biastilde_{t, \train}$ and putting these together, we have the following bound for $\Biastilde_{t} = \Biastilde_{t, \eval} + \Biastilde_{t, \train}$,
\begin{equation}
  \Biastilde_t = O\left ( \sum_{i=1}^t \left \|\widehat \pi_t(X_i) - \pi(X_i) \right \| \sum_{a = 0}^1  \left \|\widehat \mu_t^a(X_i) - \mu^a(X_i) \right \| \right ),
\end{equation}
which completes the proof.
\end{proof}

\begin{proposition}[General \asympcs{}s for time-varying causal effects under sequential cross fitting]
  \label{proposition:general-asympcs-time-varying-cross-fit}
    Consider the cross-fit estimator as defined in \eqref{eq:cross_fit_estimator}:
    \begin{equation}
      \widehat \psi_t^\times := \frac{\sum_{i=1}^\Teval f_{\Ttrain}(Z_i^\eval) + \sum_{i=1}^{\Ttrain} f_{\Teval}(Z_i^\train)}{t},
    \end{equation}
    and suppose we have access to a variance estimator $\widehat \Var_t(\widehat f)$ such that
    \[ \widehat \Var_t(\widehat f) - \widetilde \Var (\widebar f) = o(1). \]
    Suppose that $\Biastilde_t$ and $\EPtilde_t$ are both $o(\sqrt{t \log \log t})$, and that the conditions of~\cref{corollary:lyapunov-asympcs} hold but with $\infseq Yt1$ replaced by $(\widebar f(Z_t))_{t=1}^\infty$. Then,
    \[
      \widehat \psi_t^\times \pm \sqrt{\frac{2(t \rho^2 \widehat \Var_t(\widebar f) + 1)}{t^2 \rho^2} \log \left ( \frac{\sqrt{t \rho^2\widehat \Var_t(\widebar f) + 1}}{\alpha} \right )}
    \]
    forms a $(1-\alpha)$-\asympcs{} for $\widetilde \psi_t := \frac{1}{t}\sum_{i=1}^t \psi_i$.
\end{proposition}
\begin{proof}
  Writing out the centered cross-fit estimator on the ``sum scale'' $t(\widehat \psi_t^\times - \widetilde \psi)$ using the decomposition of Lemma~\ref{lemma:estimatorDecomp-time-varying}, we have
  \begin{align}
    t(\widehat \psi_t^\times - \widetilde \psi_t) &= \SAtilde_t + \underbrace{\EPtilde_{t, \eval} + \EPtilde_{t,\train}}_{\EPtilde_t=o(\sqrt{t \log \log t})} + \underbrace{\Biastilde_{t, \eval} + \Biastilde_{t, \train}}_{\Biastilde_t=o(\sqrt{t \log \log t})}. \nonumber
  \end{align}
  Therefore, we have that
  \[ \widehat \psi_t - \widetilde \psi_t = \frac{1}{t} \sum_{i=1}^t (\widebar f(Z_i) - \psi_i) + o(\sqrt{\log \log t / t}). \]
  Applying~\cref{corollary:lyapunov-asympcs} to $(\widebar f(Z_t))_{t=1}^\infty$ above, we have that
  \[
    \widetilde C_t^{(\times\star)} := \widehat \psi_t^\times \pm \sqrt{\frac{2(t \rho^2 \widehat \Var_t(\widebar f) + 1)}{t^2 \rho^2} \log \left ( \frac{\sqrt{t \rho^2\widehat \Var_t(\widebar f) + 1}}{\alpha} \right )} + o(\sqrt{\log \log t/ t})
  \]
  forms a nonasymptotic $(1-\alpha)$-\cs{} for $\widetilde \psi_t := \frac{1}{t}\sum_{i=1}^t \psi_i$, meaning ${\PP \left (\exists t: \widetilde \psi_t \notin \widetilde C_t^{(\times \star)}\right )\leq \alpha}$. Consequently,
  \[
    \widehat \psi_t^\times \pm \sqrt{\frac{2(t \rho^2 \widehat \Var_t(\widebar f) + 1)}{t^2 \rho^2} \log \left ( \frac{\sqrt{t \rho^2\widehat \Var_t(\widebar f) + 1}}{\alpha} \right )}
  \]
  forms a $(1-\alpha)$-\asympcs{} for $\widetilde \psi_t$, which completes the proof.
\end{proof}

\subsubsection{Proof of Theorem~\ref{theorem:lyapunov_ate_randomized}}
\begin{proof}
By Lemma~\ref{lemma:as_converge_EPtilde} combined with Assumption~\ref{assumption:sup-regression-estimators-converge}, we have that $\EPtilde_t = o(\sqrt{t\log \log t})$. In a randomized experiment, Assumption~\ref{assumption:doubly-robust-average-bias-convergence} holds by design, and thus by Lemma~\ref{lemma:boundingBiastilde}, we have that $\Biastilde_t = o(\sqrt{t \log \log t})$. Invoking Proposition~\ref{proposition:general-asympcs-time-varying-cross-fit}, we obtain the desired result.
\end{proof}


\subsection{Proof of \cref{theorem:type-I-error}}\label{proof:type-I-error}

Before diving into the proof of \cref{theorem:type-I-error}, we will introduce some notation. 
Let $(S(v))_{v \in \RR^+}$ be the (continuous-time) process obtained via $S(V_t) = S_t$ for each $t \in \NN$ (with $S(0) = 0$) and remaining piecewise constant in between the integers. Let $\widehat S(v)$ be defined analogously but with $\widehat V_t$ instead of $V_t$.
Our first lemma will allow $\widehat S(v)$ to be written in terms of $S(v)$ up to a smaller term in the argument.

  \begin{lemma}\label{lemma:V_t-polynomial-rate}
    Suppose $\widehat V_t - V_t = o(V_t^{\eta})$ for some $0< \eta < 1$, or in other words, $\widehat \sigma_t^2$ is a consistent estimator for $\widetilde \sigma_t^2$ at any polynomial rate in $V_t$. Then,
    \begin{equation}
      \widehat S(v) = S(v + o(v^{\eta})) ~~\text{as $v \to \infty$.}
    \end{equation}
  \end{lemma}

  \begin{proof}
    The result follows from the fact that $\widehat V_t - V_t = o(V_t^{\eta})$ and that $V_t \to \infty$ as $t \to \infty$.
  \end{proof}

  The next lemma approximates $\widehat S(v)$ by a Wiener process using the strong approximation results of \citet{strassen1967and} combined with \cref{lemma:V_t-polynomial-rate}.
  \begin{lemma}\label{lemma:interpolated-strong-approx}
    After potentially enlarging the probability space, there exists a Wiener process $(W(v))_{v \in \RR^+}$ such that 
    \begin{equation}
      |\widehat S(v) - W(v)| = o(v^{\kappa/2})
    \end{equation}
    for some $0<\kappa < 1$.
  \end{lemma}

  \begin{proof} By direct calculation, we have
    \begin{align}
      |\widehat S(v) - W(v)| &\overset{(i)}{=} |S(v + o(v^{\eta})) - W(v)|\\
                             &= \underbrace{|S(v + o(v^{\eta})) - W(v + o(v^{\eta}))|}_{(\star)} + \underbrace{|W(v + o(v^{\eta})) - W(v)|}_{(\dagger)},
    \end{align}
    where $(i)$ follows from \cref{lemma:V_t-polynomial-rate}. We will subsequently bound $(\star)$ using \citet[Theorem 4.4]{strassen1967and}, and $(\dagger)$ via the Borel-Cantelli lemma combined with elementary properties of the Wiener process.
    \paragraph*{Bounding $(\star)$}
    By \citet[Theorem 4.4]{strassen1967and}, we have that on a sufficiently rich probability space, $S(v) - W(v) = o \left ( v^{3/8} \log v \right )$, and hence
    \begin{align}
      (\star) &:= |S(v + o(v^{\eta})) - W(v + o(v^{\eta}))| \\
              &= o \left ( (v+o(v^{\eta}))^{3/8} \log (v + o(v^{\eta})) \right ) \label{eq:strassen-type-I-error-proof}\\
      &= o \left ( v^{3/8} \log v \right ) \label{eq:approx-rate-without-eta},
    \end{align}
    where \eqref{eq:approx-rate-without-eta} follows from the assumption that $\eta < 1$.
    \paragraph*{Bounding $(\dagger)$}
    By elementary properties of the Wiener process, we have $W(v + o(v^{\eta})) - W(v) = \sqrt{o(v^{\eta})}Z$ where $Z \sim N(0, 1)$ is an independent standard Gaussian random variable.
    It then remains to show that $\sqrt{o(v^{\eta})} Z = o(v^{\kappa/2})$. Indeed, notice that
    \begin{equation}
      \PP \left (\frac{v^{\eta/2}}{v^{\kappa/2}} Z \geq \eps \right ) \leq \exp \left \{ -\frac{\eps^2 v^{\kappa} }{2  v^{\eta}} \right \} = \exp \left \{ -\frac{\eps^2 v^{\kappa-\eta}}{2} \right \},
    \end{equation}
    and hence by the Borel-Cantelli lemma,
    \begin{equation}
      \frac{v^{1/2q}}{v^{\kappa/2}} Z \to 0 ~~\text{almost surely.}
    \end{equation}
    for any $\eta < \kappa \leq 1$.
  \end{proof}

\begin{proof}[Proof of \cref{theorem:type-I-error}]
  Recall that $\rho_m = \sqrt{c/\widehat V_m d_m}$ where $d_m$ is any increasing sequence $\to \infty$ as $m \to \infty$. Define $(\widehat V'(t))_{t \in \RR^+}$ as the continuous-time process obtained by setting $\widehat V'(t) = \widehat V_t$ for each $t \in \NN$ and piecewise-constant between the integers. Writing down the limit supremum in \cref{theorem:type-I-error},
  \begin{align}
    & \limsup_{m \to \infty}\PP\left (\exists t \in \NN^{\geq m} : |S_t| \geq \sqrt{\frac{2(t\widehat \sigma_t^2\rho_m^2 + 1)}{\rho_m^2}\log \left ( \frac{\sqrt{t\widehat \sigma_t^2 \rho_m^2 +1}}{\alpha} \right )} \right )\\
    =\ &\limsup_{m \to \infty}\PP\left (\exists t \in \RR^{\geq m} : |S'(t)| \geq \sqrt{\frac{2(\widehat V'(t)\rho_m^2 + 1)}{\rho_m^2}\log \left ( \frac{\sqrt{\widehat V'(t) \rho_m^2 +1}}{\alpha} \right )} \right ),
  \end{align}
  where the equality follows from the fact that $S'(t)$ and $\widehat V'(t)$ are piecewise-constant between integers.
  Note that since $V_t \to \infty$ as $t \to \infty$, we have that for any real sequence $(v_m)_{m =1}^\infty$ such that $v_m \to \infty$,
  \begin{align}
    &\limsup_{m \to \infty}\PP\left (\exists t \in \RR^{\geq m} : |S'(t)| \geq \sqrt{\frac{2(\widehat V'(t)\rho_m^2 + 1)}{\rho_m^2}\log \left ( \frac{\sqrt{\widehat V'(t) \rho_m^2 +1}}{\alpha} \right )} \right )\\
    =\ &\limsup_{m \to \infty}\PP\left (\exists v \geq v_m : |\widehat S(v)| \geq \sqrt{\frac{2(v\rho_m^2 + 1)}{\rho_m^2}\log \left ( \frac{\sqrt{v\rho_m^2 +1}}{\alpha} \right )} \right ).
  \end{align}
  Re-writing $v$ as $sv_md_m$ for some $s \in \RR^+$, the above probability can be written as
  \begin{align}
    &\PP\left (\exists v \geq v_m : |\widehat S(v)| \geq \sqrt{\frac{2(v\rho_m^2 + 1)}{\rho_m^2}\log \left ( \frac{\sqrt{v\rho_m^2 +1}}{\alpha} \right )} \right )\\
=\ &\PP\left (\exists sv_m d_m \geq v_m : |\widehat S(sv_md_m)| \geq\sqrt{\frac{2(sv_md_m\rho_m^2 + 1)}{\rho_m^2}\log \left ( \frac{\sqrt{sv_md_m\rho_m^2 +1}}{\alpha} \right )} \right )\\
    =\ &\PP\left (\exists sv_m d_m \geq v_m : |\widehat S(sv_md_m)| \geq\sqrt{\frac{2(s\cancel{v_md_m}c / \cancel{v_md_m} + 1)}{c/v_md_m}\log \left ( \frac{\sqrt{s\cancel{v_md_m}c / \cancel{v_md_m} +1}}{\alpha} \right )} \right ).
  \end{align}
  Applying \cref{lemma:interpolated-strong-approx}, we can approximate $\widehat S(sv_md_m)$ by $W(sv_md_m)$ up to a $o((sv_md_m)^{\kappa/2})$, yielding
  \begin{align}
    &\PP\left (\exists sv_m d_m \geq v_m : |\widehat S(sv_md_m)| \geq\sqrt{\frac{2(sc + 1)}{c/v_md_m}\log \left ( \frac{\sqrt{sc +1}}{\alpha} \right )} \right )\\
    =\ &\PP\left (\exists s \geq \frac{1}{d_m} : \left |W(sv_md_m) + o((sv_md_m)^{\kappa/2}) \right | \geq \sqrt{\frac{2(sc + 1)}{c/v_md_m}\log \left ( \frac{\sqrt{sc +1}}{\alpha} \right )} \right )\\
    =\ &\PP\left (\exists s \geq \frac{1}{d_m} : \left | \sqrt{v_md_m}W(s) + o((sv_md_m)^{\kappa/2}) \right | \geq \sqrt{\frac{2(sc + 1)}{c/v_md_m}\log \left ( \frac{\sqrt{sc +1}}{\alpha} \right )} \right ) \label{eq:wiener-process-scaling}\\
    =\ &\PP\left (\exists s \geq \frac{1}{d_m} : \left | W(s) + o \left ( \frac{(s v_m d_m)^{\kappa/2}}{(v_md_m)^{1/2}} \right ) \right | \geq \sqrt{\frac{2(sc + 1)}{c}\log \left ( \frac{\sqrt{sc +1}}{\alpha} \right )} \right ),
  \end{align}
  where~\eqref{eq:wiener-process-scaling} follows from the fact that $(W(cs))_{s \in \RR^+} \overset{d}{=} (\sqrt{c} W(s))_{s \in \RR^+}$ for any $c > 0$. Recalling the limit supremum of \cref{theorem:type-I-error} and taking the limsup of the above, we have
  \begin{align}
    & \limsup_{m \to \infty}\PP\left (\exists t \in \NN^{\geq m} : |S_t| \geq \sqrt{\frac{2(t\widehat \sigma_t^2\rho_m^2 + 1)}{\rho_m^2}\log \left ( \frac{\sqrt{t\widehat \sigma_t^2 \rho_m^2 +1}}{\alpha} \right )} \right )\\
    =\ &\limsup_{m\to\infty}\PP\left (\exists s \geq \frac{1}{d_m} : \left | W(s) + o \left ( s^{\kappa/2}(v_md_m)^{\frac{\kappa-1}{2}} \right ) \right | \geq \sqrt{\frac{2(sc + 1)}{c}\log \left ( \frac{\sqrt{sc +1}}{\alpha} \right )} \right ) \label{eq:type-I-final-limsup}
  \end{align}
  It remains to show that \eqref{eq:type-I-final-limsup} converges to $\alpha$. We will do so by separately showing that \eqref{eq:type-I-final-limsup} is upper- and lower-bounded by $\alpha$. Indeed for the upper bound, we have that
  \begin{align}
    &\sup_{s \geq 0} \left \{ \left | W(s) + o \left ( s^{\kappa/2}(v_md_m)^{\frac{\kappa-1}{2}} \right ) \right |  - \sqrt{\frac{2(sc + 1)}{c}\log \left ( \frac{\sqrt{sc +1}}{\alpha} \right )} \right \} \\
      \xrightarrow{d}\ &\sup_{s \geq 0} \left \{ |W(s)| - \sqrt{\frac{2(sc + 1)}{c}\log \left ( \frac{\sqrt{sc +1}}{\alpha} \right )} \right \}~~\text{as } m \to \infty,
  \end{align}
  and hence by classical results for boundary-crossing of Wiener processes, we have that
  \begin{align}
    \eqref{eq:type-I-final-limsup} &\equiv \limsup_{m\to\infty}\PP\left (\exists s \geq \frac{1}{d_m} : \left | W(s) + o \left ( s^{\kappa/2}(v_md_m)^{\frac{\kappa-1}{2}} \right ) \right | \geq \sqrt{\frac{2(sc + 1)}{c}\log \left ( \frac{\sqrt{sc +1}}{\alpha} \right )} \right )\\
    &\leq\limsup_{m\to\infty}\PP\left (\exists s \geq 0 : \left | W(s) + o \left ( s^{\kappa/2}(v_md_m)^{\frac{\kappa-1}{2}} \right ) \right | \geq \sqrt{\frac{2(sc + 1)}{c}\log \left ( \frac{\sqrt{sc +1}}{\alpha} \right )} \right )\\
    &=\PP\left (\exists s \geq 0 : |W(s)| \geq \sqrt{\frac{2(sc + 1)}{c}\log \left ( \frac{\sqrt{sc +1}}{\alpha} \right )} \right )\\
    &=\alpha.
  \end{align}
  Now, for the lower bound on \eqref{eq:type-I-final-limsup}, we have that
  \begin{align}
    \eqref{eq:type-I-final-limsup} &\equiv\limsup_{m\to\infty}\PP\left (\exists s \geq \frac{1}{d_m} : \left | W(s) + o \left ( s^{\kappa/2}(v_md_m)^{\frac{\kappa-1}{2}} \right ) \right | \geq \sqrt{\frac{2(sc + 1)}{c}\log \left ( \frac{\sqrt{sc +1}}{\alpha} \right )} \right ) \\
    &\geq \underbrace{\limsup_{m \to \infty}\PP\left (\exists s \geq 0 : \left |  W(s) + o \left ( s^{\kappa/2}(v_md_m)^{\frac{\kappa-1}{2}} \right ) \right | \geq \sqrt{\frac{2(sc + 1)}{c}\log \left ( \frac{\sqrt{sc +1}}{\alpha} \right )} \right )}_{(i)} \label{eq:type-I-error-later-Ws} \\
    &\quad- \underbrace{\limsup_{m \to \infty}\PP \left (\sup_{s \in [0, 1/d_m)} \left | W(s) + o \left ( s^{\kappa/2}(v_md_m)^{\frac{\kappa-1}{2}} \right ) \right | - \sqrt{\frac{2(sc + 1)}{c}\log \left ( \frac{\sqrt{sc +1}}{\alpha} \right )} \geq 0 \right )}_{(ii)}\label{eq:type-I-error-initial-Ws}
  \end{align}
  where the inequality~\eqref{eq:type-I-error-later-Ws} results from breaking the boundary-crossing event into the cases $s \in [0, 1/d_m)$ and $s \in [1/d_m, \infty)$. From the same argument as the upper bound, we have that $(i) = \alpha$, and hence it remains to show that $(ii) \to 0$. Indeed, the local modulus of continuity of the Wiener process states that
  \begin{equation}
    \limsup_{s \to 0} \frac{|W(s)|}{\sqrt{2s \log \log (1/s)}} = 1~~\text{almost surely,}
  \end{equation}
  and thus $(ii) \to 0$, so $\eqref{eq:type-I-final-limsup} \geq \alpha$. This completes the proof of \cref{theorem:type-I-error}.

\end{proof}

\subsection{Proof of \cref{proposition:delta-method}}\label{proof:delta-method}

\begin{proof}
  By Taylor's theorem, for each $t$, let $\widetilde \theta_t$ lie between $\widehat \theta_t$ and $\theta$ so that with probability one,
  \begin{align}
    g(\widehat \theta_t) - g(\theta) &= g'(\widetilde \theta_t)(\widehat \theta_t - \theta) \\
                                     &= (\widehat \theta_t - \theta) g'(\theta)  + \underbrace{(\widehat \theta_t - \theta)}_{(\star)} \underbrace{( g'(\widetilde \theta_t) - g'(\theta) )}_{(\dagger)} \\
    &= \frac{1}{t}\sum_{i=1}^t g'(\theta)\phi(Z_i) + o \left ( \sqrt{\log t / t} \right )
  \end{align}
  where we have used the fact that $(\star) = O(\sqrt{\log \log t / t}) + o(\sqrt{\log t / t})$ by the law of the iterated logarithm and $(\dagger) = o(1)$ by the almost-sure continuous mapping theorem, completing the proof.
\end{proof}

\subsection{Proof of \cref{proposition:asympcs-delayed-start}}\label{proof:asympcs-delayed-start}

First, we need the following lemma that calculates the probability of a mixture exponential Brownian motion exceeding $\eps > 0$ for any $t \geq 1$.

\begin{lemma}[A maximal (in)equality for mixture exponential Brownian motion with a delayed start]\label{lemma:exponential-brownian-motion-start-at-one}
  Define the continuous-time process $(M^1(t))_{t\geq 0}$ given by
  \begin{equation}
    M^1(t) := \frac{\exp \left \{ \frac{1}{2t} W(t)^2\right \}}{\sqrt{t}}
  \end{equation}
  where $(W(t))_{t\geq 0}$ is a Wiener process. Then,
  \begin{equation}
    \PP \left ( \exists t \geq 1 : M^1(t) > \eps \right ) = 2 \left [ 1- \Phi(a) + a\phi(a) \right ],
  \end{equation}
  where $a := \sqrt{2 \log \eps}$.
\end{lemma}

\begin{proof}
  The proof proceeds in 3 steps. First, we construct a conditional mixture martingale based on geometric Brownian motion akin to the one found in the proof of \cref{theorem:lindeberg-martingale-asympcs} but with the Wiener process replaced by itself plus independent standard Gaussian noise. Second, we show via iterated expectation and Step 1 that
\begin{equation}
  \EE \left [ \PP ( \exists t \geq 1 : M^1(t) \geq \eps \mid W(1)) \right ] = \EE \left [ \left (\eps^{-1} \exp \left \{ W(1)^2 / 2 \right \} \right ) \land 1 \right ].
\end{equation}
Finally, we evaluate the expectation on the right-hand side by noting that $W(1) \sim N(0, 1)$.
  \paragraph*{Step 1: Showing that $(M^{1-}(s))_{s \geq 0}$ forms a conditional nonnegative martingale given $Z$}
  Let $Z$ be a standard Gaussian random variable independent of $(W(s))_{s \geq 0}$. Then since geometric Brownian motion is a martingale, it is obvious that the process $(M^{1-}(s; \lambda))_{s \geq 0}$ given by
  \begin{align}
    M^{1-}(s ; \lambda) := \exp \left \{ \lambda [W(s) + Z] - s\lambda^2/2\right \}
  \end{align}
  forms a conditional nonnegative martingale given $Z$ with $M^{1-}(0; \lambda) \equiv \exp \left \{ \lambda Z \right \}$. Formally, $M^{1-}(s;\lambda)$ forms a nonnegative martingale with respect to the filtration $(\Fcal_s)_{s\geq0}$ given by $\Fcal_s := \sigma(W(s), Z)$. By Fubini's theorem, we have that for any probability distribution $F(\lambda)$, the mixture
  \begin{equation}
    \int_{\lambda \in \RR} M^{1-}(s; \lambda)dF(\lambda)
  \end{equation}
  also forms a conditional nonnegative martingale given $Z$ but now starting at $\int_{\lambda} M^{1-}(0; \lambda)dF(\lambda)$. In particular, using the techniques found in the proof of \cref{theorem:lindeberg-martingale-asympcs}, we have that for $dF(\lambda) := \frac{1}{\sqrt{2\pi}} \exp \left \{ \lambda^2 / 2 \right \}d\lambda$,
  \begin{align}
    M^{1-}(s) := \int_{\lambda \in \RR} M^{1-}(s;\lambda) dF(\lambda) = (s+1)^{-1/2} \exp \left \{ \frac{(W(s) + Z)^2}{2(s+1)} \right \}
  \end{align}
  forms a conditional nonnegative martingale starting at $M^{1-}(0) \equiv \exp \left \{ Z^2/2 \right \}$ conditional on $Z$.
  \paragraph*{Step 2: Showing that $\EE [\PP(\exists t \geq 1 : M^1(t) \geq \eps \mid W(1))] = \EE [(\eps^{-1} \exp \left \{ W(1)^2/2\right \}) \land 1 ]$}
  Notice that if we let $s := t - 1$, then $M^1(t)$ can be written as
  \begin{align}
    M^1(t) &:= \frac{\exp \left \{ \frac{1}{2t} W(t)^2\right \}}{\sqrt{t}}\\
    &= \left ( s + 1 \right )^{-1/2}\exp \left \{ \frac{W(s + 1)^2}{2(s + 1)} \right \},
  \end{align}
  and notice that by definition of the Wiener process, we have that $(W(s + 1))_{s\geq0} = (W'(s) + W(1))_{s\geq 0}$ where $(W'(s))_{s \geq 0}$ is a Wiener process independent of $W(1)$. Since $W(1) \sim N(0, 1)$, we have that
  \begin{equation}
    (M^1(t))_{t \geq 1} \mid W(1) \overset{d}{=} (M^{1-}(t-1))_{t \geq 1}\mid Z 
  \end{equation}
  and hence
  \begin{align}
    \PP \left ( \exists t \geq 1 : M^1(t) \geq \eps \right ) &= \EE \left [ \PP \left ( \exists t \geq 1 : M^1(t) \geq \eps \mid W(1) \right ) \right ]\\
                                                              &= \EE \left [ \PP \left ( \exists t \geq 1 : M^{1-}(t-1) \geq \eps \mid Z \right ) \right ] \\
                                                                                   &= \EE \left [ \PP \left ( \exists s \geq 0 : M^{1-}(s) \geq \eps \mid Z \right )  \right ]\\
    &= \EE \left [ \left (\eps^{-1} \exp \left \{ Z^2/2 \right \} \right ) \land 1 \right ]
  \end{align}
  which completes the proof of Step 2.

  \paragraph*{Step 3: Integrating out the distribution of $W(1)$ to obtain $\PP(\exists t \geq 1 : M^1(t) \geq \eps)$}
  Writing out the desired probability $\PP(\exists t \geq 1: M^1(t) \geq \eps)$ as the final expectation from Step 2, we have
  \begin{align}
    \PP(\exists t \geq 1 : M^1(t) > \eps) &= \EE \left [ \left (\eps^{-1} \exp \left \{ W(1)^2/2 \right \} \right ) \land 1 \right ] \\
                                          &= \int_{z \in \RR} \left ( \eps^{-1} \exp \left \{ w^2 / 2 \right \} \right ) \land 1 \cdot \dd \PP(W(1) \leq w) \\
                                          &= \int_{|w| > \sqrt{2\log \eps}} \frac{1}{\sqrt{2\pi}} \exp \left \{ -w^2/2 \right \} \dd w\\
                                          &\quad + \int_{|w| \leq \sqrt{2\log \eps}} \left ( \eps^{-1} \exp \left \{ w^2 / 2 \right \} \right ) \frac{1}{\sqrt{2\pi}} \exp \left \{ -w^2/2 \right \}\dd w \\
    &= 2 \left ( 1-\Phi(\sqrt{2 \log \eps}) \right ) + \frac{2}{\eps \sqrt{2 \pi}} \sqrt{2 \log \eps}  \\
    &= 2 \left ( 1-\Phi(a) + a\phi(a) \right ),
  \end{align}
  where $a := \sqrt{2 \log \eps}$. This completes the proof of \cref{lemma:exponential-brownian-motion-start-at-one}.
\end{proof}

The proof of \cref{lemma:exponential-brownian-motion-start-at-one} is similar in spirit to the derivation of \citet[Eq. (20)]{robbins1970statistical} but for continuous-time Wiener processes. Moreover, we only relied on a simple application of Ville's inequality and iterated expectation. Given the above result, we are ready to prove \cref{proposition:asympcs-delayed-start}.

\begin{proof}[Proof of \cref{proposition:asympcs-delayed-start}]
  The proof is a straightforward application of the ideas in the proof of \cref{theorem:type-I-error} (some of the relevant notation can be found therein) combined with \cref{lemma:exponential-brownian-motion-start-at-one}. Writing out the crossing probability of $|S_t|$ for any $t \geq m$ for a fixed $m$, we have that
  \begin{align}
    &\PP \left ( \exists t \geq m : |S_t| \geq \sqrt{ t \widehat \sigma_t^2\left [  a^2 + \log \left (\frac{t \widehat \sigma_t^2}{m \widehat \sigma_m^2} \right ) \right ]} \right ) \\
    =\ &\PP \left ( \exists t \geq m : |S'(t)| \geq \sqrt{\widehat V_t' \left [ a^2 + \log \left (\widehat V_t'/\widehat V_m' \right )  \right ]} \right ) \\
    =\ &\PP \left  ( \exists v \geq v_m : |\widehat S(v)| \geq \sqrt{v \left [ a^2 + \log(v / v_m)  \right ]}  \right ) \\
    =\ &\PP \left  ( \exists sv_m \geq v_m : |\widehat S(sv_m)| \geq \sqrt{sv_m \left [ a^2 + \log(sv_m / v_m)  \right ]}  \right ) \\
    =\ &\PP \left  ( \exists s \geq 1 : \left |W(sv_m) + o((sv_m)^{\kappa/2}) \right | \geq \sqrt{sv_m \left [ a^2 + \log s  \right ]}  \right ) \\
    =\ &\PP \left  ( \exists s \geq 1 : \left |\sqrt{v_m}W(s) + o((sv_m)^{\kappa/2}) \right | \geq \sqrt{sv_m \left [ a^2 + \log s  \right ]}  \right ) \\
    =\ &\PP \left  ( \exists s \geq 1 : \left |W(s) + o( s^{\kappa/2} v_m^{\frac{\kappa -1}{2}} ) \right | \geq \sqrt{s \left [ a^2 + \log s  \right ]}  \right ).
  \end{align}
  Now, since $\sup_{s \geq 1} \left \{ |W(s) + o(s^{\kappa/2} v_m^{(\kappa-1)/2})| - \sqrt{s [a^2 + \log s]} \right \} \xrightarrow{d} \sup_{s \geq 1} \left \{ |W(s)| - \sqrt{s [a^2 + \log s]} \right \}$ as $m \to \infty$, we have that
  \begin{align}
    &\limsup_{m \to \infty} \PP \left ( \exists t \geq m : |S_t| \geq \sqrt{ t \widehat \sigma_t^2\left [  a^2 + \log \left (\frac{t \widehat \sigma_t^2}{m \widehat \sigma_m^2} \right ) \right ]} \right ) \\
    =\ &\PP (\exists s \geq 1 : |W(s)| \geq \sqrt{s [a^2 + \log s]}).
  \end{align}
  Writing out the event $ \left \{ M^1(s) > \exp \left \{ a^2 / 2 \right \} \right \}$ with $M^1(s)$ defined as in \cref{lemma:exponential-brownian-motion-start-at-one}, we note that it is equivalent to the event in the above probability statement:
  \begin{align}
    &\left \{ \exists s \geq 1 : M^1(s) > \exp \left \{ a^2 / 2 \right \} \right \}\\
    =\ &\left \{ \exists s \geq 1: \frac{\exp \left \{ \frac{1}{2s} W(s)^2 \right \}}{\sqrt{s}} > \exp \left \{ a^2 / 2 \right \} \right \}\\
    =\ &\left \{ \exists s \geq 1: \frac{1}{2s} W(s)^2 - \log\sqrt{s} > a^2 / 2 \right \} \\
    =\ &\left \{ \exists s \geq 1: |W(s)| > \sqrt {s [a^2 + \log s]}  \right \},
  \end{align}
  and hence by \cref{lemma:exponential-brownian-motion-start-at-one}, we have that
  \begin{equation}
    \limsup_{m \to \infty} \PP \left ( \exists t \geq m : |S_t| \geq \sqrt{ t \widehat \sigma_t^2\left [  a^2 + \log \left (\frac{t \widehat \sigma_t^2}{m \widehat \sigma_m^2} \right ) \right ]} \right ) = 2(1-\Phi(a) + a\phi(a)),
  \end{equation}
  as desired. This completes the proof.
\end{proof}




\section{Additional discussions}

\subsection{One-sided asymptotic confidence sequences}\label{section:one-sided}
In Sections~\ref{section:universal-asympcs-mean} and~\ref{section:lindeberg-type-asympcs}, we derived universal two-sided \asympcs{}s for the means of independent random variables in the \iid{} and time-varying settings, respectively. Here, we give analogous one-sided bounds for the aforementioned settings. First, let us derive a one-sided \asympcs{} for the mean of \iid{} random variables analogous to \cref{theorem:acs}.

\begin{proposition}\label{proposition:acs-one-sided}
  Given the same setup as in~\cref{theorem:acs}, we have that
  \begin{equation}
    \widehat \mu_t - \widehat \sigma_t\sqrt{\frac{2(t\rho^2 + 1)}{t^2 \rho^2} \log \left ( 1 + \frac{\sqrt{t\rho^2 + 1}}{2\alpha} \right )}
  \end{equation}
  forms a lower $(1-\alpha)$-\asympcs{} for $\mu$ with the same rates as given in~\cref{theorem:acs}.
\end{proposition}

Notice that the $(1-\alpha)$-\asympcs{} of~\cref{proposition:acs-one-sided} resembles the $(1-2\alpha)$-\asympcs{} of~\cref{theorem:acs} but with an additional additive 1 inside the log. A similar phenomenon appears in the one- and two-sided sub-Gaussian \cs{}s of~\citet{howard2018uniform}. Recall, however, that their bounds are nonasymptotic and require much stronger assumptions (and in particular are not applicable to the observational causal inference setup of this paper).

Similar to the relationship between \iid{} (\cref{theorem:acs}) and martingale (\cref{theorem:lindeberg-martingale-asympcs}) two-sided \asympcs{}s, an analogue of~\cref{proposition:acs-one-sided} can be derived under martingale dependence with time-varying means and variances.

\begin{proposition}\label{proposition:martingale-asympcs-one-sided}
  Given the same setup and assumptions as~\cref{theorem:lindeberg-martingale-asympcs}, we have that
  \begin{equation}
    \widehat \mu_t - \sqrt{\frac{2(t \widehat \sigma_t^2 \rho^2 + 1)}{t^2 \rho^2} \log \left ( 1 +  \frac{\sqrt{t \widehat \sigma_t^2 \rho^2 + 1}}{2\alpha} \right )} 
  \end{equation}
  forms a lower $(1-\alpha)$-\asympcs{} for the time-varying average $\widetilde \mu_t := \frac{1}{t}\sum_{i=1}^t \mu_i$.
\end{proposition}

We will first prove a lemma concerning one-sided boundaries for sums of independent \emph{Gaussian} random variables, which in turn will be used to prove Propositions~\ref{proposition:acs-one-sided} and~\ref{proposition:martingale-asympcs-one-sided} shortly.

\begin{lemma}\label{lemma:one-sided-gaussian-cs}
  Suppose $\infseqt{G_t} \sim N(0, 1)$ is an \iid{} sequence of standard Gaussian random variables. Then, 
  \begin{equation}
    \PP \left ( \forall t \in \NN,\ \widetilde \mu_t \geq \underbrace{\frac{1}{t}\sum_{i=1}^t(\sigma_i G_i + \mu_i) -\sqrt{\frac{2(t\rho^2\widetilde \sigma_t^2  + 1)}{(t\rho)^2} \log \left (1 + \frac{\sqrt{t\rho^2\widetilde \sigma_t^2  + 1}}{2\alpha} \right ) }}_{L_t^\star} \right )  \geq 1-\alpha.
  \end{equation}
  In other words, $L_t^\star$ forms a nonasymptotic lower $(1-\alpha)$-\cs{} for $\widetilde \mu_t$.
\end{lemma}

\begin{proof}
The proof begins similarly to that of \cref{theorem:acs} but with a modified mixing distribution, and proceeds in four steps. First, we derive a sub-Gaussian nonnegative supermartingale (NSM) indexed by a parameter $\lambda \in \RR$ identical to that of \cref{theorem:acs}. Second, we mix this NSM over $\lambda$ using a \emph{folded} Gaussian density (rather than the classical Gaussian density used in the proof of \cref{theorem:acs}), and justify why the resulting process is also an NSM\@. Third, we derive an implicit lower \cs{} for $\infseq{\widetilde \mu^\star}{t}{1}$. Fourth and finally, we compute a closed-form lower bound for the implicit \cs{}.

\paragraph*{Step 1: Constructing the $\lambda$-indexed NSM} Similar to the proof of \cref{theorem:acs}, let $\infseqt {G_t}$ be an infinite sequence of \iid{} standard Gaussian random variables, and let $S_t := \sum_{i=1}^t \sigma_i G_i$. Then, we have that for any $\lambda \in \RR$,
\begin{equation}
    M_t(\lambda) := \exp \left \{ \lambda S_t - t \widetilde \sigma^2_t \lambda^2 / 2 \right \},
\end{equation}
forms an NSM with respect to the filtration given by $\Fcal_t := \sigma(G_1^t)$.

\paragraph*{Step 2: Mixing over $\lambda \in (0, \infty)$ to obtain a mixture NSM} Let us now construct a one-sided sub-Gaussian mixture NSM\@. First, note that the mixture of an NSM with respect to a probability density is itself an NSM \citep{robbins1970statistical,howard_exponential_2018} and is a simple consequence of Fubini's theorem. For our purposes, we will consider the density of a \emph{folded} Gaussian distribution with location zero and scale $\rho^2$. In particular, if $\Lambda \sim N(0, \rho^2)$, let $\Lambda_+ := |\Lambda|$ be the folded Gaussian. Then $\Lambda_+$ has a probability density function $f_{\rho^2}^+(\lambda)$ given by

\begin{equation}
  \label{eq:folded-gaussian-pdf}
  f_{\rho^2}^+ (\lambda) := \1(\lambda > 0)\frac{2}{\sqrt{2\pi \rho^2}} \exp \left \{ \frac{-\lambda^2}{2\rho^2} \right \}.
\end{equation}
Note that $f_{\rho^2}^+$ is simply the density of a mean-zero Gaussian with variance $\rho^2$, but truncated from below by zero, and multiplied by two to ensure that $f_{\rho^2}^+(\lambda)$ integrates to one.

Then, since mixtures of NSMs are themselves NSMs, the process $\infseq Mt0$ given by
\begin{equation}
    M_t := \int_{\lambda} M_t(\lambda) f_{\rho^2}^+(\lambda) d\lambda
\end{equation}
is an NSM. We will now find a closed-form expression for $M_t$. Some of the algebraic steps are the same as those in the proof of \cref{theorem:acs}, but we repeat them here for completeness. Writing out the definition of $M_t$, we have

\begin{align}
  M_t &:= \int_{\lambda \in \RR} \exp \left \{ \lambda S_t - t\widetilde \sigma^2_t \lambda^2/2 \right \} f_{\rho^2}^+(\lambda) d\lambda\\
    &= \int_{\lambda} \1(\lambda > 0)\exp\left \{ \lambda S_t - t\widetilde \sigma_t^2 \lambda^2/2 \right \} \frac{2}{\sqrt{2\pi \rho^2}} \exp \left \{ \frac{-\lambda^2}{2\rho^2} \right \}d\lambda \\
    &= \frac{2}{\sqrt{2\pi \rho^2}}\int_{\lambda}\1(\lambda > 0) \exp\left \{ \lambda S_t - t\widetilde \sigma_t^2 \lambda^2/2 \right \} \exp \left \{ \frac{-\lambda^2}{2\rho^2} \right \}d\lambda \\
    &= \frac{2}{\sqrt{2\pi \rho^2}} \int_\lambda \1(\lambda > 0)\exp \left \{ \lambda S_t - \frac{\lambda^2 (t\rho^2\widetilde \sigma_t^2 + 1)}{2 \rho^2}\right \} d\lambda \\
    &= \frac{2}{\sqrt{2\pi \rho^2}} \int_\lambda \1(\lambda > 0)\exp \left \{ \frac{-\lambda^2 (t\rho^2\widetilde \sigma_t^2 + 1) + 2\lambda \rho^2 S_t }{2\rho^2} \right \} d\lambda \\
    &= \frac{2}{\sqrt{2\pi \rho^2}} \int_\lambda\1(\lambda > 0) \underbrace{\exp \left \{ \frac{-a(\lambda^2 - \frac{b}{a} 2\lambda) }{2\rho^2} \right \}}_{(\star)} d\lambda,
\end{align}
where we have set $a:= t\rho^2\widetilde \sigma^2_t + 1$ and $b := \rho^2 S_t$. Completing the square in $(\star)$, we have that
\begin{align}
  \exp \left \{ \frac{-a(\lambda^2 - \frac{b}{a} 2\lambda) }{2\rho^2} \right \} &= \exp \left \{ \frac{-\lambda^2 + 2\lambda \frac{b}{a} + \left ( \frac{b}{a} \right )^2 - \left ( \frac{b}{a} \right )^2 }{2 \rho^2 /a} \right \}\\
                                                                                &= \exp \left \{ \frac{-(\lambda - b/a)^2}{2\rho^2/a} + \frac{a \left ( b/a \right )^2}{2\rho^2} \right \} \\
                                                                                &= \exp \left \{ \frac{-(\lambda - b/a)^2}{2\rho^2/a} \right \} \exp \left \{  \frac{b^2}{2a\rho^2} \right \}.
\end{align}
Plugging this back into our derivation of $M_t$ and multiplying the entire quantity by $a^{-1/2} / a^{-1/2}$, we have
\begin{align}
  M_t &= \frac{2}{\sqrt{2\pi \rho^2}} \int_\lambda\1(\lambda > 0) \underbrace{\exp \left \{ \frac{-a(\lambda^2 + \frac{b}{a} 2\lambda) }{2\rho^2} \right \}}_{(\star)} d\lambda\\
  &= \frac{2}{\sqrt{2\pi \rho^2}} \int_\lambda\1(\lambda > 0) \exp \left \{ \frac{-(\lambda - b/a)^2}{2\rho^2/a} \right \} \exp \left \{  \frac{b^2}{2a\rho^2} \right \} d\lambda\\
  &= \frac{2}{\sqrt{a}}\exp \left \{  \frac{b^2}{2a\rho^2} \right \}  \underbrace{\int_\lambda\1(\lambda > 0) \frac{1}{\sqrt{2\pi \rho^2/a}}\exp \left \{ \frac{-(\lambda - b/a)^2}{2\rho^2/a} \right \} d\lambda}_{(\star \star)}.
\end{align}
Now, notice that $(\star \star) = \PP(N(b/a, \rho^2 / a) \geq 0)$, which can be rewritten as $\Phi(b/\rho\sqrt{a})$, where $\Phi$ is the CDF of a standard Gaussian. Putting this all together and plugging in $a = t\rho^2\widetilde \sigma^2_t + 1$ and $b = \rho^2 S_t$, we have the following expression for $M_t$,
\begin{align}
  M_t &= \frac{2}{\sqrt{a}} \exp \left \{ \frac{b^2}{2a\rho^2} \right \} \Phi\left( \frac{b}{\rho \sqrt{a}} \right ) \nonumber\\
  &= \frac{2}{\sqrt{t\rho^2\widetilde \sigma^2_t + 1}} \exp \left \{ \frac{\rho^4 S_t^2}{2(t\rho^2\widetilde \sigma^2_t + 1)\rho^2} \right \}  \Phi\left( \frac{\rho^2 S_t}{\rho \sqrt{t\rho^2\widetilde \sigma_t^2 + 1}} \right ) \nonumber\\
  &= \frac{2}{\sqrt{t\rho^2\widetilde \sigma^2_t + 1}} \exp \left \{ \frac{\rho^2 S_t^2}{2(t\rho^2\widetilde \sigma^2_t + 1)} \right \} \Phi\left( \frac{\rho^2 S_t}{\sqrt{t\rho^2\widetilde \sigma^2_t + 1}} \right ).\label{eq:one-sided-mixture-NSM}
\end{align}

\paragraph*{Step 3: Deriving a $(1-\alpha)$-lower \cs{} $(L_t')_{t=1}^\infty$ for $(\widetilde \mu_t)_{t=1}^\infty$} Now that we have computed the mixture NSM $\infseq Mt0$, we apply Ville's inequality to it and ``invert'' a family of processes --- one of which is $\infseqt{M_t}$ --- to obtain an \emph{implicit} lower \cs{} (we will further derive an \emph{explicit} lower \cs{} in Step 4).

First, let $\infseq m t1$ be an arbitrary real-valued process --- i.e.~not necessarily equal to $\infseq {\mu} t1$ --- and define their running average $\widetilde m_t := \frac{1}{t} \sum_{i=1}^t m_i$. Define the partial sum process in terms of $\infseq {\widetilde m}t1$,
\[ S_t(\widetilde m_t) := S_t + t\widetilde \mu_t - t\widetilde m_t \]
and the resulting nonnegative process,
\begin{equation}
  \label{eq:e-process-mean-so-far}
M_t(\widetilde m_t) := \frac{2}{\sqrt{t\rho^2\widetilde \sigma^2_t + 1}} \exp \left \{ \frac{\rho^2 S_t(\widetilde m_t)^2}{2(t\rho^2\widetilde \sigma^2_t + 1)} \right \} \Phi\left( \frac{\rho S_t(\widetilde m_t)}{\sqrt{t\rho^2\widetilde \sigma^2_t + 1}} \right ).
\end{equation}
Notice that if $\widetilde m_t = \widetilde \mu_t$, then $S_t(\widetilde \mu_t) = S_t = \sum_{i=1}^t \sigma_i G_i$ and $M_t(\widetilde \mu_t) = M_t$ from Step 2. Importantly, $(M_t(\widetilde \mu_t))_{t=0}^\infty$ is an NSM\@. Indeed, by Ville's inequality, we have
\begin{equation}
  \label{eq:ville-lower-cs-mean-so-far}
  \PP(\exists t : M_t(\widetilde \mu_t ) \geq 1/\alpha ) \leq \alpha.
\end{equation}
We will now ``invert'' this family of processes to obtain an implicit lower boundary given by
\begin{equation}
  \label{eq:lower-cs-mean-so-far-implicit}
  L_t' := \inf \{\widetilde \mu_t : M_t(\widetilde \mu_t) < 1/\alpha \},
\end{equation}
and justify that $(L_t')_{t=1}^\infty$ is indeed a lower $(1-\alpha)$-\cs{} for $\widetilde \mu_t$. Writing out the probability of miscoverage at any time $t$, we have
\begin{align}
  \PP(\exists t : \widetilde \mu_t < L_t') &\equiv \PP\left (\exists t : \widetilde \mu_t < \inf_{\widetilde m_t} \{ M_t(\widetilde m_t) < 1/\alpha \} \right) \\
  &= \PP\left (\exists t : M_t(\widetilde \mu_t) \geq 1/\alpha \right) \\
  &\leq \alpha,
\end{align}
where the last line follows from Ville's inequality applied to $(M_t(\widetilde \mu_t))_{t=0}^\infty$. In particular, $L_t'$ forms a $(1-\alpha)$-lower \cs{}, meaning
\[ \PP(\forall t,\ \widetilde \mu_t \geq L_t') \geq 1-\alpha. \]

\paragraph*{Step 4: Obtaining a closed-form lower bound $(\widetilde L_t)_{t=1}^\infty$ for $(L_t')_{t=1}^\infty$}
The lower \cs{} of Step 3 is simple to evaluate via line- or grid-searching, but a closed-form expression may be desirable in practice, and for this we can compute a sharp lower bound on $L_t'$.

First, take notice of two key facts:
\begin{enumerate}[(a)]
  \item When $\widetilde m_t = S_t/t + \widetilde \mu_t $, we have that $S_t(\widetilde m_t) = 0$ and hence $M_t(\widetilde m_t) < 1$, and
  \item $S_t(\widetilde m_t)$ is a strictly decreasing function of $\widetilde m_t \leq S_t / t + \widetilde \mu_t$, and hence so is $M_t(\widetilde m_t)$.
\end{enumerate}
 Property (a) follows from the fact that $\Phi(0) = 1/2$, and that $\sqrt{t \rho^2\widetilde \sigma_t^2 + 1} > 1$ for any $\rho > 0$. Property (b) follows from property (a) combined with
 the definitions of $S_t(\cdot)$,
\[ S_t(\widetilde m_t) := S_t + t \widetilde \mu_t - t\widetilde m_t \]
and of $M_t(\cdot)$,
\[ M_t(\widetilde m_t) := \frac{2}{\sqrt{t\rho^2\widetilde \sigma_t^2 + 1}} \exp \left \{ \frac{\rho^2 S_t(\widetilde m_t)^2}{2(t\rho^2\widetilde \sigma_t^2 + 1)} \right \} \Phi\left( \frac{\rho S_t(\widetilde m_t)}{\sqrt{t\rho^2\widetilde \sigma_t^2 + 1}} \right ), \]
In particular, by facts (a) and (b), the infimum in \eqref{eq:lower-cs-mean-so-far-implicit} must be attained when $S_t(\cdot) \geq 0$. That is,
\begin{equation}
  \label{eq:sum_positive_at_arginf}
  S_t(L_t') \geq 0.
\end{equation}
Using \eqref{eq:sum_positive_at_arginf} combined with the inequality $1 - \Phi(x) \leq \exp \{ -x^2 / 2\}$ (a straightforward consequence of the Cram\'er-Chernoff technique), we have the following lower bound on $M_t(L_t')$:
\begin{align}
  M_t(L_t') &= \frac{2}{\sqrt{t\rho^2\widetilde \sigma_t^2 + 1}} \exp \left \{ \frac{\rho^2 S_t(L_t')^2}{2(t\rho^2\widetilde \sigma_t^2 + 1)} \right \} \Phi\left( \frac{\rho S_t(L_t')}{\sqrt{t\rho^2\widetilde \sigma_t^2 + 1}} \right ) \\
            &\geq \frac{2}{\sqrt{t\rho^2\widetilde \sigma_t^2  + 1}} \exp \left \{ \frac{\rho^2 S_t(L_t')^2}{2(t\rho^2\widetilde \sigma_t^2  + 1)} \right \} \left ( 1 - \exp\left \{ -\frac{\rho^2 S_t(L_t')^2}{2(t\rho^2\widetilde \sigma_t^2  + 1)} \right \} \right ) \\
            &= \frac{2}{\sqrt{t\rho^2\widetilde \sigma_t^2  + 1}} \left ( \exp \left \{ \frac{\rho^2 S_t(L_t')^2}{2(t\rho^2\widetilde \sigma_t^2  + 1)} \right \}  - 1 \right )\\
            &=: \widecheck M_t(L_t').
\end{align}
Finally, the above lower bound on $M_t(L_t')$ implies that $1/\alpha \geq M_t(L_t') \geq \widecheck M_t(L_t')$ which yields the following lower bound on $L_t'$:
\begin{align}
  \widecheck M_t(L_t') \leq 1/\alpha &\iff \frac{2}{\sqrt{t\rho^2\widetilde \sigma_t^2  + 1}} \left ( \exp \left \{ \frac{\rho^2 S_t(L_t')^2}{2(t\rho^2\widetilde \sigma_t^2  + 1)} \right \}  - 1 \right ) \leq 1/\alpha \\
  &\iff \exp \left \{ \frac{\rho^2 S_t(L_t')^2}{2(t\rho^2\widetilde \sigma_t^2  + 1)} \right \} \leq 1 + \frac{\sqrt{t\rho^2\widetilde \sigma_t^2  + 1}}{2\alpha}\\
  &\iff \frac{\rho^2 S_t(L_t')^2}{2(t\rho^2\widetilde \sigma_t^2  + 1)} \leq \log \left (1 + \frac{\sqrt{t\rho^2\widetilde \sigma_t^2  + 1}}{2\alpha} \right ) \\
  &\iff S_t(L_t') \leq \sqrt{\frac{2(t\rho^2\widetilde \sigma_t^2  + 1)}{\rho^2} \log \left (1 + \frac{\sqrt{t\rho^2\widetilde \sigma_t^2  + 1}}{2\alpha} \right ) }\\
  &\iff \sum_{i=1}^t\sigma_i G_i + t\widetilde \mu_t - t L_t' \leq \sqrt{\frac{2(t\rho^2\widetilde \sigma_t^2  + 1)}{\rho^2} \log \left (1 + \frac{\sqrt{t\rho^2\widetilde \sigma_t^2  + 1}}{2\alpha} \right ) }\\
  &\iff t L_t'  \geq \sum_{i=1}^t (\sigma_i G_i + \mu_i) -\sqrt{\frac{2(t\rho^2\widetilde \sigma_t^2  + 1)}{\rho^2} \log \left (1 + \frac{\sqrt{t\rho^2\widetilde \sigma_t^2  + 1}}{2\alpha} \right ) }\\
  &\iff L_t' \geq \underbrace{\frac{1}{t}\sum_{i=1}^t(\sigma_i G_i + \mu_i) -\sqrt{\frac{2(t\rho^2\widetilde \sigma_t^2  + 1)}{(t\rho)^2} \log \left (1 + \frac{\sqrt{t\rho^2\widetilde \sigma_t^2  + 1}}{2\alpha} \right ) }}_{L_t^\star},
\end{align}
and hence $\PP \left ( \forall t \in \NN,\ \widetilde \mu_t \geq L_t^\star \right ) \geq 1-\alpha$.
\end{proof}

\begin{proof}[Proof of~\cref{proposition:acs-one-sided}]
  In this case, the data $\infseqt{Y_t}$ are \iid{}, and hence we have that $\sigma_1=\sigma_2 = \cdots = \sigma$ and $\mu_1 = \mu_2 = \cdots = \mu$. First, notice that if we define $\beta = \rho \sigma$, then we can write $L_t^\star$ as
  \begin{equation}
    L_t^\star := \frac{\sigma}{t}\sum_{i=1}^t(G_i + \mu) - \sigma\sqrt{\frac{2(t\beta^2  + 1)}{(t\beta)^2} \log \left (1 + \frac{\sqrt{t\beta^2  + 1}}{2\alpha} \right ) }.
  \end{equation}
  Now, by the strong approximation of \citet{strassen1964invariance}, we have that
  \begin{equation}
    L_t^\star = \frac{1}{t}\sum_{i=1}^t Y_i - \sigma\sqrt{\frac{2(t\beta^2  + 1)}{(t\beta)^2} \log \left (1 + \frac{\sqrt{t\beta^2  + 1}}{2\alpha} \right ) } + \eps_t'
  \end{equation}
  where $\eps_t' = o(\sqrt{\log \log t /t})$. Writing the above in terms of an empirical standard deviation $\widehat \sigma_t$, we have by the proof of~\cref{lemma:strongGaussianApproxSampleAvg} that
 \begin{equation}
   L_t^\star = \frac{1}{t}\sum_{i=1}^t Y_i - \widehat \sigma_t\sqrt{\frac{2(t\beta^2  + 1)}{(t\beta)^2} \log \left (1 + \frac{\sqrt{t\beta^2  + 1}}{2\alpha} \right ) } + \eps_t
 \end{equation} 
 where $\eps_t = o \left ( \sqrt{\log \log t / t} \right )$ as well, and hence
 \begin{equation}
   \frac{1}{t}\sum_{i=1}^t Y_i - \widehat \sigma_t \sqrt{\frac{2(t\beta^2  + 1)}{(t\beta)^2} \log \left (1 + \frac{\sqrt{t\beta^2  + 1}}{2\alpha} \right )}
 \end{equation}
 forms a lower $(1-\alpha)$-\asympcs{} for $\mu$ with approximation rate $\eps_t$. This completes the proof.\footnote{While we wrote this final bound in terms of $\beta$, we leave the statement of the original result in terms of $\rho$ to maintain consistency with other boundaries throughout the paper. The change from $\beta$ to $\rho$ and back is entirely cosmetic, and does not affect the interpretation of the final result.}

\end{proof}

\begin{proof}[Proof of~\cref{proposition:martingale-asympcs-one-sided}]
Similar to the proof of~\cref{theorem:lindeberg-martingale-asympcs}, \cref{lemma:strong-approx-martingale} yields the following strong invariance principle
\begin{equation}
  \sum_{i=1}^t Y_i = \sum_{i=1}^t \sigma_i (G_i + \mu_i) + o \left ( V_t^{3/8} \log V_t \right ).
\end{equation}
Therefore, with probability at least $(1-\alpha)$,
\[ \forall t \geq 1, \  \widetilde \mu_t \geq \frac{1}{t}\sum_{i=1}^t Y_i - \sqrt{\frac{2(t\rho^2\widetilde \sigma_t^2  + 1)}{(t\rho)^2} \log \left (1 + \frac{\sqrt{t\rho^2\widetilde \sigma_t^2  + 1}}{2\alpha} \right ) } + o\left ( V_t^{3/8} \log V_t / t \right ). \]
In particular, we have that
\begin{equation}
  \label{eq:asympcs-widetilde-sigma-lower}
  \left ( \widehat \mu_t \pm \sqrt{\frac{2(t\widetilde \sigma_t^2 \rho^2 + 1)}{t^2 \rho^2} \log \left (1 +  \frac{\sqrt{t \widetilde \sigma_t^2 \rho^2 + 1}}{2\alpha} \right )}\right )
\end{equation}
forms a $(1-\alpha)$-AsympCS for $\widetilde \mu_t$. The derivation of an analogous lower \asympcs{} in terms of the empirical variance $\widehat \sigma_t^2$ proceeds similarly to Step 3 of the proof of \cref{theorem:lindeberg-martingale-asympcs}.
In particular, we get that
\[ \widehat \mu_t - \widetilde \boundary_t^\star := \widehat \mu_t - \sqrt{\frac{2(t\widehat \sigma_t^2 \rho^2 + 1)}{t^2\rho^2} \log \left (1 + \frac{\sqrt{t \widehat \sigma_t^2 \rho^2 + 1}}{2\alpha} \right )} + o\left ( \frac{\sqrt{V_t \log V_t}}{t} \right ) \]
forms a nonasymptotic $(1-\alpha)$-CS for $\widetilde \mu_t$, meaning $ \PP\left (\forall t \in \NN,\ \widetilde \mu_t \geq \widehat \mu_t - \widetilde \boundary_t^\star \right) \leq \alpha$. Combined with \cref{assumption:lindeberg_variance_estimator_consistent}, we have that
\[ \widehat \mu_t - \widetilde \boundary_t := \widehat \mu_t - \sqrt{\frac{2(t\widehat \sigma_t^2 \rho^2 + 1)}{t^2\rho^2} \log \left (1 + \frac{\sqrt{t \widehat \sigma_t^2 \rho^2 + 1}}{2\alpha} \right )}  \]
forms a $(1-\alpha)$-AsympCS for $\widetilde \mu_t$ since $\widetilde \boundary_t \asymp \sqrt{V_t \log V_t} / t$.
This completes the proof.

\end{proof}

\subsection{\texorpdfstring{Optimizing Robbins' normal mixture for $(t, \alpha)$}{Optimizing Robbins' normal mixture}}
\label{section:optimizingMixture}
In this section, we outline how one can choose $\rho$ to optimize the boundary $\widebar \boundary_t$ in Theorem~\ref{theorem:acs} for a specific time $t^\star$ and type-I error level $\alpha \in (0, 1)$.\footnote{We will discuss choosing $\rho$ for the two-sided \asympcs{} in~\cref{theorem:acs} but for the one-sided \asympcs{}s of~\cref{section:one-sided}, we suggest repeating the same argument but with $\alpha$ replaced by $2\alpha$.} We will outline both the (computationally inexpensive) exact solution, and the closed-form approximate solution. Note that the derivations that follow are essentially the same as those in \citet[Section 3.5]{howard2018uniform} but we repeat them here to keep our results self-contained.

\paragraph*{The exact solution}
Let $W_{-1}$ be the lower branch of the Lambert $W$ function \citep{corless1996lambertw}. Then,
\begin{equation}
\label{eq:exactRhoOpt}
     \argmin_{\rho > 0} \widebar \boundary_{t^\star}(\alpha) = \sqrt{\frac{-W_{-1}\left (-\alpha^2 \exp\left \{- 1\right \}\right ) - 1}{t^\star}}.
\end{equation}
\begin{proof}
Consider the boundary in Theorem~\ref{theorem:acs} at time $t$,
\[ \widebar \boundary_t(\alpha) := \sqrt{\frac{2(t\rho^2 + 1)}{t^2 \rho^2} \log \left ( \frac{\sqrt{t\rho^2 + 1}}{\alpha} \right )}. \]
Defining $x := \rho^2$ and after some simple algebra, notice that
\[ \argmin_{\rho > 0} \widebar \boundary_t(\alpha) = \sqrt{\argmin_{x > 0} f(x)}, \]
\[ \text{where }~~ f(x) := \frac{tx + 1}{t^2 x} \log \left ( \frac{tx + 1}{\alpha^2} \right ). \]
Notice that $\lim_{x \to 0} f(x) = \lim_{x \to \infty} f(x) = \infty$ and thus if we find that $df / dx = 0$ has exactly one positive solution, we know that it must be the minimizer of $f$.

To that end, it is straightforward to show that
\[ \frac{df}{dx} = -\frac{1}{t^2 x^2} \log \left ( \frac{tx + 1}{\alpha^2} \right )+ \frac{1}{tx}. \]
Setting the above to 0, we obtain
\[ \alpha^2 \exp \left \{ tx \right \} = tx + 1, \]
which, after some algebra, can be rewritten as
\begin{equation}
  -\alpha^2\exp \left \{ -1 \right \} = -(tx + 1) \exp \left \{ -(tx + 1) \right \}
\end{equation}
Notice that if we rewrite $y := -(tx + 1)$, we have that $y = W_{-1}\left (-\alpha^2 \exp\left \{ - 1 \right \} \right )$ where $W_{-1}$ is the lower branch of the Lambert $W$ function. Furthermore, $y = W_{-1}(z)$ only has a solution if $z \geq -e^{-1}$, requiring that $\alpha^2 \leq 1$, which we have trivially by the definition of $\alpha \in (0, 1)$. 
In summary, we have that
\begin{equation}
  \argmin_{\rho > 0} \widebar \boundary_{t^\star}(\alpha) = \sqrt{\frac{-W_{-1}\left (-\alpha^2 \exp\left \{ - 1\right \}\right ) - 1}{t^\star}}.
\end{equation}
This completes the proof.
\end{proof}

\paragraph*{An approximate solution}

We can derive a closed-form approximation to \eqref{eq:exactRhoOpt} by considering the Taylor series expansion to the Lambert $W$ function \citep{corless1996lambertw},
\begin{equation}
    W_{-1}(z) = \log(-z) - \log(-\log(-z)) + o(1).
\end{equation}
Replacing $W_{-1}(z)$ by $\log(-z) - \log(-\log(-z))$ in \eqref{eq:exactRhoOpt}, we obtain the following approximate solution,
\begin{equation}
  \label{eq:approxRhoOpt}
  \rho'(t^\star) := \sqrt{\frac{- 2\log \alpha + \log(-2 \log \alpha + 1)}{t^\star}}.
\end{equation}
In practice, we find that using \eqref{eq:approxRhoOpt} over \eqref{eq:exactRhoOpt} has negligible downstream effects on the resulting \cs{}s, but both are inexpensive to compute. Moreover, notice that $\rho'(t^\star)$ is quite similar to $\sqrt{2\log(1/\alpha) / t^\star}$, which is precisely what one would choose when sharpening a sub-Gaussian confidence interval based on the Cram\'er-Chernoff technique for a fixed sample size $t^\star$.

\subsection{Time-uniform convergence in probability is equivalent to almost sure convergence}
\label{section:TimeUniformInProb_AlmostSure}

In Theorems~\ref{theorem:acs}, \ref{theorem:csate_randomized}, and \ref{theorem:csate_observational}, we justified the asymptotic validity of our \asympcs{}s by showing that the approximation error 
\begin{equation}
\label{eq:almostSureConvergence}
    \varepsilon_t \xrightarrow{a.s.} 0
\end{equation} at a particular rate. At first glance, this may seem like a slightly stronger statement than required since we only need the approximation error $\varepsilon_t$ to vanish \emph{time-uniformly in probability}:
\begin{equation}
\label{eq:timeUniformlyInProb}
    \sup_{k \geq t} |\varepsilon_k| \xrightarrow{p} 0. 
\end{equation}
It turns out that \eqref{eq:almostSureConvergence} and \eqref{eq:timeUniformlyInProb} are equivalent. This is not new, but we present a proof for completeness.

\begin{proposition}
    Let $(X_n)_{n=1}^\infty$ be a sequence of random variables. Then, 
    \[ X_n \xrightarrow{a.s.} 0 \iff \sup_{k \geq n} |X_n| \xrightarrow{p} 0. \]
\end{proposition}

\begin{proof}
    First, we prove $(\implies)$. By the continuous mapping theorem, $|X_n| \xrightarrow{a.s.} 0$. Thus,
    \begin{align}
        1 = \PP\left (\lim_n |X_n| = 0 \right) &\leq \PP \left( \limsup_n |X_n| = 0 \right) \leq  1.
    \end{align}
    In other words, $\sup_{k \geq n} | X_k| \xrightarrow{a.s.} 0$, which implies $\sup_{k \geq n} |X_k| \xrightarrow{p} 0$.
    \\
    
    \noindent Now, consider $(\impliedby)$. Suppose for the sake of contradiction that $\PP\left ( \lim_n |X_n| = 0 \right ) < 1$. Then with some probability $\delta > 0$, we have that $\lim_n |X_n| \neq 0$, meaning there exists some $\epsilon > 0$ such that $|X_k| > \epsilon$ for some $k \geq n$ no matter how large $n$ is. In other words,
    \begin{align}
        \delta &< \PP\left ( \lim_{n \rightarrow \infty} \sup_{k \geq n} | X_k | > \epsilon \right ) \\
        &\leq \PP \left ( \sup_{k \geq n} |X_k| > \epsilon \right ) ~~\text{for any $n \geq 1$.}
    \end{align}
    In particular, $\PP\left ( \sup_{k \geq n} |X_k| > \epsilon \right ) \not\rightarrow 0$, which is equivalent to saying $\sup_{k \geq n}|X_k| \overset{p}{\not\rightarrow} 0$, a contradiction. This completes the proof.
\end{proof}

\subsection{Comparing \asympcs{}s to group-sequential repeated confidence intervals}\label{section:group-sequential}
Another approach to sequential inference is via so-called \emph{group-sequential trials} and \emph{repeated confidence intervals} (see the now-classical text of \citet{jennison1999group}). In brief, group-sequential trials allow the analyst to peek at \ci{}s at certain prespecified times $(t_1, t_2, ..., t_K)$. They differ from the ``anytime-valid'' approach of \cs{}s and (and hence \asympcs{}s) in that they do not permit continuous monitoring (i.e.~updating of inferences for each new data point collected) and require a (fixed, data-independent) maximum sample size $t_K$. In this way, they can be thought of as fixed-time CLT-based \ci{}s that permit a fixed number of peeks prior to $t_K$ for early stopping. \asympcs{}s by contrast can be continuously monitored indefinitely, allowing the study to continue for as long as needed by the analyst.

\begin{figure}[!htbp]
  \centering
  \includegraphics[width=\columnwidth]{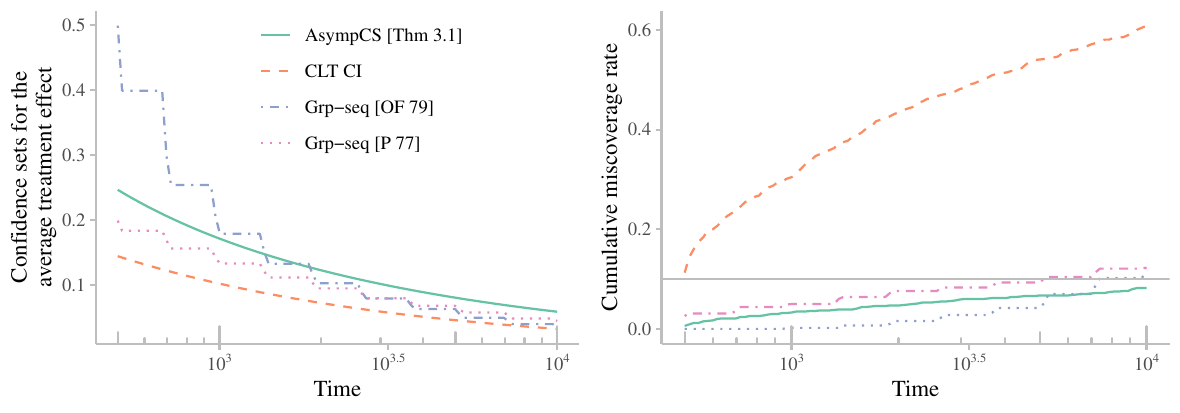}
  \caption{A comparison of 90\% confidence sets for the average treatment effect in a simulated sequential experiment using (a) the \asympcs{} of \cref{theorem:csate_randomized}, (b) the CLT-based \ci{}, and the group-sequential repeated \ci{}s of (c) \citet{o1979multiple}, and (d) \citet{pocock1977group}. The group-sequential methods used a start time of $t_1 = 500$ and a maximal sample size of $t_{10} = 10,000$ with $10$ logarithmically-spaced peeking times. The jagged drops in their widths correspond to \ci{}s being updated at each peek. Notice in the left-hand side plot the group-sequential repeated \ci{}s lie somewhere in between the \asympcs{} and the CLT \ci{} (especially those using the Pocock boundary). Furthermore, notice in the right-hand side plot that the cumulative miscoverage rate of the \asympcs{} lies below $\alpha = 0.1$ while the group-sequential \ci{}s slightly exceed it. Nevertheless, those of the CLT-based \ci{}'s miscoverage rates greatly exceed the other three, diverging quickly beyond $\alpha = 0.1$ and exceeding $0.6$ by $t_{10} = 10,000$.}
  \label{fig:group-sequential}
\end{figure}
Here, we provide a simulation comparing the widths and cumulative miscoverage rates of \asympcs{}s to two popular repeated \ci{}s from the group-sequential literature --- namely those of \citet{pocock1977group} and \citet{o1979multiple} and display the classical CLT-based \ci{} alongside them for reference (\cref{fig:group-sequential}).
Perhaps unsurprisingly, group-sequential methods (and especially those using the boundary of \citet{pocock1977group}) tend to lie somewhere in between CLT \ci{}s and \asympcs{}s. This added tightness over \asympcs{}s of course comes at the cost of less flexibility, especially in regards to continuous monitoring and unbounded time horizons discussed above.

\subsection{The Lyapunov-type condition implies the Lindeberg-type condition}\label{section:lyapunov-implies-lindeberg}

Here, we give a brief proof of the fact that the Lyapunov-type condition discussed in the paragraph following \cref{theorem:lindeberg-martingale-asympcs} indeed implies \cref{assumption:lindeberg-type_condition}.

\begin{proof}
  
  Suppose that the following Lyapunov-type condition holds with some $\delta > 0$:
  \begin{equation}
    \sum_{t=1}^\infty\frac{\EE \left ( |X_t|^{2 + \delta} \mid Y_1^{t-1} \right )}{\sqrt{V_t}^{2+\delta}} < \infty,
  \end{equation}
  and we want to show that \cref{assumption:lindeberg-type_condition} holds for some $\kappa \in (0, 1)$. Indeed, set $\kappa := \frac{1+\delta/2}{1+\delta}$. Note that $1 \leq V_t^{-\kappa \delta/2} |X_t|^{\delta}$ whenever $X_t^2 > V_t^{\kappa}$, and hence
  \begin{align}
    \sum_{t=1}^\infty \frac{\EE \left ( X_t^2 \1 \left \{ X_t^2 > V_t^\kappa   \right \} \mid Y_1^{t-1}\right )}{V_t^\kappa} &\leq \sum_{t=1}^\infty\frac{\EE \left ( V_t^{-\kappa \delta/2} |X_t|^{2 + \delta} \1 \left \{ X^2_t > V_t^\kappa \right \} \mid Y_1^{t-1} \right )}{V_t^\kappa} \\
                                                                                                                             & = \sum_{t=1}^\infty\frac{\EE \left ( |X_t|^{2+\delta} \1 \left \{ X^2_t > V_t^\kappa \right \} \mid Y_1^{t-1} \right )}{V_t^{\kappa(1+\delta/2)}}\\
                                                                                                                             & = \sum_{t=1}^\infty\frac{\EE \left ( |X_t|^{2+\delta} \mid Y_1^{t-1} \right )}{\sqrt{V_t}^{2+\delta}} < \infty,
  \end{align}
  which completes the proof.

\end{proof}

\subsection{On martingale \asympcs{}s for running average treatment effects}\label{section:martingale-time-varying-trt-effects}

In \cref{remark:predictable_nuisance_estimation}, we mentioned that it is possible to derive \asympcs{}s for the running average treatment effect $\widetilde \psi_t$ in randomized experiments without sample splitting. The following proposition is a corollary of \cref{theorem:lindeberg-martingale-asympcs} applied to a particular sequence of estimated influence functions. In what follows, we use $R_t$ for the outcome of subject $t$ instead of $Y_t$ so that it is not confused with the appearance of $Y_t$ in \namecref{assumption:lindeberg_variance_does_not_vanish}s~\ref{assumption:lindeberg_variance_does_not_vanish}, \ref{assumption:lindeberg-type_condition}, and \ref{assumption:lindeberg_variance_estimator_consistent} which we will refer to.

\begin{proposition}\label{proposition:martingale-time-varying-trt-effect}
  Let $\infseqt{Z_t} \equiv \infseqt{X_t, A_t, R_t}$ be independent triplets as in \cref{theorem:lyapunov_ate_randomized} and let $\widehat \mu_t^a$ be an estimator for $\mu^a$ derived from $Z_1^{t-1}$; $a \in \{0, 1\}$. If \namecref{assumption:lindeberg_variance_does_not_vanish}s~\ref{assumption:lindeberg_variance_does_not_vanish}, \ref{assumption:lindeberg-type_condition}, and \ref{assumption:lindeberg_variance_estimator_consistent} hold but with $Y_t$ replaced by $f_t(Z_t)$ everywhere given by
  \begin{equation}
     f_t(Z_t) := \left \{ \widehat \mu^1_t(X_t) - \widehat \mu^0_t(X_t) \right \} + \left ( \frac{A_t}{\pi(X_t)} - \frac{1-A_t}{1-\pi(X_t)} \right )\left \{ R_t - \widebar \mu^a(X_t) \right \}.
  \end{equation}
  Then for any $\rho > 0$,
  \begin{equation}
    \widetilde C_t := \frac{1}{t}\sum_{i=1}^t f_i(Z_i) \pm \sqrt{\frac{2(t\widehat \sigma_t^2 \rho^2 + 1)}{t^2\rho^2} \log \left ( \frac{\sqrt{t\widehat \sigma_t^2 \rho^2 + 1}}{\alpha} \right )}
  \end{equation}
  forms a $(1-\alpha)$-\asympcs{} for the running average treatment effect $\widetilde \psi_t$.
\end{proposition}
For certain practical purposes, however, we still recommend sequential sample-splitting or cross fitting as described in \cref{section:csate} since the conditions of \cref{proposition:martingale-time-varying-trt-effect} are somewhat less transparent.

\subsection{Explicit connections between ``delayed start'' boundaries and other works}\label{section:explicit-connections-bibaut}

In \cref{section:type-I-error}, we elaborated on certain connections between \citet{bibaut2022near} and the results in \cref{proposition:asympcs-delayed-start} and its corollary in \eqref{eq:delayed-start-corollary}. In \namecref{claim:bv1}s~\ref{claim:bv1} and~\ref{claim:bv2} we explicitly show that $\widetilde C_t^{\mathrm{DS}\star}$ is almost surely tighter than the bound found in version 1 and identical to that found in version 2 of \citet{bibaut2022near}. In what follows, we let $C_t^\bvone$ and $C_t^\bvtwo$ denote the bounds found in versions 1 and 2 of \citet{bibaut2022near}, respectively.

\begin{claim}[$\widetilde C_t^{\DS \star}$ is almost surely tighter than $C_t^\bvone(\lambda)$]\label{claim:bv1}
  Let $\infseqt{X_t}$ be \iid{} random variables with unknown mean $\mu$ and known variance $\sigma$. Let $\widehat \mu_t := \frac{1}{t}\sum_{i=1}^t X_i$ be the sample mean and $\alpha \in (0, 1)$ the desired type-I error level. For any fixed $\lambda > 0$, consider $C_t^\bvone(\lambda)$ from version 1 of \citet[Example 1]{bibaut2022near} given by
  \begin{equation}\label{eq:bv1}
    C_t^\bvone(\lambda) := \left [ \widehat \mu_t \pm \sigma t^{-1} \cdot (t+\lambda)^{1/2} (\log(t+\lambda) - \log m - 2\log \widetilde \alpha(\alpha))^{1/2} \right ],
  \end{equation}
  where $\widetilde \alpha(\alpha)$ is the unique solution to $2\widetilde \alpha(\alpha) \sqrt{-\log \widetilde \alpha(\alpha) / \pi} + 2\left [1-\Phi\left (\sqrt{-2 \log \widetilde \alpha(\alpha)}\right ) \right ] = \alpha$.\footnote{Note that in version 1 of \citet{bibaut2022near}, Algorithm 1 includes slightly different constants (specifically some multiples of 2 are moved around), but what we have written here matches their code at \href{https://github.com/nathankallus/UniversalSPRT}{github.com/nathankallus/UniversalSPRT} which we believe contains the correct bound. Note that in version 2 of \citet{bibaut2022near}, their constants agree with what is written here.} Then, $C_t^\bvone(\lambda)$ in \eqref{eq:bv1} can be equivalently re-written as
  \begin{equation}\label{eq:bv1-rewrite}
    C_t^\bvone(\lambda) \equiv \left [ \widehat \mu_t \pm \sigma \sqrt{\frac{t+\lambda}{t^2} \left ( a^2 + \log \left ( \frac{t + \lambda}{m} \right ) \right )} \right ],
  \end{equation}
  where $a \geq 0$ is the unique solution to $2(1-\Phi(a) + a\phi(a)) = \alpha$, and thus $\widetilde C_t^{\DS\star}$ is a strict subset of $C_t^\bvone(\lambda)$ for any $\lambda > 0$.
\end{claim}

\begin{proof}
  Begin by defining $b \geq 0$ as
  \begin{equation}
    b := \sqrt{-2\log \widetilde \alpha(\alpha)} \quad\text{and thus}\quad\widetilde \alpha(\alpha) = \exp \left \{-b^2 / 2 \right \},
  \end{equation}
  noting that the transformation $x \mapsto \sqrt{-2\log x}$ is a bijection for $x \in [0, 1]$. We will now show that $b$ solves $2(1-\Phi(b) + b\phi(b)) = \alpha$ if and only if $\widetilde \alpha$ solves the equation given in the definition of $C_t^\bvone$. Indeed, writing out the equation that $\widetilde \alpha$ solves, we have
  \begin{align}
    &2 \widetilde \alpha \sqrt{-\log \widetilde \alpha / \pi} + 2 \left [ 1-\Phi \left ( \sqrt{-2 \log \widetilde \alpha} \right ) \right ] = \alpha \\
    \iff &2 \exp \{-b^2 / 2\} \sqrt{b^2 / (2\pi)} + 2 \left [ 1 - \Phi \left ( b \right ) \right ] = \alpha \\
    \iff &2 b \underbrace{\frac{1}{\sqrt{2\pi}} \exp \left \{ -b^2 / 2 \right \}}_{\equiv \phi(b)} + 2[1-\Phi(b)] = \alpha\\
    \iff &2 \left [ 1-\Phi(b) + b\phi(b) \right ] = \alpha.
  \end{align}
  Now, it remains to show that $C_t^\bvone(\lambda)$ can be written in the form given in \eqref{eq:bv1-rewrite}. Indeed, writing out the boundary in \eqref{eq:bv1}, we have
  \begin{align}
    &\sigma t^{-1} \cdot (t + \lambda)^{1/2} \left ( \log (t + \lambda) - \log m - 2\log \widetilde \alpha \right )^{1/2}\\
    =\ &\sigma t^{-1}\cdot  \sqrt{(t+\lambda) \left ( \log (t + \lambda) - \log m + b^2 \right ) }\\
    =\ &\sigma\sqrt{\frac{t+\lambda}{t^2} \left ( b^2 + \log \left ( \frac{t+\lambda}{m} \right ) \right )},
  \end{align}
  and since $b$ solves $2[1-\Phi(b)+b\phi(b)] = \alpha$ as demonstrated above, this completes the proof of the claim.
\end{proof}

\begin{claim}[$\widetilde C_t^{\DS\star}$ is equivalent to $C_t^\bvtwo$]\label{claim:bv2}
  Consider $\widetilde C_t^{\DS\star}$ as before and consider $C_t^\bvtwo$ as in version 2 of \citet[p.~4]{bibaut2022near}:
  \begin{equation}
    C_t^\bvtwo := \left [ \widehat \mu_t \pm \frac{\sigma}{t} \cdot \sqrt{ t\cdot \left ( -2\log \widetilde \alpha + \log \left ( \frac{t}{m} \right ) \right )} \right ],
  \end{equation}
  where $\widetilde \alpha$ solves $h_1(-\log \widetilde \alpha) = \alpha$ and $h_1(x)$ is given by $h_1(x) = 2\exp \{-x\} \sqrt{x/\pi} + 2(1-\Phi(\sqrt{2x}))$; $x \geq 0$. Then $C_t^\bvtwo = \widetilde C_t^{\DS\star}$.
\end{claim}
\begin{proof}
  Similar to the proof of \cref{claim:bv1}, we define $b := \sqrt{-2 \log \widetilde \alpha}$ and show that $b$ solves $2[1-\Phi(b) + b\phi(b)] = \alpha$ if and only if $\widetilde \alpha$ solves $h_1(-\log \widetilde \alpha) = \alpha$. Indeed,
  \begin{align}
    &h_1(-\log \widetilde \alpha) = \alpha\\
    \iff &2\exp \{\log \widetilde \alpha\} \sqrt{-\log \widetilde \alpha / \pi} + 2(1-\Phi(\sqrt{-2 \log \widetilde \alpha})) = \alpha\\
    \iff &2 \exp \{ -b^2 / 2\} \sqrt{b / (2\pi)} + 2(1-\Phi(\sqrt{b^2}))\\
    \iff &2 b \underbrace{\frac{1}{\sqrt{2\pi}}\exp \{ -b^2 / 2\}}_{\equiv \phi(b)} + 2(1-\Phi(b))\\
    \iff &2 \left [ 1-\Phi(b) + b \phi(b) \right ].
  \end{align}
  Finally, noting that $b^2 = -2\log \widetilde \alpha$, we have that
  \begin{equation}
    C_t^\bvtwo \equiv \left [ \widehat \mu_t \pm \sigma \sqrt{\frac{1}{t} \left ( b^2 + \log(t/m) \right )} \right ],
  \end{equation}
  and since $b$ solves $2[1-\Phi(b) + b\phi(b)] = \alpha$, we have that $C_t^\bvtwo = \widetilde C_t^{\DS \star}$ for each $t$, which completes the proof of the claim.
\end{proof}





\subsection{A brief review of efficient estimators}\label{section:efficient-estimators-review}
For a detailed account of efficient estimation in semiparametric models, we refer readers to \citet{bickel1993efficient}, \citet{van2002semiparametric}, \citet{van2003unified}, \citet{tsiatis2007semiparametric} and \citet{kennedy2016semiparametric}, but we provide a brief overview of their fundamental relevance to estimation of the ATE here.

A central goal of semiparametric efficiency theory is to characterize the set of \emph{influence functions} of $\psi$ --- the summands found in sample averages forming consistent and asymptotically normal regular estimators of $\psi$. Of particular interest is finding the \emph{efficient influence function} (EIF) as it is the one with the smallest variance (which itself acts as a semiparametric analogue of the Cramer-Rao lower bound), hence providing a benchmark for constructing optimal estimators, at least in an asymptotic local minimax sense. 
In the case of $\psi$, the (uncentered) EIF is given by
\begin{equation}\label{eq:efficientInfluenceFunction-brief-review}
    f(z) \equiv f(x, a, y) := \left \{ \mu^1(x) - \mu^0(x) \right \} + \left ( \frac{a}{\pi(x)} - \frac{1-a}{1-\pi(x)} \right )\left \{ y - \mu^a(x) \right \},
\end{equation}
where $\mu^a(x) := \EE \left ( Y \mid X = x, A= a \right )$ is the regression function among those treated at level $a \in \{0, 1\}$ and $\pi(x) := \PP\left ( A = 1 \mid X = x \right )$ is the propensity score (i.e.~probability of treatment) for an individual with covariates $x$. In particular, this means that no estimator of $\psi$ based on $t$ observations can have asymptotic mean squared error smaller than $\Var(f(Z))/t$ without imposing additional assumptions.

Of course, the exact values of $\eta := (\mu^1, \mu^0, \pi)$ are not known in general --- $\pi$ is only known in randomized experiments, but not observational studies, while $(\mu^1, \mu^0)$ are typically not known in either. As such, we will need to replace these ``nuisance functions'' $\eta$ with data-dependent estimates $\widehat \eta$, but using the same data to construct estimators $\widehat \eta$ and $\widehat f$ for both $\eta$ and $f$ (sometimes referred to as ``double-dipping'') complicates the analysis of the downstream estimator of $\psi$. A clever and simple remedy used throughout the semiparametric literature is to split the sample, whereby a random subset of the data are used to estimate $\eta$, while the remaining data are used to construct $\widehat f$, greatly simplifying downstream analysis \citep{robins2008higher,zheng2010asymptotic,chernozhukov2017double}.

In a randomized experiment, the joint distribution of $(X, Y)$ is unknown but the conditional distribution of $A \mid X = x$ is known to be Bernoulli($\pi(x)$) by design. In this case, our statistical model for $Z$ is a proper semiparametric model, and hence there are infinitely many influence functions, all of which take the form,
\begin{equation}
     \widebar f(z) \equiv \widebar f(x, a, y) := \left \{ \widebar \mu^1(x) - \widebar \mu^0(x) \right \} + \left ( \frac{a}{\pi(x)} - \frac{1-a}{1-\pi(x)} \right )\left \{ y - \widebar \mu^a(x) \right \}, 
\end{equation}
where $\widebar \mu^a : \RR^d \mapsto \RR$ is any function. However, when the joint distribution of $(X,A,Y)$ is left completely unspecified (such as in an observational study with unknown propensity scores), our statistical model for $\PP$ is nonparametric, and hence there is only one influence function, the EIF given in \eqref{eq:efficientInfluenceFunction-brief-review}.

Not only does the EIF $f(z)$ provide us with a benchmark against which to compare estimators, but it hints at the first step in deriving the most efficient estimator. Namely, $\frac{1}{t}\sum_{i=1}^t f(Z_i)$ is a consistent estimator for $\psi$ with asymptotic variance equal to the efficiency bound, $\Var(f)$ by construction. However, $f(Z)$ depends on possibly unknown nuisance functions $\eta := (\mu^1, \mu^0, \pi)$. A natural next step would be to simply estimate $\eta$ from the data $\infseq Zt1$. Crucially, it is possible to ensure that only a negligible amount of additional estimation error is incurred by replacing $\eta$ by a data-dependent estimate $\widehat \eta_t$ --- the essential technique here being sample splitting and cross fitting \citep{robins2008higher, zheng2010asymptotic, chernozhukov2017double}.

\subsection{On the sharpness of \asympcs{}s using efficient influence functions}
\label{section:unimprovability}
Consider the \asympcs{} of Theorem~\ref{theorem:csate_observational},
\begin{equation}
  \label{eq:observational-cs-recall}
  \widehat \psi_t^\times \pm \underbrace{\sqrt{\widehat \Var_{t}(\widehat f)}}_{(iii)} \cdot \underbrace{\sqrt{\frac{2(t\rho^2 + 1)}{t^2 \rho^2} \log \left ( \frac{\sqrt{t\rho^2 + 1}}{\alpha} \right ) }}_{(i)} ~~~ \text{with rate}~ \underbrace{o\left ( \sqrt{\frac{\log \log t}{t}} \right )}_{(ii)}.
\end{equation}
It is natural to wonder whether \eqref{eq:observational-cs-recall} can be tightened. In a certain sense, \eqref{eq:observational-cs-recall} inherits optimality from its three main components: $(i)$ Robbins' normal mixture boundary, $(ii)$ the approximation error rate, and $(iii)$ the estimated standard deviation $\sqrt{\widehat \Var_{t}(\widehat f)}$ of the efficient influence function $f$.

\paragraph*{Term \textit{(i)}} Starting with the width, we have that in the case of \iid{} Gaussian data $G_1, G_2, \dots \simiid \Ncal(\mu, \sigma^2)$, Robbins' normal mixture confidence sequence \citep{robbins1970statistical} is obtained by first showing that
\[ M_t(\mu) := \exp \left \{ \frac{\rho^2 (\sum_{i=1}^t(G_i-\mu))^2 }{2(t \rho^2 + 1)} \right \} (t \rho^2 + 1)^{-1/2} \]
is a nonnegative martingale starting at one, and hence by Ville's inequality \citep{ville1939etude},
\[ \PP(\exists t \geq 1 : M_t(\mu) \geq 1/\alpha) \leq \alpha. \]
The resulting confidence sequence $\widebar C_t^\Ncal$ at each time $t$ is defined as the set of $m$ such that $M_t(m) < 1/\alpha$, i.e. $\widebar C_t^\Ncal := \{ m \in \RR : M_t (m) < 1/\alpha \}$ and consequently,
\[ \PP(\exists t \geq 1 : \mu \notin \widebar C_t^\Ncal) = \PP(\exists t \geq 1 : M_t(\mu) \geq 1/\alpha) \leq \alpha. \]
This inequality is extremely tight, since Ville's inequality almost holds with equality for nonnegative martingales. Technically, the paths of the martingale need to be continuous for equality to hold, which can only happen in continuous time (such as for a Wiener process). However, any deviation from equality only holds because of this ``overshoot'' and in practice, the error probability is almost exactly $\alpha$. This means that the normal mixture confidence sequence $\widebar C_t^\Ncal$ cannot be uniformly tightened: any improvement for some times will necessarily result in looser bounds for others.
For a precise characterization of this optimality for the (sub)-Gaussian case, see \citet[Section 3.6]{howard2018uniform}, or \citet{ramdas2020admissible} for a more general discussion of admissible confidence sequences.

\paragraph*{Term \textit{(ii)}} The error incurred from almost-surely approximating a sample average $\frac{1}{t}\sum_{i=1}^t f(Z_i)$ of influence functions by Gaussian random variables is a direct consequence of \citet{strassen1964invariance} and improvements under $q > 2$ finite absolute moments by \citet*{komlos1975approximation, komlos1976approximation} and \citet{major1976approximation}, which are unimprovable without additional assumptions. Further approximation errors result from using $\widehat \Var_t(\widehat f)$ to estimate $\Var(f)$, where almost-sure law of the iterated logarithm rates appear, and are themselves unimprovable.

\paragraph*{Term \textit{(iii)}} Using the approximations mentioned in $(ii)$ permits the use of Robbins' normal mixture confidence sequence in $(i)$. However, a factor of $\sqrt{\widehat \Var_t(\widehat f)}$ necessarily appears in front of the width as an estimate of the standard deviation $\sqrt{\Var(f)}$ of the efficient influence function $f$ discussed in \cref{section:csate}. Importantly, $\sqrt{\Var(f)}$ corresponds to the semiparametric efficiency bound, so that no estimator of $\psi$ can have asymptotic mean squared error smaller than $\Var(f(Z))/t$ without imposing additional assumptions \citep{van2002semiparametric}.


\subsection{Multivariate asymptotic confidence sequences}\label{section:multivariate}
We have thus far focused on univariate \asympcs{}s since even this simple setting encompasses several areas of application including some of our main causal inference-related motivations found in \cref{section:csate}. One may nevertheless wonder if the notion of an \asympcs{} can be generalized to $\RR^d$ for $d \geq 2$, and if so, whether constructions thereof are possible. In this section, we provide a definition for multivariate \asympcs{}s and construct explicit examples of them, relying on certain multivariate strong invariance principles due to \citet{einmahl1987strong} and nonasymptotic confidence sets for means of Gaussian vectors due to \citet{manole2023martingale} and \citet{chugg2023time}. In what follows, let $\nu$ be the Lebesgue measure on $\RR^d$.
\begin{definition}[Multivariate asymptotic confidence sequences]\label{definition:multivariate-asympcs}
  Let $\Tcal$ be a totally ordered infinite set including a minimum value $t_0 \in \Tcal$. We say that the sequence of $\RR^d$-valued random sets $(C_t)_{t \in \Tcal}$ is a $(1-\alpha)$-\emph{asymptotic confidence sequence} (AsympCS) for a sequence of parameters $(\theta_t)_{t\in \Tcal}$ taking values in $\RR^d$ if there exists a nonasymptotic $(1-\alpha)$-CS $(C_t^\star)_{t\in \Tcal}$ for $(\theta_t)_{t \in \Tcal}$ meaning that
  \begin{equation}
    \PP \left ( \forall t \in \Tcal ,\ \theta_t \in C_t^\star \right ) \geq 1-\alpha,
  \end{equation}
  so that the normalized measure of the symmetric difference between $(C_t^\star)_{t\in \Tcal}$ and $(C_t)_{t\in \Tcal}$ a.s.~vanishes:
  \begin{equation}
    \frac{\nu(C_t \symdiff C_t^\star) }{\nu(C_t^\star)} \to 0~~\text{almost surely.}
  \end{equation}
\end{definition}

In the univariate case, any \asympcs{} satisfies \cref{definition:multivariate-asympcs} (in fact, \cref{definition:multivariate-asympcs} is still more general since it encompasses sets other than intervals --- such as disjoint unions thereof --- but we ignore this technicality since most confidence sets we are interested in are in fact compact and connected). Now, let us use a strong invariance principle due to \citet{einmahl1987strong} combined with the nonasymptotic confidence sequences for means of sub-Gaussian random vectors due to \citet{manole2023martingale} to derive multivariate \asympcs{}s for means of \iid{} random vectors with certain finite moments.
\begin{proposition}[Multivariate \asympcs{}s for means of \iid{} random vectors]\label{proposition:multivariate-asympcs}
  Let $\infseqt{X_t}$ be \iid{} random vectors with mean $\mu \in \RR^d$ and covariance matrix $\Sigma \in \RR^{d \times d}$ so that
  \begin{equation}\label{eq:multivariate-finite-moments}
    \EE \| X \|_2^{2+\delta} < \infty
  \end{equation}
  for some $\delta > 0$. Let $\widehat \mu_t := \frac{1}{t} \sum_{i=1}^t X_i$ be the sample mean and $\widehat \gamma_t \equiv \gamma(\widehat \Sigma_t)$ the maximum eigenvalue of the empirical covariance matrix $\widehat \Sigma_t := \frac{1}{t} \sum_{i=1}^t (X_i - \widehat \mu_t) (X_i - \widehat \mu_t)^\top$. Then,
  \begin{equation}\label{eq:multivariate-manole}
    C_t := \left \{ \mu' \in \RR^d : \left \| \widehat \mu_t - \mu' \right \|_2 <   \sqrt{\widehat \gamma_t 16 \left [ 2 \log ( \log_2( t) + 1) + \log(1/\alpha) + d \log 5 \right ]/ t} \right \}
  \end{equation}
  forms a multivariate $(1-\alpha)$-\asympcs{} for $\mu$. 
The set $C_t$ can be alternatively written as
\begin{equation}
  \ball(\widehat \mu_t, r_t),\quad\text{where}~~r_t := \sqrt{\widehat \gamma_t  16  \left [ 2 \log ( \log_2( t) + 1) + \log(1/\alpha) + d \log 5 \right ]/t}
\end{equation}
and $\ball(a, r)$ is the ball in $\RR^d$ centered at $a \in \RR^d$ with radius $r > 0$.
\end{proposition}
Similar to the relationship between \cref{proposition:acsLIL} and \cref{theorem:acs}, one could have replaced $\sqrt{16[2 \log (\log_2(t) + 1) + \log(1/\alpha) + d\log 5]}$ with any other time-uniform $(1-\alpha)$ boundary for means of 1-sub-Gaussian random vectors. For example, replacing the bound of \citet{manole2023martingale} with that of \citet{chugg2023time} in Step 2 of the proof of \cref{proposition:multivariate-asympcs}, we have that for any $\rho > 0$,
\begin{equation}\label{eq:multivariate-chugg}
  C_t := \left \{ \mu' \in \RR^d : \| \widehat \mu_t - \mu'\|_2 < \sqrt{\frac{\widehat \gamma_t9d}{2} \cdot \frac{1+t \rho^2 }{t^2 \rho^2} \cdot \left [ 2 + \log \left ( \frac{\sqrt{1+ t \rho^2}}{\alpha} \right )  \right ]} \right \}
\end{equation}
forms a $(1-\alpha)$-\asympcs{} for $\mu$ for the same reason that \cref{theorem:lindeberg-martingale-asympcs} reduces to \cref{theorem:acs} for an appropriate tuning parameter $\rho > 0$.
The bounds in \eqref{eq:multivariate-manole} and \eqref{eq:multivariate-chugg} can be thought of as multivariate analogues of those found in \cref{proposition:acsLIL} and \cref{theorem:acs}, respectively.

\begin{proof}[Proof of Proposition~\ref{proposition:multivariate-asympcs}]
  Like many of the proofs found throughout this paper, multivariate \asympcs{}s are constructed from three key ingredients: strongly approximating partial sums of Gaussian vectors, a nonasymptotic time-uniform concentration result, and a justification for nuisance parameter estimation. We deal with these three ingredients in Steps 1, 2, and 3, respectively, ultimately combining them in Step 4. 
  \paragraph*{Step 1: Multivariate strong invariance via \citet{einmahl1987strong}}
  Using the strong invariance principle of \citet{einmahl1987strong} combined with the condition in \eqref{eq:multivariate-finite-moments}, we have that on a sufficiently rich probability space, there exist \iid{} multivariate Gaussian random vectors $\infseqt{Y_t}$ with mean zero and covariance matrix $\Sigma := \Cov(X)$ so that for some $q > 2$,
  \begin{equation}
    (\widehat \mu_t - \mu) - \widebar Y_t = o(t^{1/q - 1})
  \end{equation}
  coordinate-wise almost-surely where $\widebar Y_t := \frac{1}{t} \sum_{i=1}^t Y_i$.
  \paragraph*{Step 2: Nonasymptotic time-uniform concentration for Gaussian random vectors}
  We will rely on a result due to \citet[Corollary 23]{manole2023martingale} which states that if for some $\gamma > 0$,
  \begin{equation}
    \sup_{u \in \RR^d : \| u \|_2 = 1} \log \EE \exp \left \{ \lambda \left \langle u, Y \right \rangle
      \right \} \leq \gamma \lambda^2/2;\quad \lambda \in \RR,
  \end{equation}
  then with probability at least $(1-\alpha)$,
  \begin{equation}
    \forall t \geq 1,\ \| \widebar Y_t \| < \sqrt{\gamma 16 [2 \log ( \log_2 t + 1) + \log(1/\alpha) + d\log 5] / t}.
  \end{equation}
  Now, notice that the \iid{} random variables $\infseqt{Y_t}$ from Step 1 are multivariate Gaussian and hence their moment generating function is given by
  \begin{equation}
    \EE \exp \left \{ \langle u, Y \rangle \right \} = \exp \left \{ \frac{1}{2} u^T \Sigma u \right \}; \quad u \in \RR^d.
  \end{equation}
  Consequently, we have that for any $\lambda \in \RR$,
  \begin{equation}
    \sup_{u \in \RR^d : \| u \|_2 = 1} \log \EE \exp \left \{ \lambda \left \langle u, Y \right \rangle
      \right \} = \gamma(\Sigma) \cdot \lambda^2 /2
  \end{equation}
  where $\gamma \equiv \gamma(\Sigma)$ is the largest eigenvalue value of $\Sigma$. Applying \citet[Corollary 23]{manole2023martingale} as discussed above, we have
  \begin{equation}
    \PP \left ( \exists t \geq 1 : \left \| \widebar Y_t \right \|_2 \geq  \sqrt{ \gamma(\Sigma) 16 \left [ 2 \log (\log_2 (t) + 1) + \log(1/\alpha) + d\log 5 \right ] / t}  \right ) \leq \alpha.
  \end{equation}
  \paragraph*{Step 3: Strongly consistent nuisance parameter estimation}
  Applying the strong law of large numbers coordinate-wise, we have that the sample covariance matrix $\widehat \Sigma_t$ is strongly consistent for the covariance matrix $\Sigma$ in the Frobenius norm $\| \cdot \|_\mathrm{F}$, meaning that $\|\widehat \Sigma_t - \Sigma\|_\frob = o(1)$ almost surely. By Courant-Fischer,
  \begin{equation}
    |\gamma(\widehat \Sigma_t) - \gamma(\Sigma)| \leq \| \widehat \Sigma_t - \Sigma \|_\frob = o(1),
  \end{equation}
  and hence we have that $\gamma(\widehat \Sigma_t) \to \gamma(\Sigma)$ almost surely.
  \paragraph*{Step 4: Proving that $\ball(\widehat \mu_t, r_t)$ is a $(1-\alpha)$-\asympcs{} for $\mu$} Combining Steps 2 and 3, there exists a radius $r_t^\star \equiv r_t + o \left ( \sqrt{\log \log t / t} \right )$ so that
    $\PP \left ( \exists t \geq 1 : \left \| \widehat \mu_t - \mu \right \|_2 \geq r_t^\star \right ) \leq \alpha$,
    or in other words, we have with probability at least $(1-\alpha)$,
    \begin{equation}
      \forall t\geq 1,\quad \mu \in \ball \left ( \widehat \mu_t , r_t^\star \right ).
    \end{equation}
    Therefore, it remains to show that $\nu \left (\ball (\widehat \mu_t, r_t^\star) \symdiff \ball (\widehat \mu_t, r_t) \right ) / \nu \left ( \ball(\widehat \mu_t , r_t^\star) \right ) \to 0$ where
    \begin{equation}
      r_t :=  \sqrt{\gamma(\widehat \Sigma_t)16 \left [ 2 \log(\log_2(t) + 1) + \log(1/\alpha) + d \log 5 \right ] / t}.
    \end{equation}
    Indeed, writing out the aforementioned normalized symmetric difference, we notice that
  \begin{align}
    &\frac{\nu \left ( \ball \left (\widehat \mu_t, r_t^\star \right ) \symdiff{}\ \ball \left (\widehat \mu_t, r_t \right ) \right )}{\nu \left ( \ball \left (\widehat \mu_t, r_t^\star \right )  \right )}\\
    =\ & \frac{\left \lvert \cancel{\frac{\pi^{d/2}}{\Gamma(d/2 + 1)}} \cdot \left (r_t^\star \right )^d - \cancel{\frac{\pi^{d/2}}{\Gamma(d/2 + 1)}} \cdot r_t^d \right \rvert}{\cancel{\frac{\pi^{d/2}}{\Gamma(d/2 + 1)}} \cdot \left ( r_t^\star \right )^d}\\
    =\ & \left \lvert 1 - \left ( \frac{r_t }{r_t + o(\sqrt{\log \log t / t})} \right )^d \right \rvert \\
    =\ & \left \lvert 1 - \left ( \frac{C}{C + o(1)} \right )^d \right \rvert = o(1)
  \end{align}
  almost surely for some constant $C > 0$ by the continuous mapping theorem and since $r_t \asymp \sqrt{\log \log t / t}$. This completes the proof of \cref{proposition:multivariate-asympcs}.
\end{proof}


\end{document}